\documentclass{article}
\usepackage[english]{babel}
\usepackage{amsmath,amssymb,latexsym,theorem}
\usepackage{authblk}
\usepackage{color}
\usepackage{epsfig}
\numberwithin{equation}{section}

\newcommand{\assign}{:=}
\newcommand{\nin}{\not\in}

\newcommand{\tmtextit}[1]{\text{{\itshape{#1}}}}
\newcommand{\tmop}[1]{\ensuremath{\operatorname{#1}}}
\newenvironment{itemizedot}{\begin{itemize} }{\end{itemize}}
\newenvironment{itemizeminus}{\begin{itemize} }{\end{itemize}}
\newenvironment{proof}{\noindent\textbf{Proof\ }}{\hspace*{\fill}$\Box$\medskip}
\newtheorem{theorem}{Theorem}[section]

\newtheorem{lemma}[theorem]{Lemma}
\newtheorem{proposition}[theorem]{Proposition}
\newtheorem{example}[theorem]{Example}
\newtheorem{remark}[theorem]{Remark}

\begin{document}

\title{On the Number of Normalized Ground State Solutions for a class of Elliptic Equations with general  nonlinearities and potentials}
\author{Hichem Hajaiej}
\affil{Department of Mathematics, California State University at Los Angeles, Los Angeles, CA 90032, USA}
\author{Eliot Pacherie}
\affil{NYUAD Research Institute, New York University Abu Dhabi \\ PO Box 129188, Abu Dhabi, UAE}
\author{Linjie Song}
\affil{Department of Mathematical Sciences, Tsinghua University, Beijing 100084, China}

\maketitle

\begin{abstract}
We provide a precise description of the set of normalized ground state solutions (NGSS) for the class of elliptic equations:
   $$
   -\Delta u - \lambda u + V (| x |) u - f (| x |, u) = 0,\quad\text{in}\quad \mathbb{R}^n,\ n\geq 1.
   $$
    In particular, we show that under suitable assumptions on $V$ and $f$, the NGSS is unique for all the masses except for at most a finite number. Moreover, we prove that when unique, the NGSS $u_c$ is a smooth function of the mass $c.$ Our method is as follow: using the NGSS for a given mass $c$, we construct an exhaustive list of potential candidates to the minimization problem for masses close to $c$, and we develop a strategy how to pick the right one. In particular, if there is a unique NGSS for a given mass $c_0,$ then this uniqueness property is inherited for all the masses $c$ close to $c_0.$ Our method is general and applies to other equations provided that some key properties hold true. 

\end{abstract}
 
\section{Introduction and presentation of the results} \label{sec. intro}

\subsection{The minimization problem and known results}

In this paper, we provide suitable assumptions on the nonlinearities $f$ and potentials $V$ for which NGSS of Equation \eqref{solarstone} below is unique for all masses $c,$ except at most for a finite number. Our study covers the classical nonlinearities as well as the non-autonomous ones. More precisely, we will focus our attention on the equation:

\begin{equation}
  \tmop{Eq}_{\lambda} (u) = - \Delta u + V(|x|) u - \lambda u - f (|x|, u) = 0,
  \label{solarstone}
\end{equation}
on $\mathbb{R}^n$ with $\lambda < \Lambda_0 \leqslant 0$ where $\Lambda_0 := \inf \sigma(-\Delta + V)$ is the infimum of the spectrum of the operator $-\Delta + V$ on $L^2(\mathbb{R}^n)$. The function $V$ is a potential and $f$ is the nonlinear part of the equation. Throughought the whole paper, we make the following basic assumptions on $V$ and $f$:
\begin{itemizeminus}
  \item $(F 1)$ The function $f (r, t)$ is $C^1$ with respect to $t \in
  \mathbb{R}$ and $C^1$ with respect to $r \in \mathbb{R}^+$, $f (r, 0) =
  f_t (r, 0) = 0$, and $f (r, - t) = - f (r, t)$.

  \item $(F 2)$ $V$ is $C^1$ and $V (| x |) \rightarrow 0$ when $ |x| \rightarrow + \infty $.
\end{itemizeminus}

We say that \eqref{solarstone} is autonomous if $V = 0$ and $f$ is independent of
its first variable.

\

We define the energy
\[ E (u) \assign  \int_{\mathbb{R}^n} \frac{1}{2}| \nabla u |^2 + \frac{1}{2}V (| x |) |
   u |^2 - F (| x |, u) \]
where $F_t (r, t) = f (r, t)$ and the mass
\[ Q (u) \assign \frac{1}{2} \int_{\mathbb{R}^n} | u |^2 . \]
Our goal is to study the set of minimizers of the energy at a fixed mass, that
is, for $c$ in some interval $]c_{\ast},c^{\ast}[$ and
\[ m (c) \assign \inf_{u \in X_c} E (u), \quad X_c \assign \{u \in H^1, Q(u) = c\}.\]
We want to explore the properties of the sets
\[ S_c \assign \{ u \in X_c, E (u) = m (c) \}, \]
and
\[ \mathcal{T}_c \assign \{\lambda, \exists u \in S_c \ s.t. \ \tmop{Eq}_{\lambda} (u) = 0\} .\]
Functions in $S_c$ are called \textbf{normalized ground state solutions} (NGSS).
The interest of the study of NGSS has drastically increased after the establishment of the relationship between the sign of the derivative of $Q$ (w.r.t. $\lambda$) and the orbital stability of standing waves of the Cauchy problem from which \eqref{solarstone} derives. In fact, in \cite{GSS}, Grillakis-Shatah-Strauss showed that $u_\lambda e^{-i\lambda t}$ is stable when $\partial_\lambda Q(u_\lambda)>0$ while if $\partial_\lambda Q(u_\lambda)<0,$  $u_\lambda e^{-i\lambda t}$ is unstable ($u_\lambda$ is a solution to \eqref{solarstone}). The study of the uniqueness of NGSS is salient as this property, together with the uniqueness of the branch of solutions of the PDE (\eqref{solarstone} in our setting) provide two crucial information: the stability of the standing waves, and the characterization of all the solutions of the PDE as NGSS with distinct masses.


Since the breakthrough result of \cite{GSS}, the study of NGSS  has gained a lot of interest as it provides physicists and engineers with precious guidance on the realization of the model they want to study. Many papers addressed the existence, non-existence, or multiplicity. However, the literature remained quite silent regarding the uniqueness of NGSS. So far, only a handful of papers dealt with this issue. Contrarily to the study of the uniqueness of positive solutions of \eqref{solarstone}, where many valuable contributions addressing various and numerous nonlinearities, the uniqueness of NGSS is a more challenging and delicate problem where novel ideas and techniques are needed to make significant progress. A tangible proof of the delicacy of the study of the uniqueness of NGSS is provided by the counter-example exhibited in \cite[Appendix A]{Hajaiej-Song-unique} where the authors showed the uniqueness of NGSS for all masses $c>0$ except for one mass $c_1$ for a pure type power nonlinearity. Moreover, they proved that for $c_1,$ there are exactly two NGSS. This counter-example can be extended to prove the non-uniqueness for a finite number of masses by modifying the nonlinearity. The counter-example in \cite{Hajaiej-Song-unique} was constructed by a nonliearity $f$ which satisfies $f(x,-t)\neq-f(x,t)$ (and thus not $(F 1)$).
To the best of our knowledge, the only situations (except \cite{Hajaiej-Song-unique}) where the uniqueness of NGSS has been established is when the mass $Q$ (or $m(c)$) is explicitly given as a function of $c$ and the unique positive radial decreasing solution to \eqref{solarstone}. These cases only cover situations where $f$ is a pure power nonlinearity and the potential is null (see identity (1.8) in Theorem 1.2 of \cite{HZ} for example). In \cite{Lewin-Nodari-2020}, the authors addressed  the double power nonlinearity. This certainly was a more delicate study, but even in this case, explicit expressions of $Q$ can be obtained for small and large masses $c.$ (see (3.10) and (3.18) of \cite{Lewin-Nodari-2020} respectively). In their paper, some partial uniqueness results have been proven in their context ($f(|x|,u)=u^p-u^q,$ $p>q$, see Theorem \ref{mainthing0} for more details). Other results were conjectured (conjecture 2,3 in \cite{Lewin-Nodari-2020}) and investigated numerically. Their analysis of the mass $Q$  (as a function of $c$) is an interesting first step towards a better understanding of general nonlinearities. However, a self-contained approach is clearly needed to discuss the uniqueness  of NGSS for \eqref{solarstone} with $f$ and $V$ general. In such situations, there is no hope to obtain the expression of the mass as a function of the unique positive solution of \eqref{solarstone}.
More recently, Hajaiej and Song \cite{Hajaiej-Song-unique} 
have developed two general methods to show the uniqueness of NGSS. Their Method 1 states that for a given $c>0,$ $|S_c|=1,$ is equivalent to the differentiability of $m$ at $c$. It is easy to show that $m$ is strictly decreasing for large classes of $f$ and $V$, therefore $|S_c|=1$ for almost all $c>0,$ in these cases. This method was based on the fact that the branch of solutions to \eqref{solarstone} is unique and that NGSS are all of Morse index 1. In this work, we are in a much more general setting where the number of branches can be infinite. We also provide more qualitative and quantitative properties on the possible masses where NGSS are possibly not unique. 
Let us point out that Method 1 of \cite{Hajaiej-Song-unique} covers the non-local operators and that only Method 2 of \cite{Hajaiej-Song-unique} addresses the critical and super-critical nonlinearities as it does not  use the finiteness of the energy functional. Therefore, Method 1 and Method 2 of \cite{Hajaiej-Song-unique} as well as the new approach presented here can be seen as complementary.\\

 In this work, we explore more properties on the NGSS such as their possible number for a given mass, the continuity and differentiability of $u_c$ (the NGSS when it is unique), $m(c)$ and $\lambda(c)$ with respect to the mass $c.$ These are highly important features, especially that Method 1 of \cite{Hajaiej-Song-unique} established the equivalence between the uniqueness of NGSS, and the differentiability of $m$ at a given $c$. Some of the differentiability properties of $m$ and $\lambda$ were touched upon very superficially in a few previous works. We will discuss this in details in the next paragraph and explain the importance of our contribution.\\

In this paper, we always assume that this minimization problem satisfies some classical properties (see properties $(A)$
), and we will deal with two cases. The first one is when $f$ is real analytic w.r.t. its second variable. We also need the assumption $(N)$ to ensure that any ground state solution (GSS, defined as the minimizer of $E-\lambda Q$ constrained on the Nehari manifold) is indeed a NGSS. A full study of the relationship between GSS and NGSS can be found in \cite{Hajaiej-Song-unique}. If $\lambda \in \mathcal{T}_c$, then all ground states for this $\lambda$ also belong to $S_c$ (see \cite[Appendix C]{Hajaiej-Song-unique} for a very general equation by only assuming $(N)$). Then, we will show that the set $\mathcal{T}_c$ contains
exactly one element for all $c$ in the existence range except at most for a finite number of values, and the set $S_c$ contains finitely many elements for all $c$. We will say that $c$ is a {\bf "bad" mass} if $\mathcal{T}_c$ contains more than one element. Between two "bad" masses, defining $\lambda(c)$ as the unique element of $\mathcal{T}_c$, $\lambda(c)$ is continuous and strictly decreasing, which implies that the elements in $\mathcal{T}_c$ form an open interval. The union of $S_c$ for $c$ between two bad masses consists of many finite  $C^\alpha$ branches for some $\alpha \in (0,1)$ except at most for a finite number of isolated elements. Additionally except for finite values, these branches are $C^1$. This improves the regularity of $\lambda(c)$. Indeed, $\lambda(c)$ is $C^\alpha$ between two "bad" masses and $C^1$ except for a finite number of values. If we replace $(N)$ by a weak uniqueness hypothesis $(C)$, all the results above hold and the number of branches is exactly one. (See more details in Theorems \ref{mainthing0} and \ref{mainthing2}). Let us point out that we don't require the real analyticity of $f$ in the second case. As an alternative, we suppose conditions $(C)$ and $(D)$, and give similar characterizations of $\mathcal{T}_c$ and $S_c$. Particularly, the union of $S_c$ for $c$ between two "bad" masses consists exactly of one $C^1$ branch. (See more details in Theorems \ref{mainthing0} and \ref{mainthing1}).\\

For the sake of clarity, the paper will be divided in two parts. In the first one (Sections \ref{sec. intro} and \ref{behind}), we show our theorems and prove them using our key properties such as $(A)$, $(B)$, $(C)$, $(N)$ that are written below. In the next part (Sections \ref{sec. conditions} and \ref{sono}) we will give some classes  of functions $f$ and $V$ for which properties $(A)$, $(B)$, $(C)$ and $(N)$ are satisfied. After that, we provide concrete hypotheses on $f$ and $V$, and exhibit concrete examples such that all the properties in the first part are satisfied. As a foretaste, examples include the following important classes of functions. Precise assumptions will be provided in Section 4.
\begin{itemizeminus}
	\item $V\equiv 0, f(x,u) = \Sigma_{i=1}^Na_i(x)|u|^{p_i-2}u$, $a_i(x) > 0$, $p_i \in ]2,2+\frac{4}{n}], p_i \in \mathbb{N}$.
	\item $V\equiv 0, f(u) = |u|^{p-2}u - |u|^{q-2}u, p,q \in \mathbb{N}$.
	\item $\Lambda_0$ is the first eigenvalue of $-\Delta + V$ and $f = |u|^{p-2}u$, $p \in ]2,2+\frac{4}{n}]$ and more conditions on $V$ should be given to satisfy some hypotheses.
\end{itemizeminus}
We believe that our method should apply to sum of pure powers whose exponents are not necessarily integers. We will discuss this later on.
\\

For $u$ solution of $\tmop{Eq}_{\lambda} (u) = 0$, we define the linearized
operator around $u$ as
\[ L_u (\varphi) \assign - \Delta \varphi + V \varphi - \lambda \varphi - f_t
   (x, u_{\lambda}) (\varphi) . \]
Finally, we define the scalar product between two real functions as
\[ \langle f, g \rangle \assign \int_{\mathbb{R}^n} f g. \]

\subsubsection{Properties of the minimization problem}

We denote by $(A)$ the set of the following properties on a set of mass of the form $]c_{\ast},c^{\ast}[$:
\begin{itemizeminus}
	\item $(A 1)$ For any $c \in ]c_{\ast},c^{\ast}[$, $m(c) > -\infty$ is continuous and the set $S_c$
	contains at least one element. Furthermore, every element $u \in S_c$ is
	positive, radial  decreasing, and satisfies $\tmop{Eq}_{\lambda} (u) = 0$ for some $\lambda <
	\Lambda_0 \leq 0$. Moreover, $L_{u}$ has exactly one simple negative eigenvalue.
	
	\item $(A 2)$ For any $c \in ]c_{\ast},c^{\ast}[$ and $u \in S_c$ with $\tmop{Eq}_{\lambda} (u)
	= 0$, we have $\tmop{Ker} L_u \cap L^2_{\tmop{rad}} = \{ 0 \}$.  
	
	\item $(A 3)$ For $c \in ]c_{\ast},c^{\ast}[,$ any sequence $(v_k)_{k \in \mathbb{N}} \in (H^1
	(\mathbb{R}^n))^{\mathbb{N}}$ with $Q (v_k) \rightarrow c, E (v_k)
	\rightarrow m (c)$, up to a subsequence (and translation in autonomous
	cases), $v_k$ converges to $v$ strongly in $H^1 (\mathbb{R}^n)$ to a limit
	$v \in S_c$.
\end{itemizeminus}
We also denote by $(B)$,$(C)$ and $(N)$ the following properties (on a similar set):
\begin{itemizeminus}
	\item $(B)$ There exist $c_0, c^0 \in ]c_{\ast},c^{\ast}[$ such that for any $c \in ]c_{\ast},c^{\ast}[ \backslash \left[ c_0, c^0 \right]$, $S_c$ contains at most one element.
	\item $(C)$ For any $c \in ]c_{\ast},c^{\ast}[$ and $u, v \in S_c$ with $u \neq v$, if
	$\tmop{Eq}_{\lambda_1} (u) = 0 = \tmop{Eq}_{\lambda_2} (v)$, then $\lambda_1
	\neq \lambda_2$.
 \item $(N)$ Define $\Phi_\lambda = E - \lambda Q$, $
		\mathcal{N}_{\lambda} = \{u \in H^1 \setminus \{0\}: D_u\Phi_\lambda(u)(u) = 0\}
		$. For any $u \in H^1 \setminus \{0\}$, there exists a unique $t(u) = t(\lambda,u) > 0$ such that $t(u)u \in \mathcal{N}_{\lambda}$, and $\Phi_\lambda(t(u)u) = \max_{t > 0}\Phi_\lambda(tu)$.
\end{itemizeminus}
Finally, define $\lambda_{min}(c) := \min\mathcal{T}_c$ and $\lambda_{max}(c) := \max\mathcal{T}_c$, then $(D)$ is the set of the following two properties:
	\begin{itemizeminus}
		\item $(D 1)$ For any $\epsilon_k \to 0^+$ and $u_{c+\epsilon_k} \in S_{c+\epsilon_k}$ with $\tmop{Eq}_{\lambda_{min}(c+\epsilon_k)}u_{c+\epsilon_k} = 0$, up to a subsequence, there exists $u_c \in S_c$ such that $\tmop{Eq}_{\lambda_{min}(c)}u_c = 0$, $u_{c+\epsilon_k} \to u_c$ strongly in $H^1$ and
		\begin{align} \label{eq right assum}
			\limsup_{k \to \infty}\|\frac{u_{c+\epsilon_k} - u_c}{\epsilon_k}\|_{H^1} < \infty.
		\end{align}
	\item $(D 2)$ For any $\epsilon_{k} \to 0^+$ and $u_{c-\epsilon_k} \in S_{c-\epsilon_k}$ with $\tmop{Eq}_{\lambda_{max}(c-\epsilon_k)}u_{c-\epsilon_k} = 0$, up to a subsequence, there exists $u_c \in S_c$ such that $\tmop{Eq}_{\lambda_{max}}(c)u_c = 0$, $u_{c-\epsilon_k} \to u_c$ strongly in $H^1$ and
	\begin{align} \label{eq left assum}
		\limsup_{k \to \infty}\|\frac{u_{c-\epsilon_k} - u_c}{\epsilon_k}\|_{H^1} < \infty.
	\end{align}
	\end{itemizeminus}

\begin{remark}
	We give some remarks on these properties.
	\begin{itemizeminus}
		\item $(A 2)$ is a non-degeneracy result. Together with $(A 3)$, they imply that for any fixed $\lambda$, $\tmop{Eq}_{\lambda} (u) = 0$ possesses only finitely many NGSS in $S_c$. In fact, by \cite[Appendix C]{Hajaiej-Song-unique}, if $u \in S_{c_1}, v \in S_{c_2}, \tmop{Eq}_{\lambda} (u) = \tmop{Eq}_{\lambda} (v) = 0$, then $c_1 = c_2$. Then arguing by contradiction we can deduce the claim above. If one tries to remove $(A 2)$, we may end up with infinitely many NGSS in many situations.
		\item $(C)$ is a weak uniqueness result. Indeed, if $u \in S_c$ and $\tmop{Eq}_{\lambda} (u) = 0$, then \cite[Appendix C]{Hajaiej-Song-unique} yields that any GSS at $\lambda$ has mass $c$ and belongs to $S_c$, thus, $(C)$ means that the GSS at $\lambda$ is unique.
        \item $(N)$ is used to show that any NGSS is indeed a GSS, and other GSS at the same $\lambda$ are NGSS (see \cite[Appendix C]{Hajaiej-Song-unique}).

	\end{itemizeminus}
\end{remark}

\subsection{Main results}


We recall the definitions of
\[ \tmop{Eq}_{\lambda} (u) = - \Delta u + V (| x |) u - \lambda u - f (| x |,
   u), \]
\[ Q (u) = \frac{1}{2} \int_{\mathbb{R}^n} u^2, \]
\[ m (c) = \{ \inf E (u), u \in H^1 (\mathbb{R}^n), Q (u) = c \}, \]
\[ S_c = \{ u \in H^1 (\mathbb{R}^n), E (u) = m (c), Q (u) = c \}, \]
\[ \mathcal{T}_c = \{ \lambda, \exists u \in S_c, \tmop{Eq}_{\lambda} (u) = 0
   \}, \]
\[ \mathcal{P}= \{  c \in] c_{\ast}, c^{\ast} [, | \mathcal{T}_c |
   \geqslant 2  \} \]
and we define
\[ \mathcal{M} \assign \{ c \in] c_{\ast}, c^{\ast} [, | S_c | \geqslant 2 \}
   . \]
We always have $\mathcal{P} \subset \mathcal{M}$ and if we assume $(C)$, then
$\mathcal{P}=\mathcal{M}$. 
We now state the main results of this paper.

\subsubsection{Finitude of $S_c$}

\begin{theorem} \label{fleurs}
  Assume that $(A)$ and one of the following properties holds for some $c \in]
  c_{\ast}, c^{\ast} [$:
  \begin{itemizedot}
    \item $f (r, u)$ is analytic w.r.t $u$
    
    \item $\langle L^{- 1}_u (u), u \rangle \neq 0$ for all $u \in S_c$,
  \end{itemizedot}
  then $S_c$ is a finite set.
\end{theorem}


\subsubsection{Accumulation points of $\mathcal{P}$}

\begin{theorem} \label{mainthing0}\label{finite of P}
  Assume that $(A)$ and one of the following properties:
  \begin{itemizedot}
    \item $f (r, u)$ is analytic w.r.t $u$ 
    and $(N)$ holds
    
    \item $(C)$ and $(D)$ hold.
  \end{itemizedot}
  Then, any accumulation point of $\mathcal{P}$ must be either $c_{\ast}$ or
  $c^{\ast}$. In particular, if we also assume $(B)$, then $\mathcal{P}$ is a
  finite set. If $(C)$ holds, then $\mathcal{P}=\mathcal{M}$ and thus
  $\mathcal{M}$ is also a finite set.
\end{theorem}

We recall that accumulations points of a set $X \subset ] c_{\ast},
c^{\ast} [$ is the set of points $c \in [c_{\ast}, c^{\ast}]$ such
that for any $\varepsilon > 0$, the set $X \cap ([c - \varepsilon, c +
\varepsilon] \backslash c)$ is non empty. In particular, isolated points of
$X$ are not accumulation points, and $X$ has no accumulation points if and
only if it is a finite set. This result state that the set $\mathcal{P}$ of
``bad'' masses (that is $| \mathcal{T}_c | > 1$) is finite on any compact of
$] c_{\ast}, c^{\ast} [$. We can say more about
$\mathcal{T}_c$ between two of these ``bad'' masses, see the next result.

\subsubsection{Properties of $\mathcal{T}_c$}

\begin{theorem} \label{Planitia}
  Under the assumptions of Theorem \ref{finite of P}, for any $c_0, c^0 \in ]
  c_{\ast}, c^{\ast} [$ with $c_0 < c^0$, we define by $c_1 < \ldots
  < c_k$ the elements of $\mathcal{P} \cap ] c_0, c^0 [$
  (this can be empty, then $k = 0$) and $c_{k + 1} = c^0$. If we assume $(B)$,
  we can take $c_0 = c_{\ast}$ and $c^0 = c^{\ast}$.
  
  Then, taking $\lambda (c)$ the unique element of $\mathcal{T}_c$ for $c \in]
  c_j, c_{j + 1} [, j \in \{ 0 \ldots k \}$, the function $c \rightarrow
  \lambda (c)$ is continuous and strictly decreasing on $] c_j, c_{j + 1} [$
  and $\max \mathcal{T}_{c_{j + 1}} < \min \mathcal{T}_{c_j}$. Furthermore, we
  have
  \[ \{ \lambda \in \mathcal{T}_c : c \in] c_j, c_{j + 1} [\} =] \max
     \mathcal{T}_{c_{j + 1}}, \min \mathcal{T}_{c_j} [. \]
  Finally, for any $\tilde{c} \in] c_j, c_{j + 1} [$ we have the following
  equivalence:
  \[ \tilde{c} \nin \mathcal{P} \Leftrightarrow \tmop{the} \tmop{two}
     \tmop{limits} \lim_{c \rightarrow \tilde{c}^{\pm}} \frac{m (c) - m
     (\tilde{c})}{c - \tilde{c}} \tmop{exist} \tmop{and} \tmop{are}
     \tmop{equal} . \]
\end{theorem}

We will show in Proposition \ref{properties of m(c)} that $\min \mathcal{T}_c$ and
$\max \mathcal{T}_c$ are well defined and not $- \infty$ or $0$ for any $c
\in] c_{\ast}, c^{\ast} [$. This result state that the Lagrange multipliers
$\lambda$ of elements in the sets $S_c$ are regrouped in continuous branches.
$c \rightarrow \lambda (c) \in \mathcal{T}_c$ is single valued, continuous
strictly decreasing between two ``bad'' masses, and is multivalued and
discontinuous exactly at these finitely many ``bad'' masses. Furthermore, the
``bad'' masses are caracterise by points of discontinuity of the derivative of
$c \rightarrow m (c)$, the minimal energy at mass $c$.

\subsubsection{Properties of $S_c$}

Here, $k, c_j$ and $\lambda (c)$ are defined as in Theorem \ref{Planitia}. We state
more properties of minimizers between ``bad'' masses.

\begin{theorem}[Properties of $S_c$, I] \label{mainthing1}
  Suppose that properties $(A)$, $(C)$ and $(D)$ hold, $|S_c| = 1$ for any $c \in ] c_0, c^0[ \setminus \{c_1, \cdots, c_{k}\}$ and the set
  \[ \mathcal{M} = \{ c \in ]c_{\ast},c^{\ast}[, | S_c | \geqslant 2 \} \]
  is such that $\mathcal{M} \cap ] c_0, c^0[ = \{c_1,\cdots,c_k\}$.

  Furthermore, for any $j \in [0 \ldots k]$ and $c \in] c_j, c_{j + 1} [$, the unique element of
  the set $S_c$, denoted $u_c$, satisfies
  \[ c \rightarrow u_c \in C^1 (] c_j, c_{j + 1} [, H^1 (\mathbb{R}^n)) . \]

  Moreover, for any $j \in [0 \ldots k]$, $\lambda(c)$ is $C^1$, $m(c)$ is $C^2$ and $m''(c) = \lambda'(c) < 0$ on $] c_j, c_{j + 1} [$.
\end{theorem}

\begin{theorem}[Properties of $S_c$, II]
 \label{mainthing2}
	Suppose that properties $(A)$, $(N)$ hold and that $f(x,u)$ is real analytic w.r.t. $u$. 	Then, for any $j \in [0 \ldots k]$ , there exist finite ($K_j \geq 1$) branches $c \to u_{c,i} \in C^{\alpha} (] c_j, c_{j + 1} [, H^1 (\mathbb{R}^n)), i \in [1 \ldots K_j]$ for some $\alpha \in (0,1)$, such that except at most for a finite number of isolate elements $\{u \in S_c, c \in ] c_j, c_{j + 1} [\}$ consists of these $K_j$ branches. Furthermore, $c \to u_{c,i}$ is $C^1$ except at most for a finite number of values (that will be called bad values). The existence of a bad value $\tilde{c}$ is equivalent to the existence of $u \in S_{\tilde{c}}$ such that $\langle L^{-1}_u(u),u \rangle = 0$. Also, at a bad value $\tilde{c}$,
	$$\lim_{c \to \tilde{c}}\|\frac{u_{c}-u_{\tilde{c}}}{c-\tilde{c}}\|_{H^1} = \infty.$$
	Moreover, for any $j \in [0 \ldots l]$, $\lambda(c)$ is $C^1$, $m(c)$ is $C^2$ and $m''(c) = \lambda'(c) < 0$ on $] \tilde{c}_j, \tilde{c}_{j + 1} [$ except at these finite bad values.
	
	Let us further assume that the property $(C)$ holds in the rest of this theorem. The set
	\[ \mathcal{M} = \{ c \in ]c_{\ast},c^{\ast}[, | S_c | \geqslant 2 \} \]
	is such that $\mathcal{M} \cap ] c_0, c^0[ = \{c_1,\cdots,c_k\}$. 
	For any $j \in [0 \ldots k]$ and $c \in] c_j, c_{j + 1} [$, the unique element of
	the set $S_c$, denoted $u_c$, satisfies that for some $\alpha \in (0,1)$,
	\[ c \rightarrow u_c \in C^{\alpha} (] c_j, c_{j + 1} [, H^1 (\mathbb{R}^n)) . \]
	Furthermore, $c \to u_{c}$ is $C^1$ except at most for a finite number of values.
\end{theorem}

 For the map $\lambda \mapsto c(\lambda), \lambda \in \mathcal{T}$, note that $c'(\lambda) = \langle L^{-1}_u(u),u \rangle$. Hence, $\langle L^{-1}_u(u),u \rangle = 0$ is equivalent to the fact that $\lambda'(c)$ is infinite.



\




\subsubsection{Remarks on these results}

As previously mentioned, articles on the uniqueness of NGSS are very rare when the nonlinearity is inhomogeneous. This also seems to be the case for the regularity of $m(c), \lambda(c)$ and the map $c \to u_c$. 

Concerning the uniqueness, the most general study was performed in \cite{Hajaiej-Song-unique} so far. Our resluts contain two main improvements compared to the ones in \cite{Hajaiej-Song-unique} where it is shown that $\mathcal{M}$ has Lebesgues measure $0$. In this paper, we prove that $\mathcal{M}$ is a finite set when restricted to any compact set, and is a finite set if $(B)$ holds. Furthermore, outside of $\mathcal{M}$, it was only known previously that $c \rightarrow m(c)$ was $C^1$, we show a much stronger property that $c \rightarrow u_c$ is $C^1$, where $u_c$ is the only element of $S_c$ when the uniqueness result holds.

In \cite{Lewin-Nodari-2020}, the authors addressed a specific case $V = 0, f = u^p- u^q$, where equation
(\ref{solarstone}) admits a unique radial positive solution for any $\lambda \in (-\mu_\ast,0)$ (see \cite{Lewin-Nodari-2020} for the explicit expression of $\mu_\ast$) and this solution is non-degenerate. They proved that $\mathcal{M}$ is finite, i.e. $|S_c| = 1$ except possibly at finitely many points in the existence range, and at those bad values, $S_c$ is also finite. Their proofs relied on the uniqueness, the non-degeneracy of the positive solution, and the analyticity of $Q$ w.r.t. $\lambda$. Our results contain three main improvements. Firstly, we do not need the analyticity of $Q$ in all cases. In fact, we were able to prove that $\mathcal{M}$ is finite with a weak uniqueness condition $(C)$ and without using that $Q$ is real analytic w.r.t. $\lambda$. Moreover, when we assume the analyticity of $f$ with respect to its second variable, we do not need their hypotheses on the uniqueness and the non-degeneracy of positive solutions.  Additionally, when the uniqueness result holds this can be seen as a special situation of our general setting. Finally, our method not only permits to find the derivative of $m(c)$, but it also enables us to find its second derivative, as well as the derivative of $\lambda(c)$. They only computed the left and right derivatives of $m(c).$ Finally, they state that even if $p,q$ are not integer, then $c \rightarrow Q(u_c)$ is still an analytic function of $c$. In all our results, we could replace the hypothese that $f$ is analytic with the fact that $c \rightarrow Q(u_c)$ is analytic. However, the analycity of $Q(u_c)$ was given without proof in \cite{Lewin-Nodari-2020} and we were not able to show this fact.

Regarding the regularity of $\lambda(c),$ in \cite{Stefanov-2019} the author considered  pure power functions. He showed that $\lambda(c)$ is continuous, $m(c)$ is $C^1$, strictly concave down and thus twice differentiable almost everywhere on an interval $(a,b)$ by assuming that for each $c \in (a, b)$, $u_c \in S_c$, then $\lim_{\epsilon \to 0}\|u_{c+\epsilon}-u_c\|_{L^2} = 0$. Such an assumption was not discussed at all in \cite{Stefanov-2019}, it does not seem easy to find concrete conditions such that his assumption $(1.12)$ is satisfied. We show similar results under general hypotheses of $V$ and $f$. We also make another major contribution: Assuming that $\lambda'(c)$ exists, in \cite{Stefanov-2019}, the author proved that $c \to u_c$ is differentiable as an $L^2$-valued mapping. For the power nonlinearities that are real analytic, we show that $c \to u_c$ is differentiable as an $H^1$-valued mapping without assuming that $\lambda'(c)$ exists. Even for nonlinearities that are not real analytic, we obtain the result using the general condition $(D)$ that is weaker than the existence of $\lambda'(c)$.  

We thus have made essential improvements in several areas relative to the results in \cite{Hajaiej-Song-unique,Lewin-Nodari-2020,Stefanov-2019}. Moreover, we provide a unified and general approach that turns out to be very effective to deal with various aspects related to the understanding of the uniqueness of NGSS.


We now give some additional remarks about our main results.

\begin{remark}
	In theorem \ref{mainthing2}, if we replace $(N)$ by $(C)$, we also have $\mathcal{M} \cap ] c_0, c^0[$ is finite for any $c_0 > 0$. Furthermore, for any $j \in [0 \ldots k]$ and $c \in] c_j, c_{j + 1} [$, the unique element of
	the set $S_c$, denoted by $u_c$, satisfies that for some $\alpha \in (0,1)$,
	\[ c \rightarrow u_c \in C^{\alpha} (] c_j, c_{j + 1} [, H^1 (\mathbb{R}^n)) . \]
	In fact, condition $(N)$ will be used in Proposition \ref{ext 2}. If we replace it by $(C)$, as in Proposition \ref{ext 1} we can obtain similar results there.
\end{remark}

\begin{remark}
	It is unclear how to show properties $(D)$ under concrete hypotheses of $f$ and $V$. However, it is worth noticing that in the most relevant physical applications the nonlinearity is real analytic w.r.t. its second variable and we don't use properties $(D)$. Furthermore, $(D)$ can be deduced from the existence of $\lambda'(c)$, which holds almost everywhere. In fact, we have
	\begin{align}
	\limsup_{k \to \infty}\|\frac{u_{c+\epsilon_k} - u_c}{\epsilon_k}\|_{H^1} & = \limsup_{k \to \infty}\left( \|\frac{u_{\lambda(c+\epsilon_k)} - u_{\lambda(c)}}{\lambda(c+\epsilon_k)-\lambda(c)}\|_{H^1}|\frac{\lambda(c+\epsilon_k)-\lambda(c)}{\epsilon_k}|\right) \nonumber \\
	& = \lambda'(c)(\partial_\lambda u_{\lambda})_{|\lambda = \lambda(c)} < \infty.	
	\end{align}
	Finally, we can replace properties $(D)$ by assuming that $\langle L^{-1}_u(u),u \rangle \neq 0$ for any $u \in S_c$.
\end{remark}

\begin{remark}
	Under the setting of this paper, it is hopeful to show that, if $\langle L^{-1}_{u_\lambda}(u_\lambda),u_\lambda \rangle \neq 0$ for any $u_\lambda$ which is non-degenerate in $H^1_{rad}$, then the NGSS is unique. Here we don't need the uniqueness of GSS. This result has been proved in \cite{Hajaiej-Song-unique} under a setting where GSS is unique for any $\lambda < \Lambda_0$. Moreover, \cite{Hajaiej-Song-unique} uses the condition $\partial_\lambda u_\lambda(0) < 0$ to ensure that the assumption $\langle L^{-1}_{u_\lambda}(u_\lambda),u_\lambda \rangle \neq 0$ holds true.
\end{remark}

\begin{remark}
	Assuming $(B)$, and even without the hypotheses $(C)$, the branches constructed in Theorem \ref{mainthing2} are unique on $(c_0,c_1)$ and $(c_k,c^0)$.
\end{remark}

\begin{remark}
	In Theorem \ref{mainthing2} without the property $(C)$, the case $K_j >1$ is possible but is very unlikely to hold. $K_j=1$ means that two branches contain solutions to \eqref{solarstone} with the same energy and mass, at the same $\lambda$, without being equal. With $f$ analytic this seems quite difficult to have. For this reason, we can expect $K_j =1$ except in very degenerate cases.
\end{remark}




\subsection{Extension of non-degenerate solutions and
applications}\label{sand}

We summarize here the main steps of the proof of Theorems \ref{fleurs}-\ref{mainthing2}. The
results of this subsection are proven in section \ref{behind}.

\subsubsection{Two preliminary lemmas}

We start with a first result that does not require more hypotheses than $(F 1)$ and $(F 2)$. It is the fact that radial solutions of the equation decay exponentially fast in position.
\begin{lemma}
	\label{bb}Consider $u_{\lambda}$ a radial solution of $\tmop{Eq}_{\lambda}
	(u_{\lambda}) = 0$ for $\lambda < 0$ with finite energy and $u'(r) < 0$. Then, $u_{\lambda}
	\in C^2 (\mathbb{R}^n, \mathbb{R})$ and there exists $K, a > 0$ such that
	\[ (| u_{\lambda} | + | \nabla u_{\lambda} | + | \nabla^2 u_{\lambda} |) (x)
	\leqslant K e^{- a | x |} . \]
	In particular, there exists $K > 0$ such that for any $\varphi \in
	L^{\infty} (\mathbb{R}^n)$,
	\[ | f_t (x, u_{\lambda}) (\varphi) | (x) \leqslant K e^{- a | x |} |
	\varphi | (x) . \]
\end{lemma}

We then state a coercivity result on minimizers if $(A)$ holds.

\begin{lemma}
	\label{ukr}Suppose that the properties $(A)$ hold. If $\langle u, L^{- 1}_u (u) \rangle = 0$, for any $c \in ]c_{\ast},c^{\ast}[$ and
	$u \in S_c$, there exists $\kappa > 0$ such that if a function $\varphi \in
	H^1_{\tmop{rad}} (\mathbb{R}^n, \mathbb{R})$ satisfies
	\[ \langle \varphi, u \rangle = \langle \varphi, L^{- 1}_u (u) \rangle = 0,
	\]
	then
	\[ \langle L_u (\varphi), \varphi \rangle \geqslant \kappa \| \varphi
	\|_{H^1}^2 . \]
	Furthermore, if $\langle u, L^{- 1}_u (u) \rangle \neq 0$, then we have the
	following improvement: there exists $\kappa > 0$ such that if a function
	$\varphi \in H^1_{\tmop{rad}} (\mathbb{R}^n, \mathbb{R})$ satisfies
	\[ \langle \varphi, u \rangle = 0, \]
	then
	\[ \langle L_u (\varphi), \varphi \rangle \geqslant \kappa \| \varphi
	\|_{H^1}^2 . \]
\end{lemma}

This results use the fact that $L_u$ has no kernel and exactly one negative eigenvalue. Therefore, we can expect to have coercivity of $L_u$ under one orthogonality condition. This result state that this is indeed the case for function orthogonal to $u$ if  $\langle u, L^{- 1}_u (u) \rangle \neq 0$, and if this quantity is $0$, then we still have the coercivity but with a second orthogonality, on $ L^{- 1}_u(u)$.

\subsubsection{Extension of the branch}

We recall that $\tmop{Eq}_{\lambda} (u) = - \Delta u + V u - \lambda u - f (x,
u)$. The first part of the proof consists in proving that a non-degenerate
solution of $\tmop{Eq}_{\lambda} (u_{\lambda}) = 0$ generates a branch of
solutions of $\tmop{Eq}_{\lambda + \varepsilon} (U_{\lambda + \varepsilon}) =
0$ for $| \varepsilon |$ small. This is the following result.

\begin{proposition}[Construction of a local branch]
  \label{constbranch}Consider $\lambda < \Lambda_0 \leqslant 0$ and $u_{\lambda}$ a radial
  solution of $\tmop{Eq}_{\lambda} (u_{\lambda}) = 0$ that satisfies
  $\tmop{Ker} L_{u_{\lambda}} \cap L^2_{\tmop{rad}} = \{ 0 \}$ and $u'(r) < 0$. Then, there
  exists $\varepsilon_0 > 0$ (depending on $\lambda$ and $u_{\lambda}$) and a
  function
  \[ \varepsilon \rightarrow U_{\lambda + \varepsilon} \in C^1 (] -
     \varepsilon_0, \varepsilon_0 [, H^1_{\tmop{rad}}) \]
  with $U_{\lambda} = u_{\lambda}$ such that, for all $\varepsilon \in] -
  \varepsilon_0, \varepsilon_0 [$, we have
  \[ \tmop{Eq}_{\lambda + \varepsilon} (U_{\lambda + \varepsilon}) = 0. \]
  In particular, $\varepsilon \rightarrow E (U_{\lambda + \varepsilon}) \in
  C^1 (] - \varepsilon_0, \varepsilon_0 [, \mathbb{R})$ and $\varepsilon
  \rightarrow Q (U_{\lambda + \varepsilon}) \in C^1 (] - \varepsilon_0,
  \varepsilon_0 [, \mathbb{R})$. If $f$ is analytic w.r.t its second variable, then these two functions are also analytic for $\varepsilon_0$ small enough.
\end{proposition}

The proof of this result follows the classical Lyapounov-Schmidt approach, see
{\cite{MR1959647}} and some examples of applications in {\cite{MR4360608,delP_Kow_Mus}}. That is, $U_{\lambda + \varepsilon}$ is defined by the
implicit problem $\tmop{Eq}_{\lambda + \varepsilon} (U_{\lambda +
\varepsilon}) = 0$ next to a solution $u_{\lambda}$ for the case $\varepsilon
= 0$. The non-degeneracy of $L_{u_{\lambda}}$ ($\tmop{Ker} L_{u_{\lambda}}
\cap L^2_{\tmop{rad}} = \{ 0 \}$) means that we can use a result like the implicit
function theorem (here in a Hilbert space) to construct $U_{\lambda +
\varepsilon}$.

Remark that with the properties $(A)$, any element of $S_c$ for any $c \in ]c_{\ast},c^{\ast}[$
satisfies the hypotheses of Proposition \ref{constbranch}. However this result
can be applied to other solutions of the equation, not necessarly NGSS. Furthermore, if $u_{\lambda} \in S_{\tilde{c}}$ for some
$\tilde{c}> 0$, then the branch $U_{\lambda + \varepsilon}$ constructed from
it in Proposition \ref{constbranch} is not necessarly formed of NGSS (that is, it is possible that for $\varepsilon \neq 0$, $U_{\lambda
+ \varepsilon} \nin S_c$ whatever $c \in ]c_{\ast},c^{\ast}[$ is, despite the fact that
$U_{\lambda} = u_{\lambda} \in S_{\tilde{c}}$). However, we will show later
that NGSS with mass $c$ close to $\tilde{c}$ must be part of
a branch generated this way by an element of $S_c$, see Proposition
\ref{simo}.

\

Before proving that in a rigorous way, let us show some properties of this branch of solutions.

\subsubsection{Local uniqueness for the extended branch}

We define the linearized operator around $U_{\lambda + \varepsilon}$ as
follows:
\[ L_{U_{\lambda + \varepsilon}} (\varphi) \assign - \Delta \varphi + V
   \varphi - (\lambda + \varepsilon) \varphi - f_t (x, U_{\lambda +
   \varepsilon}) (\varphi) . \]
Let us show a coercivity result on this operator.

\begin{proposition}[Coercivity of $L_{U_{\lambda + \varepsilon}}$]
  \label{coerc2}Under the notations of Proposition \ref{constbranch}, suppose
  furthermore that $u_{\lambda}$ satisfies the following property:
  \begin{itemizeminus}
  	\item $(E 1)$ There exists $\kappa > 0$ such that, for $\varphi \in
  	H^1_{\tmop{rad}}$, if $\langle \varphi, u_{\lambda} \rangle = \langle
  	\varphi, L^{- 1}_{u_{\lambda}} (u_{\lambda}) \rangle = 0$ then
  	\[ \langle L_{u_{\lambda}} (\varphi), \varphi \rangle \geqslant \kappa \|
  	\varphi \|_{H^1}^2 . \]
  \end{itemizeminus}
  In that case, up to reducing the value of $\varepsilon_0$, there exists
  $\kappa > 0$ such that for any $\varphi \in H^1_{\tmop{rad}}$ and $|
  \varepsilon | \leqslant \varepsilon_0$ we have
  \[ \langle L_{U_{\lambda + \varepsilon}} (\varphi), \varphi \rangle
  \geqslant \kappa \| \varphi \|_{H^1}^2 - \frac{1}{\kappa} \langle
  \varphi, U_{\lambda + \varepsilon} \rangle^2 - \frac{1}{\kappa} \langle
  \varphi, L^{- 1}_{U_{\lambda + \varepsilon}} (U_{\lambda + \varepsilon})
  \rangle^2 . \]
  If instead we suppose the following property:
  \begin{itemizeminus}
  	\item $(E 2)$ There exists $\kappa > 0$ such that, for $\varphi \in
  	H^1_{\tmop{rad}}$, if $\langle \varphi, u_{\lambda} \rangle = 0$ then
  	\[ \langle L_{u_{\lambda}} (\varphi), \varphi \rangle \geqslant \kappa \|
  	\varphi \|_{H^1}^2 . \]
  \end{itemizeminus}
  In that case, up to reducing the value of $\varepsilon_0$, there exists
  $\kappa > 0$ such that for any $\varphi \in H^1_{\tmop{rad}}$ and $|
  \varepsilon | \leqslant \varepsilon_0$ we have
  \[ \langle L_{U_{\lambda + \varepsilon}} (\varphi), \varphi \rangle
  \geqslant \kappa \| \varphi \|_{H^1}^2 - \frac{1}{\kappa} \langle
  \varphi, U_{\lambda + \varepsilon} \rangle^2 . \]
  Finally, in both cases, we have the following identities:
  \[ \partial_{\varepsilon} (Q (U_{\lambda + \varepsilon}))_{| \varepsilon = 0
  } = \langle L^{- 1}_{u_{\lambda}} (u_{\lambda}), u_{\lambda}
  \rangle, \]
  \[ L_{u_{\lambda}} (\partial_{\varepsilon} U_{\lambda + \varepsilon |
  	\varepsilon = 0 }) = u_{\lambda} \]
  and
  \[ \partial_{\varepsilon} (E (U_{\lambda + \varepsilon}) - \lambda Q
  (U_{\lambda + \varepsilon}))_{| \varepsilon = 0 } = 0. \]
  The second identity implies in particular that $L^{- 1}_{u_{\lambda}}
  (u_{\lambda}) = \partial_{\varepsilon} U_{\lambda + \varepsilon |
  	\varepsilon = 0 }$.
\end{proposition}

As for the previous proposition, with the properties $(A)$, in Lemma \ref{ukr}
we proved that any element of $S_c$ for any $c \in ]c_{\ast},c^{\ast}[$ satisfies either $(E 1)$ and
$\langle u_{\lambda}, L^{- 1}_{u_{\lambda}} (u_{\lambda}) \rangle = 0$ or $(E
2)$ and $\langle u_{\lambda}, L^{- 1}_{u_{\lambda}} (u_{\lambda}) \rangle \neq
0$. This result shows that if $u_{\lambda}$ has only one negative eigenvalue and
an empty kernel, then so does the element of the branch $U_{\lambda +
\varepsilon}$ generated by it. This can be interpreted as the continuity of
the spectrum of $L_{U_{\lambda + \varepsilon}}$ with respect to $\varepsilon$.

Using this result, we will show a local uniqueness result on this branch among
solutions of $\tmop{Eq}_{\lambda}$, which is as follows:

\begin{proposition}[Local uniqueness]
  \label{uniquresult}Under the hypotheses and notations of Proposition
  \ref{constbranch}, supposing furthermore that either we have $(E 1)$
  and $\langle u_{\lambda}, L^{- 1}_{u_{\lambda}} (u_{\lambda}) \rangle = 0$,
  or we have $(E 2)$ and $\langle u_{\lambda}, L^{- 1}_{u_{\lambda}}
  (u_{\lambda}) \rangle \neq 0$, then there exists $\eta > 0$ small such that,
  if a radial function $\psi$ satisfies
  \[ \tmop{Eq}_{\lambda_0} (\psi) = 0 \]
  for some $\lambda_0 < 0$ and $\| \psi - u_{\lambda} \|_{H^1 (\mathbb{R}^n)}
  \leqslant \eta$, then there exists $\varepsilon \in \mathbb{R}$ such that
  \[ \psi = U_{\lambda + \varepsilon} . \]
\end{proposition}

Remark that here $\lambda_0$ might be far away from $\lambda$ a priori.
However we will show in the proof of Proposition \ref{uniquresult} that $\|
\psi - u_{\lambda} \|_{H^1 (\mathbb{R}^n)} \leqslant \eta$ implies in fact
that $\lambda_0$ and $\lambda$ are close if $\eta$ is small. The proof of this
result uses modulation theory, see for instance {\cite{MR1888800,MR820338}}.

This is a uniqueness result in the following sense: Among radial functions
that are solution of $\tmop{Eq}_{\lambda_0}$ for any $\lambda_0 < 0$, the only
ones that are close to $u_{\lambda}$ are the elements of the branch
$U_{\lambda + \varepsilon}$ from Proposition \ref{constbranch}.

\

The three propositions above are in some sense disconnected to the
minimization problem, as they do not require $u_{\lambda}$ to be a NGSS.
However, they can be applied to elements of $S_c$. Proposition
\ref{uniquresult} will be particularly useful with the property $(A 3)$ to
propagate the uniqueness result. The idea of using this kind of construction not
connected to the minimization to study uniqueness of the minimization problem
was first used in {\cite{Chi_Pac_3,MR4360608,MR4522936}}
(and see {\cite{SLSEDP_2021-2022____A2_0}} for an overview of these papers).

Now, we use the above results to show some properties on the elements of $\mathcal{T}_c$ and $S_c$. Let us start by proving some crucial properties of $m(c)$ and $\mathcal{T}_c$.

\subsubsection{Properties of $m(c)$ and $\mathcal{T}_c$}

We recall the definitions
\[ m (c) = \inf_{u \in X_c} E (u), \quad X_c = \{u \in H^1, Q(u) = c\}, \]
\[ S_c = \{ u \in X_c, E (u) = m (c) \}, \]
and
\[ \mathcal{T}_c = \{\lambda, \exists u \in S_c \ s.t. \ \tmop{Eq}_{\lambda} (u) = 0\} .\]

\begin{proposition}[Properties of $m(c)$ and $\mathcal{T}_c$]  \label{properties of m(c)}
	Suppose that $(A 1)$ and $(A 3)$ hold.
	\begin{itemize}
		\item[(1)] We have
		$$
		m'_{+}(c) = \min\mathcal{T}_c, \quad m'_{-}(c) = \max\mathcal{T}_c,
		$$
		and
		$$
		\lim_{c \to \tilde{c}^+}\min\mathcal{T}_c = \lim_{c \to \tilde{c}^+}\max\mathcal{T}_c = \min\mathcal{T}_{\tilde{c}},
		$$
		$$
		\lim_{c \to \tilde{c}^-}\min\mathcal{T}_c = \lim_{c \to \tilde{c}^-}\max\mathcal{T}_c = \max\mathcal{T}_{\tilde{c}}.
		$$
		Thus $m(c)$ is decreasing, Lipschitz continuous and differentiable at almost everywhere $c \in ]c_{\ast},c^{\ast}[$.
		\item[(2)]  $m(c)$ is strictly concave down on $c \in ]c_{\ast},c^{\ast}[$, it is then twice differentiable at almost everywhere $c \in ]c_{\ast},c^{\ast}[$. Furthermore, letting
		$$
		\lambda_{min}(c) = \min\mathcal{T}_c, \quad \lambda_{max}(c) = \max\mathcal{T}_c,
		$$
		then
		$$
		\lambda_{min}(c_1)\leq \lambda_{max}(c_1) < \lambda_{min}(c_2) \leq \lambda_{max}(c_2), \forall c_1 > c_2.
		$$
	\end{itemize}
 Finally, 
 if we further assume that $(D 1)$ holds, then
	\begin{align} \label{eq right res}
		\liminf_{\tilde{c} \to c^{+}}\frac{\lambda_{min}(\tilde{c})-\lambda_{min}(c)}{\tilde{c}-c} > -\infty.
	\end{align}
	Similarly, 
	if we further assume that $(D 2)$ holds, then
	\begin{align} \label{eq left res}
		\liminf_{\tilde{c} \to c^-}\frac{\lambda_{max}(\tilde{c})-\lambda_{max}(c)}{\tilde{c}-c} > -\infty.
	\end{align}
\end{proposition}	

The result (1) can be deduced using similar discussions to \cite{{Hajaiej-Song-unique}}. That $m(c)$ is strictly concave down in result (2) can be proved by borrowing an idea of \cite[Section 6.3]{Lewin-Nodari-2020}, in which $c$ is constrained on an interval $(a,b)$ and that $|\mathcal{T}_c| = 1$, $\lambda(c)$ is continous on $(a,b)$ are assumed. We obtain stronger results without these constraints. We also show the monotonicity of $\mathcal{T}_c$ w.r.t. $c \in ]c_{\ast},c^{\ast}[$ and the boundedness of Dini derivative, this is crucial for the following steps, particularly, for the proofs of $\langle L^{-1}_uu,u\rangle < 0$ and the finiteness of $\mathcal{P}$. Though some parts of the proofs are not new, we provide all the details in subsubsection \ref{sec. properties of m(c)} for the sake of consistency in the article and for the readers' convenience.

\subsubsection{Results on the sign of $\langle L^{-1}_{u_c} u_c, u_c \rangle$ for $u_c \in S_c$}

Using Proposition \ref{properties of m(c)}, we can show the negativity of $\langle L^{-1}_{u_c} u_c, u_c \rangle$ for $u_c \in S_c$ when $(D)$ holds true.

\begin{lemma} \label{< 0}
We assume that properties $(A)$ hold. 

If $(D 1)$ holds, then $\langle L^{-1}_{u_c} u_c, u_c \rangle < 0$ where $u_c$ is such that $\tmop{Eq}_{\lambda_{min}(c)}(u_c) = 0$.

If $(D 2)$ holds, then $\langle L^{-1}_{u_c} u_c, u_c \rangle < 0$ where $u_c$ is such that $\tmop{Eq}_{\lambda_{max}(c)}(u_c) = 0$.
\end{lemma}

When $f(x,u)$ is real analytic w.r.t. $u$, we can show that except for finite values of $c$ that $\langle L_u^{-1}u,u\rangle \neq 0$ holds for $u \in S_c$.

\begin{lemma} \label{=0 is finite}
	Assume that the properties $(A)$ hold and that $f(x,u)$ is real analytic w.r.t. $u$. Then, for any $c_0,c^0 \in ]c_{\ast},c^{\ast}[$, the value of $c$ in $[c_0,c^0]$ such that $\langle L_u^{-1}u,u\rangle = 0$ for some $u \in \mathcal{S}_c$ is finite.
\end{lemma}

\subsubsection{Proof of Theorem \ref{fleurs}}

We have now enough results to conclude the proof of Theorem \ref{fleurs}.

\

\begin{proof}
	We first address the case when the conditions $(C)$ and $(D)$ hold. We argue by contradiction. Suppose that $S_{\tilde{c}}$ is infinite for some $\tilde{c} > 0$ and take $(u_k)_{k \in \mathbb{N}}$ a sequence of distincts elements of $S_{\tilde{c}}$. By $(A 1)$, we denote by $\lambda_k$ their
	Lagrange multipliers (that is $\tmop{Eq}_{\lambda_k} (u_k) = 0$). By $(A
	3)$, there exists $u_{\infty} \in S_{\tilde{c}}$ such that, up to a
	subsequence, $u_k \rightarrow u_{\infty}$ strongly in $H^1 (\mathbb{R}^n)$.
	we denote by $\lambda_{\infty}$ the Lagrange multiplier of $u_{\infty}$, and
	thus $\tmop{Eq}_{\lambda_{\infty}} (u_{\infty}) = 0$. Note that $\langle L^{-1}_{u_{\infty}}u_{\infty},u_{\infty}\rangle \neq 0$ by Lemma \ref{< 0}.
	
	Since $u_{\infty} \in S_c$, by the properties $(A)$ we construct by
	Proposition \ref{constbranch} the branch $U_{\lambda_{\infty} +
		\varepsilon}$, solution of $\tmop{Eq}_{\lambda_{\infty} + \varepsilon}
	(U_{\lambda_{\infty} + \varepsilon}) = 0$ for $| \varepsilon |$ small.
	Furthermore, since $\partial_{\varepsilon} (Q (U_{\lambda_{\infty} +
		\varepsilon}))_{| \varepsilon = 0 } \neq 0$ by Proposition
	\ref{coerc2}, for $| \varepsilon |$ small enough, $Q (U_{\lambda_{\infty} +
		\varepsilon}) = \tilde{c}$ if and only if $\varepsilon = 0$.
	
	Now, we have shown that $\| u_k - u_{\infty} \|_{H^1 (\mathbb{R}^n)} = o_{k
		\rightarrow + \infty} (1)$, hence by Proposition \ref{uniquresult}, $u_k =
	U_{\lambda_{\infty} + \varepsilon_k}$ for $k$ large enough and a family of
	small $(\varepsilon_k)$. But since $Q (u_k) = \tilde{c}$, we have
	$\varepsilon_k = 0$ for $k$ large enough, which is in contradiction with the
	fact that the $(u_k)_{k \in \mathbb{N}}$ are distincts from one another.
	
	Now we focus on the case when $f(x,u)$ is real analytic w.r.t. $u$. Noticing that $Q (U_{\lambda_{\infty} +
		\varepsilon})$ is strictly monotonous using its real analyticity. We can complete the proof using similar arguments to the case when the conditions $(C)$ and $(D)$ hold. Indeed, we can show that the energy is also independent of $\lambda$ and then we can extend the branch for all values of $\lambda$ which is in contradiction with the fact that $\lambda_{max}$ is bounded.
\end{proof}

\subsubsection{Extension by NGSS and its applications}

We are now equiped to show that the functions constructed in Proposition \ref{constbranch} from suitable elements are exactly NGSS in small interval.

\begin{proposition}
	\label{simo}
	
	Assume that properties $(A)$ hold. Let $$(u_{j,min})_{j \in [1 \ldots N_1]}, N_1 \geqslant 1; \quad (u_{j,max})_{j \in [1 \ldots N_2]}, N_2 \geqslant 1$$ be all the elements in $S_{\tilde{c}}$ such that $\tmop{Eq}_{\lambda_{min}(\tilde{c})}u_{j,min} = 0$ and $\tmop{Eq}_{\lambda_{max}(\tilde{c})}u_{j,max} = 0$ respectively. Both sets are non-empty and contain at most a finite number of elements.
	
	Furthermore, there exists $\varepsilon_0 > 0$ (depending on $\tilde{c}$) such
	that all $u_{j,min}, j \in [1 \ldots N_1]$ and $u_{j,max}, j \in [1 \ldots N_2]$ can be extended to branches denotes
	$U_{j, \varepsilon, min}$ and $U_{j, \varepsilon, max}$ for $| \varepsilon | \leqslant \varepsilon_0$ following
	the construction of Proposition \ref{constbranch} (with $U_{j, 0, min} = u_{j,min}$ and $U_{j, 0, max} = u_{j,max}$).
	
	Then, there exists $\tau_0 > 0$ (depending on $\tilde{c}$) such that for any $c
	\in ]\tilde{c} , \tilde{c} + \tau_0]$, there exists $(\varepsilon_j)_{j \in
		[1 \ldots N_1]} \in [- \varepsilon_0, 0[$ such that $S_c$, the set
	of minimizers of the energy for the mass $c$ satisfies
	\[ S_c \subset \bigcup_{j \in [1 \ldots N_1]} \{ U_{j, \varepsilon_j,min} \}; \]
	for any $c
	\in [\tilde{c} - \tau_0 , \tilde{c}[$, there exists $(\varepsilon_j)_{j \in
		[1 \ldots N_2]} \in ]0,\varepsilon_0]$ such that $S_c$, the set
	of minimizers of the energy for the mass $c$ satisfies
	\[ S_c \subset \bigcup_{j \in [1 \ldots N_2]} \{ U_{j, \varepsilon_j,max} \}; \]
\end{proposition}

The core of the proof of this result is to use Propositions \ref{uniquresult}, \ref{properties of m(c)}
with $(A 3)$. By $(A 3)$ and Proposition \ref{properties of m(c)}, elements of $S_c$ must be close to elements of
$S_{\tilde{c}}$ with Lagrange multiplier being $\lambda_{min}(\tilde{c})$ or $\lambda_{max}(\tilde{c})$ if $c$ is close to $\tilde{c}$. The only possibility to do so by Proposition \ref{uniquresult} is to be of the form $U_{j, \varepsilon_j,min}$ or $U_{j, \varepsilon_j,max}$.




\

If we further assume that $(C)$ and $(D)$ hold, we can improve the above result. 

\begin{proposition}[Extension by NGSS, I] \label{ext 1}
	Assume that properties $(A)$, $(C)$ and $(D)$ hold. For any $\tilde{c} > 0$, let $u_{min}, u_{max}$ be the unique elements in $S_{\tilde{c}}$ with Lagrange multipliers $\lambda_{min}(\tilde{c}), \lambda_{max}(\tilde{c})$ respectively. Define $U_{\varepsilon, min}$, $U_{\varepsilon, max}$ the branches of Proposition \ref{constbranch} associated to them.
	\begin{itemizeminus}
		\item There exists $\tau_0 > 0$ such that for $c \in ]\tilde{c}, \tilde{c} + \tau_0]$, $|S_c| = 1$ and $\{u_c, c \in ]\tilde{c}, \tilde{c} + \tau_0]\} = \{U_{\varepsilon, min}, \varepsilon \in [-\varepsilon_0,0[\}$ where $\varepsilon_0 > 0$, $U_{-\varepsilon_0,min} = u_{\tilde{c} + \tau_0}$ and $u_c$ is the unique element in $S_c$.
		\item Similarly, There exists $\tau_0 > 0$ such that for $c \in [\tilde{c} - \tau_0, \tilde{c}[$, $|S_c| = 1$ and $\{u_c, c \in [\tilde{c} - \tau_0, \tilde{c}[\} = \{U_{\varepsilon, max}, \varepsilon \in ]0, \varepsilon_0]\}$ where $\varepsilon_0 > 0$, $U_{\varepsilon_0,max} = u_{\tilde{c} - \tau_0}$ and $u_c$ is the unique element in $S_c$.
	\end{itemizeminus}
\end{proposition}

If instead we assume that $(N)$ holds and $f$ is analytic, with have the following result.

\begin{proposition}[Extension by NGSS, II] \label{ext 2}
	Assume that properties $(A)$, $(N)$ hold and that $f$ is real analytic w.r.t. $u$.  Let $$(u_{j,min})_{j \in [1 \ldots N_1]}, N_1 \geqslant 1; \quad (u_{j,max})_{j \in [1 \ldots N_2]}, N_2 \geqslant 1$$ be elements in $S_{\tilde{c}}$ such that $\tmop{Eq}_{\lambda_{min}(\tilde{c})}u_{j,min} = 0$ and $\tmop{Eq}_{\lambda_{max}(\tilde{c})}u_{j,max} = 0$ respectively. Both sets are non-empty and contain at most a finite number of elements. Define $U_{j, \varepsilon, min}, j \in [1 \ldots N_1]$, $U_{j, \varepsilon, max}, j \in [1 \ldots N_2]$ the branches of Proposition \ref{constbranch} associated to them.
	\begin{itemizeminus}
		\item There exists $\tau_0 > 0$ such that for $c \in ]\tilde{c}, \tilde{c} + \tau_0]$, $|S_c| = L \in [1,N_1]$ and $\{u \in S_c, c \in ]\tilde{c}, \tilde{c} + \tau_0]\} = \{U_{j, \varepsilon, min}, \varepsilon \in [-\varepsilon_0,0[, j \in [1 \ldots L]\}$ where $\varepsilon_0 > 0$, $U_{j,-\varepsilon_0,min} \in S_{\tilde{c} + \tau_0}$ for any $j \in [1 \ldots L]$.
		\item Similarly, There exists $\tau_0 > 0$ such that for $c \in [\tilde{c} - \tau_0, \tilde{c}[$, $|S_c| = L \in [1,N_2]$ and $\{u \in S_c, c \in [\tilde{c} - \tau_0, \tilde{c}[\} = \{U_{j,\varepsilon, max}, \varepsilon \in ]0, \varepsilon_0], j \in [1 \ldots L]\}$, where $\varepsilon_0 > 0$, $U_{j, \varepsilon_0,max} \in S_{\tilde{c} - \tau_0}$ for any $j \in [1 \ldots L]$.
	\end{itemizeminus}	
\end{proposition}

\subsection{Conclusion}

We can now conclude the proofs of Theorems \ref{mainthing0}-\ref{mainthing2}.

\subsubsection{End of the proof of Theorem \ref{mainthing0}}

	We first address the problem when $(C)$ and $(D)$ hold. We argue by contradiction. Take $(c_k)_{k \in \mathbb{N}}$ a sequence of distincts elements such that $(c_k)_{k \in \mathbb{N}} \subset \mathcal{P}\cap [c_0,c^0]$. Up to a subsequence, assume that $c_k \to c \in [c_0,c^0]$. Note that either $c_k \to c^+$ or $c_k \to c^-$ passing to a subsequence if necessary. Assume the former holds (the proof when the latter holds is similar). Let $\lambda_k = \lambda_{min}(c_k)$ and $\tmop{Eq}_{\lambda_k}(u_k) = 0$, $\tilde{\lambda}_k = \lambda_{max}(c_k)$ and $\tmop{Eq}_{\tilde{\lambda}_k}(\tilde{u}_k) = 0$. By Proposition \ref{properties of m(c)}, we have $\lambda_k, \tilde{\lambda}_k \to \lambda := \lambda_{min}(c)$. Let $u$ be the unique element in $S_c$ such that $\tmop{Eq}_{\lambda} (u) = 0$. Then $u_k, \tilde{u}_k \rightarrow u$ strongly in $H^1$. That $\lambda_k < \tilde{\lambda}_k$ contradicts Proposition \ref{ext 1}.
	
	Now we assume that $f(x,u)$ is real analytic w.r.t $u$. Suppose on the contrary that $(c_k)_{k \in \mathbb{N}}$ is a sequence of distincts elements such that $(c_k)_{k \in \mathbb{N}} \subset \mathcal{P}\cap [c_0,c^0]$. Up to a subsequence, assume that $c_k \to c \in [c_0,c^0]$. Similar to the case when $(C)$ and $(D)$ hold, we can find a contradiction to Proposition \ref{ext 2} when $L = 1$ where $L$ is given by Proposition \ref{ext 2}. When $L > 1$, using real analyticity, solutions corresponding to the same $\lambda$ on different branches constructed by Proposition \ref{ext 2} have same masses and energies (see more details in the proof of Proposition \ref{ext 2}), thus we can get a similar contradiction as the case that $L = 1$, which completes the proof. 

\

 Adding property $(B)$, it is easy to prove that $\mathcal{P}$ is finite.

\subsubsection{End of the proof of Theorem \ref{Planitia}}

 The continuity and strict decreaseness of $\mathcal{T}_c$ are proved in Proposition \ref{properties of m(c)}. Consequently, we have $\{\lambda \in \mathcal{T}_c: c \in ] c_j,c_{j+1}[\} = ] \max\mathcal{T}_{c_{j+1}},\min\mathcal{T}_{c_j}[$. For the final equivalence in Theorem \ref{Planitia}, previous results show the direct implication, and the reciprocal is a consequence of Proposition \ref{properties of m(c)} (1): Since $m'_-(c) = \max\mathcal{T}_c$ and $m'_+(c) = \min\mathcal{T}_c$, $m$ is differentable  at $c$ $\Leftrightarrow$ $\max\mathcal{T}_c = \min\mathcal{T}_c \Leftrightarrow |\mathcal{T}_c| = 1$.

\subsubsection{End of the proof of Theorem \ref{mainthing1}}

The facts that $|S_c| = 1 \Leftrightarrow |\mathcal{T}_c| = 1$ and that $\mathcal{M} = \mathcal{P}$ are a direct consequence of $(C)$. The fact that $S_c$ for any $c \in ]c_{\ast},c^{\ast}[$ assuming $\langle L^{-1}_u(u),u \rangle \neq 0$ for any $u \in S_c$ is a finite set is Lemma \ref{fleurs}.  We now consider the proof of $c \rightarrow u_c \in C^1 (] c_j, c_{j + 1} [, H^1 (\mathbb{R}^n))$ in Theorem \ref{mainthing1}. Let $u$ be the unique element in $S_{c_j}$ with $\tmop{Eq}_{\lambda_{min}(c_j)}(u) = 0$. By Proposition
\ref{constbranch} we can construct $U_{\lambda_{min}(c_j) + \varepsilon}$ for $\varepsilon < 0$, $|\varepsilon |$ small that extend $u$, with $\varepsilon \rightarrow U_{\lambda_{min}(c_j) + \varepsilon} \in C^1 ( [- \varepsilon_0, 0], H^1_{\tmop{rad}} (\mathbb{R}^n))$. By Proposition \ref{ext 1} there exists $\tau_0 > 0$ such that $c_j+\tau_0 = Q(U_{\lambda_{min}(c_j) - \varepsilon_0})$ and for $c \in ]c_{j}, c_{j} + \tau_0]$, $|S_c| = 1$ and $\{u_c, c \in ]c_{j}, c_{j} + \tau_0]\} = \{U_{\lambda_{min}(c_j) + \varepsilon}, \varepsilon \in [-\varepsilon_0,0[\}$ where $\varepsilon_0 > 0$ and $u_c$ is the unique element in $S_c$. Repeating the process in Proposition \ref{ext 1}, $U_{\lambda_{min}(c_j) + \varepsilon}$ can be extended to $\lambda_{max}(c_{j+1})$ with $U_{\lambda_{max}(c_{j+1})} \in S_{c_{j+1}}$. Furthermore, for any $c \in ] c_j, c_{j + 1} [$, there exists a unique $\varepsilon_c \in ]\lambda_{max}(c_{j+1})-\lambda_{min}(c_j),0[$ such that $u_c = U_{\lambda_{min}(c_j) + \varepsilon_c}$. Let us show that $c
\rightarrow \varepsilon_c$ is $C^1$. Indeed, $\varepsilon_c$ is choosen in a
way such that $Q (U_{\lambda + \varepsilon_c}) = c$ (because $U_{\lambda +
\varepsilon_c} \in S_c$). Since $\varepsilon \rightarrow Q (U_{\lambda +
\varepsilon})$ is $C^1$ by Proposition \ref{constbranch}, $\partial_{\varepsilon} (Q (U_{\lambda +
\varepsilon})) \neq 0$ for any $\varepsilon \in ]\lambda_{max}(c_{j+1})-\lambda_{min}(c_j),0[$ by Proposition \ref{uniquresult} and Lemma \ref{< 0}, by the implicit
function theorem we have that $c
\rightarrow \varepsilon_c$ is $C^1$. We conclude by Proposition
\ref{constbranch}, stating that $\varepsilon_c \rightarrow U_{\lambda +
\varepsilon_c}$ is $C^1$.

\subsubsection{End of the proof of Theorem \ref{mainthing2}}

We start with the proof of $c \to u_{c,i} \in C^{\alpha} (] c_j, c_{j + 1} [, H^1 (\mathbb{R}^n)), i \in [0 \ldots K_j]$ for some $\alpha \in (0,1)$ in Theorem \ref{mainthing2}. Let $(u_i), i \in [0 \ldots N]$ be the elements in $S_{c_j}$ with $\tmop{Eq}_{\lambda_{min}(c_j)}(u_i) = 0$. By Proposition
\ref{constbranch} we can construct $U_{i,\lambda_{min}(c_j) + \varepsilon}$ for $\varepsilon < 0$, $|\varepsilon |$ small that extend $u_i$, with $\varepsilon \rightarrow U_{i,\lambda_{min}(c_j) + \varepsilon} \in C^1 ([ - \varepsilon_0, 0], H^1_{\tmop{rad}} (\mathbb{R}^n))$. By Proposition \ref{ext 2}, there exists $K_j \in [1 \ldots N]$ such that $u_i$ is isolated in $S := \{u \in S_c, c \in ]c_{\ast},c^{\ast}[\}$ for $i \in [K_j+1 \ldots N]$. Furthermore, for all $i \in [1 \ldots K_j]$, there exists $\tau_0$ such that $c_j +\tau_0 = Q(U_{i,\lambda_{min}(c_j) - \varepsilon_0})$ and for $c \in ]c_{j}, c_{j} + c_0]$, $\{u_c, c \in ]c_{j}, c_{j} + c_0]\} = \{U_{i,\lambda_{min}(c_j) + \varepsilon}, \varepsilon \in [-\varepsilon_0,0[, i \in [1 \ldots K_j]\}$ where $\varepsilon_0 > 0$. Repeating the process in Proposition \ref{ext 2} and using the real analyticity of $Q(U_{i,\lambda_{min}(c_j) + \varepsilon})$ w.r.t. $\varepsilon$, for all $i \in [1 \ldots K_j]$, $U_{i,\lambda_{min}(c_j) + \varepsilon}$ can be extended to $\lambda_{max}(c_{j+1})$ with $U_{i,\lambda_{max}(c_{j+1})} \in S_{c_{j+1}}$. Furthermore, for any $c \in ] c_j, c_{j + 1} [$, there exists a unique $\varepsilon_c \in ]\lambda_{max}(c_{j+1})-\lambda_{min}(c_j),0[$ such that $Q(U_{i,\lambda_{min}(c_j) + \varepsilon_c}) = c$. Except these $K_j$ branches, $\{u \in S_c, c \in ] c_j, c_{j + 1} [\}$ consists of at most for a finite number of isolated elements. Indeed, if there are infinitely many elements, there must be a consensation point. Then similar to the proof of Proposition \ref{ext 2} we can find a contradiction with the real analyticity. Let us show that $c \rightarrow \varepsilon_c$ is $C^\alpha$ for some $\alpha \in (0,1)$. Remark that $\varepsilon_c$ satisfies the implicit equation $c = Q (U_{i,\lambda_{min}(c_j) +
	\varepsilon_c})$ (we can fix $i \in [1 \ldots K_j]$ here since the value is the same for different $i \in [1 \ldots K_j]$). We have that $\varepsilon \rightarrow Q (U_{i,\lambda_{min}(c_j) +
	\varepsilon})$ is analytic, and thus there exists $k
\geqslant 1$ such that $\partial_{\varepsilon}^d (Q (U_{i,\lambda_{min}(c_j) +
	\varepsilon})) = 0$ for $1 \leqslant d < k$ and
$\partial_{\varepsilon}^k (Q (U_{\lambda_{min}(c_j) + \varepsilon})) \neq 0$. This is
because otherwise $Q (U_{i,\lambda_{min}(c_j) + \varepsilon})$ has the same values for all
small $\varepsilon$, and this leads to a contradiction as in the proof of
Theorem \ref{fleurs}. Therefore, nearby $c_j$,
\[ c = Q (i,U_{\lambda_{min}(c_j) + \varepsilon_c}) = \lambda_{min}(c_j) + \varepsilon_c^k
\frac{\partial_{\varepsilon}^k (Q (U_{\lambda + \varepsilon}))}{k!} +
O_{\varepsilon_c \rightarrow 0} (\varepsilon_c^{k + 1}). \]
The case at other values can be proved similarly and we deduce that $c \rightarrow \varepsilon_c$ is $C^{1 / k}$.

\

Finally, by the analycity we have that $\partial_{\varepsilon} (Q (U_{\lambda +
	\varepsilon})) = 0$ is only possible for a finite number of values of
$\varepsilon$ (see Lemma \ref{=0 is finite}). Outside of these points we have that $c \rightarrow
\varepsilon_c$ is $C^1$, and at this finite number of point we do the
procedure above. Therefore $c \rightarrow \varepsilon_c$ is $C^{\alpha}$ for
$\alpha$ the minimum of the $1 / k$ for the finite amount of point where
$\partial_{\varepsilon} (Q (U_{\lambda + \varepsilon})) = 0$. At a bad value $\tilde{c}$, from $\partial_{\varepsilon} (Q (U_{\lambda + \varepsilon})) = 0$ we deduce that $$\lim_{c \to \tilde{c}}\|\frac{u_{c}-u_{\tilde{c}}}{c-\tilde{c}}\|_{H^1} = \infty.$$
This concludes the proof of Theorem \ref{mainthing2}.

\subsubsection{Open problems and overview of the paper}

The establishment of the number of $c_k$ in Theorem \ref{Planitia} and their characterization seem to be extremely hard. We hope that this highly important aspect will be addressed in the future.

\

 Another interesting open problem is the orbital stability of ground states in the
case $c \in \mathcal{M}$. A consequence of our work is that $S_c$
is a finite set and its element are well separated. It is therefore likely
that the orbital stability holds for any element in $S_c$ separatly from each other, but this remains an open question.

\

The rest of this paper is separated in three sections. Section
\ref{behind} contains the proofs of Theorems \ref{fleurs}-\ref{mainthing2}, following step by
step the results stated in subsection \ref{sand}. Sections \ref{sec. conditions} and \ref{sono} contain the result and proof of that is properties $(A)$,
$(B)$,$(C)$ and $(N)$ are satisfied under some assumptions on $V$ and $f$.

\bigbreak
\noindent {\bf Acknowledgement:} E.P. is supported by Tamkeen under the NYU Abu Dhabi
Research Institute grant CG002.

\section{Proof of Theorems \ref{fleurs}-\ref{mainthing2}}\label{behind}

\subsection{Proof of two preliminary Lemmas}

We start with some elliptic estimates that we will use in the proof of Lemma \ref{bb} and Proposition \ref{constbranch}.

\begin{lemma}
  \label{simon}Given $\lambda < 0$ and $h \in L^2 (\mathbb{R}^n)$, the
  solution $\varphi$ of the equation
  \[ - \Delta \varphi - \lambda \varphi = h \]
  can be written as
  \[ \varphi = I_{\lambda} (| . |) \ast h \]
  where
  \[ I_{\lambda} (r) = 2 \pi \left( \frac{2 \pi}{\sqrt{| \lambda |} r}
     \right)^{\frac{n}{2} - 1} K_{\frac{n}{2} - 1} \left( \sqrt{| \lambda |} r
     \right) \]
  and $K_j$ is the modified Bessel function of the second kind. For $n = 1$ we have
  \[ I_{\lambda} (r) = \pi e^{- \sqrt{| \lambda |} r} . \]
  For $n \geqslant 2$, $K_{\frac{n}{2} - 1} \in C^{\infty} (\mathbb{R}^{+
  \ast}, \mathbb{R}^+)$ and
  \[ K_{\frac{n}{2} - 1} (r) \sim_{r \rightarrow 0} \left\{\begin{array}{l}
       - \ln (r) \tmop{if} n = 2\\
       \frac{\Gamma \left( \frac{n}{2} - 1 \right)}{2} \left( \frac{2}{r}
       \right)^{\frac{n}{2} - 1} \tmop{if} n > 2
     \end{array}\right. \]
  and
  \[ K_{\frac{n}{2} - 1} (r) \sim_{r \rightarrow + \infty} \sqrt{\frac{\pi}{2
     r}} e^{- r} . \]
  These equivalents can be differentiated to compute the equivalents of
  $K'_{\frac{n}{2} - 1} (r)$ at $r = 0$ and $r = + \infty$.
\end{lemma}

This is a classical result, see for instance {\cite{MR0167642}}. Remark in
particular that $K_{\frac{n}{2} - 1}, K'_{\frac{n}{2} - 1} \in L^1 (B (0, 1))$
despite that they are singular at $r = 0$.

\subsubsection{Proof of Lemma \ref{bb}}

\begin{proof} Consider
	\begin{equation} \label{eq1.1}
		\left\{
		\begin{array}{cc}
			-(u'' + \frac{n-1}{r}u') + V(r)u = \lambda u + f(r,u), \\
			u(r) \rightarrow 0 \ as \ r \rightarrow +\infty.
		\end{array}
		\right.
	\end{equation}
	Note that $\lim_{r \rightarrow \infty}V(r) = 0$, $\lambda < 0$, $\lim_{u \rightarrow 0}f(r,u)/u = 0$ uniformly with respect to $r$. If $u = u(r)$ is a solution satisfying that $u > 0, u'(r) < 0$, we aim to show that $u$ and $u'$ decay exponentially.
	
	Let $z = -u'/u > 0$. We first show that it is bounded as $r \rightarrow \infty$. Notice that
	$$
	z' = -\frac{u''}{u}+\frac{(u')^2}{u^2} = z^2 - \frac{n-1}{r}z + \frac{f(r,u)}{u} + \lambda - V.
	$$
	Choose $r_0$ large enough such that $(n-1)/r \leq \sqrt{1-\lambda}/2$ and $\frac{f(r,u)}{u} + \lambda - V \geq \lambda -1$ for all $r > r_0$. When $r > r_0$ and $z > 2\sqrt{1-\lambda}$, one gets that
	$$
	z' \geq \frac{1}{2}z^2 + \left(\frac{1}{2}z^2 - \frac{\sqrt{1-\lambda}}{2}z + \lambda -1\right) \geq \frac{1}{2}z^2,
	$$
	which is impossible (because then $z$ blows up at a finite value of $r$). Hence,
	$$
	\limsup_{r \rightarrow \infty}z(r) \leq 2\sqrt{1-\lambda}.
	$$
	Then by L'H\^{o}spital's rule,
	$$
	\lim_{r \rightarrow \infty}\left(\frac{u'}{u}\right)^2 = \lim_{r \rightarrow \infty}\frac{u''}{u} = -\lambda.
	$$
	Moreover, it follows that for any $\epsilon \in (0,-\lambda)$,
	$$
	-\frac{u'}{u} \geq \sqrt{-\lambda-\epsilon}
	$$
	provided that $r$ is large enough. Then we can get that
	$$
	\limsup_{r \rightarrow \infty}u(r)e^{\sqrt{-\lambda-\epsilon}r} < \infty.
	$$
Now, from Lemma \ref{simon}, We have
 $ u = I_{\lambda}(|.|) \ast (-V u -f(r,u)), $
therefore
$$ \nabla u = (\nabla I_{\lambda}(|.|)) \ast (-V u -f(r,u)). $$
Since the convolution of two integrable and exponentially decaying function is an exponential decaying function, we show the estimate on $\nabla u$. Taking an additional derivative and putting it on the term $-V u -f(r,u)$ allows us to conclude.

\end{proof}

\subsubsection{Proof of Lemma \ref{ukr}}

\begin{proof}
	We start with the first part of the lemma. Define
	$$
	\kappa := \inf\{\langle L_{u} (\varphi), \varphi \rangle: \varphi \in H^1_{\tmop{rad}}, \langle
	\varphi, u \rangle = 0, \langle \varphi, L^{- 1}_u (u) \rangle = 0, \|\varphi\|_{H^1} = 1\}.
	$$
	Using the fact that $u$ is a NGSS on $X_{c}$ and $\langle
	\varphi, u \rangle = 0$ one deduces that $\kappa \geq 0$. It is sufficient to prove that $\kappa = 0$ is impossible. Suppose on the contrary that $\kappa = 0$.

\

	Step 1: There exists a minimizer $v \in H^1_{\tmop{rad}}, \langle
	v, u \rangle = 0, \langle v, L^{- 1}_u (u) \rangle = 0, \|v\|_{H^1} = 1$ such that $\langle L_{u} (v), v \rangle = 0$.

 \
	
	Take a minimizing sequence $\varphi_k$ for $\kappa$. Up to a subsequence, assume that $\varphi_k \rightharpoonup \varphi_\infty$ weakly in $H^1$ and $L^2$, $\varphi_k \to \varphi_\infty$ strongly in $L^2_{loc}$. We first find a contradiction if $\varphi_\infty \equiv 0$. In fact, combining $f_t(x,u(x)) \to 0$ as $|x| \to \infty$ and the strong convergence in $L^2_{loc}$, one gets
	$$
	\int_{\mathbb{R}^n}(|\nabla \varphi_k|^2 + V|\varphi_k|^2-\lambda|\varphi_k|^2)dx = \langle L_{u} (\varphi_k), \varphi_k \rangle + \int_{\mathbb{R}^n}f_t(x,u(x))|\varphi_k|^2dx \to 0,
	$$
	in a contradiction with $\|\varphi_k\|_{H^1} = 1$. Hence, $\varphi_\infty \neq 0$. Furthermore, $\varphi_\infty \in H^1_{\tmop{rad}}$ and $\langle
	\varphi_\infty, u \rangle = \lim_{n \to \infty}\langle
	\varphi_k, u \rangle = 0$, $\langle
	\varphi_\infty, L^{- 1}_u (u) \rangle = \lim_{n \to \infty}\langle
	\varphi_k, L^{- 1}_u (u) \rangle = 0$. Then consider $v = \frac{\varphi_\infty}{\|\varphi_\infty\|_{H^1}}$. Combining $f_t(x,u(x)) \to 0$ as $|x| \to \infty$ and the strong convergence in $L^2_{loc}$, one gets $$\int_{\mathbb{R}^n}f_t(x,u(x))|\varphi_k|^2dx \to \int_{\mathbb{R}^n}f_t(x,u(x))|\varphi_\infty|^2dx.$$
	Thus, using the weak lower semi-continuity of norm, we obtain
	$$
	\langle L_{u} (v), v \rangle = \frac{1}{\|\varphi_\infty\|_{H^1}^2}\langle L_{u} (\varphi_\infty), \varphi_\infty \rangle \leq \frac{1}{\|\varphi_\infty\|_{H^1}^2}\lim_{n \to \infty}\langle L_{u} (\varphi_k), \varphi_k \rangle = 0.
	$$
	On the other hand, one gets $\langle L_{u} (v), v \rangle \geq 0$ using $\langle
	v, u \rangle = 0$. This shows that $v$ is a minimizer for $\kappa$.

\

	Step 2: $v = \nu L_{u}^{-1}u$ for some $\nu \neq 0$.

 \
	
	By Lagrange multiplier principle, we have
	$$
	L_{u} (v) = \mu(-\Delta v + v) + \nu u + \tau L^{-1}_u(u).
	$$
	Using $\langle L_{u} (v), v \rangle = 0$, $\langle u,v \rangle = 0$ and $\langle v, L^{- 1}_u (u) \rangle = 0$, we obtain that $\mu = 0$. Moreover, since $\langle u,L^{-1}_u(u) \rangle = 0$, from
	$$
	0 = \langle u,v \rangle = \langle L^{-1}_u(u),L_u(v) \rangle = \tau \langle L^{-1}_u(u),L^{-1}_u(u)\rangle,
	$$
	we deduce that $\tau = 0$. Hence $v = \nu L_{u}^{-1}u$, implying that $\nu \neq 0$ since $v \neq 0$.
	
	On the other hand, $\langle v,L_{u}^{-1}u\rangle = \nu \langle L_{u}^{-1}u,L_{u}^{-1}u\rangle = 0$, implying that $\nu = 0$. This is a contradiction.

\

	The proof for the second result, when $\langle u, L^{- 1}_u (u) \rangle \neq 0$, is quite similar. We define
	$$
	\kappa := \inf\{\langle L_{u} (\varphi), \varphi \rangle: \varphi \in H^1_{\tmop{rad}}, \langle
	\varphi, u \rangle = 0, \|\varphi\|_{H^1} = 1\}.
	$$
	and get a minimiser $v = \nu L_{u}^{-1}u$ satisfying $\nu \neq 0$. Then $\langle
	v, u \rangle = 0$ yields that $\langle L_{u}^{-1}u, u \rangle = 0$, which is a contradiction. The proof is complete.
\end{proof}

\subsection{Construction of the branch $U_{\lambda + \varepsilon}$}

The goal of this subsection is to show Proposition \ref{constbranch}. That is,
any non-degenerate solution of equation $\tmop{Eq}_{\lambda}$
(\ref{solarstone}) can be extended to a branch of solution for the equation
for $\tmop{Eq}_{\lambda + \varepsilon}$ with $| \varepsilon |$ small.

\subsubsection{Some usefull elliptic estimates}

For $k \in \mathbb{N}, a > 0$, we define the norm:
\[ \| f \|_{C^k (a)} \assign \sup_{| j | \leqslant k} \| \nabla^j f e^{a | x
   |} \|_{L^{\infty} (\mathbb{R}^n)} . \]
and $C^{k}(a)$ the associated space.

\begin{lemma}
  \label{doti}Given $\lambda < 0$, there exists $a_0 > 0,$ such that for any
  $0 < a \leqslant a_0$, there exists $K (a) > 0$ such that if $h \in C^0 (a)$
  and $\varphi$ is solution of $- \Delta \varphi - \lambda \varphi = h$, then
  \[ \| \varphi \|_{C^1 (a)} \leqslant K (a) \| h \|_{C^0 (a)} . \]
\end{lemma}

\begin{proof}
  First, we have $\varphi (x) = (I_{\lambda} (| . |) \ast h) (x)$ by Lemma \ref{simon} and
  \[ | I_{\lambda} (| . |) \ast h | (x) \leqslant \| h \|_{C^0 (a)} |
     I_{\lambda} (| . |) \ast e^{- a | . |} | (x) \],
  and the convolution of two locally integrable radial and exponentially
  decaying functions is a function that decay exponentially fast, at the same rate
  as the lowest one of the two functions. Therefore, by Lemma \ref{simon}, for $a
  > 0$ small enough (depending on $\lambda$), we have
  \[ | I_{\lambda} (| . |) \ast e^{- a | . |} | (x) \leqslant K (a) e^{- a | x
     |} \]
  hence $| \varphi (x) e^{a | x |} | \leqslant K (a) \| h \|_{C^0 (a)}$.
  Similarly,
  \[ | \nabla (I_{\lambda} (| . |) \ast h) | (x) = | (\nabla I_{\lambda} (| .
     |)) \ast h | (x) \leqslant \| h \|_{C^0 (a)} | (\nabla I_{\lambda} (| .
     |)) \ast e^{- a | . |} | (x) \]
  and we conclude similarly that $| \nabla \varphi (x) e^{a | x |} | \leqslant
  K (a) \| h \|_{C^0 (a)}$.
\end{proof}

\subsubsection{Inversion of the linearized operator $L_{u_{\lambda}}$}

\begin{lemma}
  \label{cerasus}Under the assumptions and notation of Proposition
  \ref{constbranch}, there exists $K > 0$ such that for any $\varphi \in
  H^2_{\tmop{rad}} (\mathbb{R}^n), h \in L^2_{\tmop{rad}} (\mathbb{R}^n)$
  solution of
  \[ L_{u_{\lambda}} (\varphi) = h, \]
  we have
  \[ \| \varphi \|_{H^2 (\mathbb{R}^n)} \leqslant K \| h \|_{L^2
     (\mathbb{R}^n)} . \]
\end{lemma}

\begin{proof}
  We argue by contradiction. Suppose that there exists $h_n \in
  L^2_{\tmop{rad}} (\mathbb{R}^n)$ with $\| h_n \|_{L^2 (\mathbb{R}^n)}
  \rightarrow 0$ when $n \rightarrow + \infty$ and $\varphi_n \in
  H^2_{\tmop{rad}} (\mathbb{R}^n)$ such that $\| \varphi_n \|_{H^2
  (\mathbb{R}^n)} = 1$ with
  \[ L_{u_{\lambda}} (\varphi_n) = h_n . \]
  Since $\| \varphi_n \|_{H^2} = 1$, up to a subsequence,
  we have $\varphi_n \rightarrow \Phi$ strongly in $H^1$ on any compact set, where $\Phi$ is a
  radial function.

  Furthermore, since $L_{u_{\lambda}} (\varphi_n)$ converge strongly to $0$ in
  $L^2 (\mathbb{R}^n)$, we deduce the convergence of $\varphi_n$ to $\Phi$ is
  strong in $H^2$ (still on compact sets). Taking $n \rightarrow + \infty$, we
  have $L_{u_{\lambda}} (\Phi) = 0$, and since $\tmop{Ker} L_{u_{\lambda}}
  \cap L^2_{\tmop{rad}} = \{ 0 \}$ from the hypotheses of Proposition
  \ref{constbranch}, we have $\Phi = 0$. This implies in particular that
  $\varphi_n$ converges weakly to $0$ in $H^1 (\mathbb{R}^n)$ (since it
  converges strongly to $0$ in $H^2$ in any compact set).

  Now remark that $L_{u_{\lambda}} (\varphi) = - \Delta \varphi - \lambda
  \varphi + R (\varphi)$ where
  \[ R (\varphi) \assign V \varphi - f_t (x, u_{\lambda}) (\varphi) \]
  and since $\lambda < 0$, by Lemma \ref{bb}, $(- \Delta - \lambda)^{- 1} R$
  is a compact operator. Since $\varphi_n \rightharpoonup 0$ weakly in $H^1
  (\mathbb{R}^n)$, we have $\langle R (\varphi_n), \varphi_n \rangle = o_{n
  \rightarrow + \infty} (1)$ \ and thus
  \[ \langle L_{u_{\lambda}} (\varphi_n), \varphi_n \rangle = K (\lambda) \|
     \varphi_n \|_{H^1 (\mathbb{R}^n)}^2 + \langle R (\varphi_n), \varphi_n
     \rangle = K (\lambda) \| \varphi_n \|_{H^1 (\mathbb{R}^n)}^2 + o_{n
     \rightarrow + \infty} (1) \]
  while
  \[ | \langle L_{u_{\lambda}} (\varphi_n), \varphi_n \rangle | = | \langle
     h_n, \varphi_n \rangle | \leqslant \| h_n \|_{L^2 (\mathbb{R}^n)} \|
     \varphi_n \|_{H^2 (\mathbb{R}^n)} = o_{n \rightarrow + \infty} (1), \]
  hence $\| \varphi_n \|_{H^1 (\mathbb{R}^n)}^2 = o_{n \rightarrow + \infty}
  (1)$ and since $L_{u_{\lambda}} (\varphi_n) = h_n$ we deduce that
  \[ \| \Delta \varphi_n \|_{L^2 (\mathbb{R}^n)} = o_{n \rightarrow + \infty}
     (1), \]
  implying that $\| \varphi_n \|_{H^2 (\mathbb{R}^n)}^2 = o_{n \rightarrow +
  \infty} (1)$, which is in contradiction with $\| \varphi_n \|_{H^2
  (\mathbb{R}^n)} = 1$.
\end{proof}

\begin{lemma}
  Under the assumptions and notation of Proposition \ref{constbranch}, given
  $h \in L^2_{\tmop{rad}} (\mathbb{R}^n)$, there exists a unique function
  $\varphi \in H^2_{\tmop{rad}} (\mathbb{R}^n)$ such that
  \[ L_{u_{\lambda}} (\varphi) = h. \]
\end{lemma}

\begin{proof}
  First, take $R > 0$ and define the spaces $L^2_{\tmop{rad}, 0} (B (0, R)),
  H^1_{\tmop{rad}, 0} (B (0, R))$ to be respectively the sets of $L^2, H^1$
  radial functions on the ball $B (0, R)$ of $\mathbb{R}^n$ that vanish at
  the boundary.

  Consider the problem $L_{u_{\lambda}} (\varphi) = h$ on $B (0, R)$ with $h
  \in L^2_{\tmop{rad}, 0} (B (0, R))$ and $\varphi \in H^1_{\tmop{rad}, 0} (B
  (0, R))$. Let us show that if $R$ is large enough, then any solution of this
  problem satisfy $\| \varphi \|_{H^1 (B (0, R))} \leqslant K \| h \|_{L^2 (B
  (0, R))}$ for a constant $K > 0$ independent of $R$.

  We argue by contradiction. If this inequality does not hold, then there
  exists $\varphi_n \in H^1_{\tmop{rad}, 0} (B (0, R_n)), h_n \in
  L^2_{\tmop{rad}, 0} (B (0, R_n)), R_n > 0$ with $R_n \rightarrow + \infty$,
  $\| h_n \|_{L^2 (B (0, R_n))} \rightarrow 0$ when $n \rightarrow + \infty$
  and $\| \varphi_n \|_{H^1 (B (0, R_n))} = 1$ such that $L_{u_{\lambda}}
  (\varphi_n) = h_n$. We reach a contradiction with arguments similar to the
  ones of the proof of Lemma \ref{cerasus}.

  \

  From this a priori estimate, by Fredholm alternative, we deduce that for
  any $R > 0$ large enough and $f \in L^2_{\tmop{rad}, 0} (B (0, R))$, there
  exists $\varphi_R \in H^1_{\tmop{rad}, 0} (B (0, R))$ such that
  $L_{u_{\lambda}} (\varphi_R) = f$. Applying this to $f = h \chi_R$ where $h
  \in L^2_{\tmop{rad}} (\mathbb{R}^n)$ and $\chi_R$ is a radial smooth cutoff
  function with $\chi_R (x) = 1$ if $| x | \leqslant R - 2$ and $\chi_R (x) =
  0$ if $| x | \geqslant R - 1$, there exists $\varphi_R \in H^1_{\tmop{rad},
  0} (B (0, R))$ such that $L_{u_{\lambda}} (\varphi_R) = h \chi_R$, and
  \[ \| \varphi_R \|_{H^1 (B (0, R))} \leqslant K \| \chi h \|_{L^2 (B (0,
     R))} \leqslant K \| h \|_{L^2 (\mathbb{R}^n)} . \]
  By a diagonal argument, we deduce the existence of $\varphi \in
  H^1_{\tmop{rad}} (\mathbb{R}^n)$ such that $L_{u_{\lambda}} (\varphi) = h$
  in $\mathbb{R}^n$, concluding the proof of this lemma.
\end{proof}

\begin{lemma}
  \label{stem}Under the assumptions and notation of Proposition
  \ref{constbranch}, for any small $a > 0$ (depending on $\lambda,
  u_{\lambda}$), there exists $K (a) > 0$ such that for $h \in
  L^2_{\tmop{rad}} (\mathbb{R}^n)$ that also satisfies
  \[ \| h \|_{C^0 (a)} < + \infty, \]
  the function $\varphi \in H^2_{\tmop{rad}} (\mathbb{R}^n)$, solution of
  $L_{u_{\lambda}} (\varphi) = h$, satisfies
  \[ \| \varphi \|_{C^1 (a)} \leqslant K (a) \| h \|_{C^0 (a)} . \]
\end{lemma}

This result implies in particular that the operator $L_{u_{\lambda}}$ has a
bounded inverse
\[ L_{u_{\lambda}}^{- 1} : C^0 (a) \rightarrow C^1 (a) \]
\begin{proof}
  First, by Lemma \ref{cerasus} we have
  \[ \| \varphi \|_{H^1 (\mathbb{R}^n)} \leqslant K \| h \|_{L^2
     (\mathbb{R}^n)} \leqslant K (a) \| h \|_{C^0 (a)} . \]
  From this and by standard elliptic estimates on the equation
  $L_{u_{\lambda}} (\varphi) = h$, we have that for any $d > 0$, there exists
  $K (a, d) > 0$ such that
  \begin{equation}
    \| \varphi \|_{L^{\infty} (B (0, d))} + \| \nabla \varphi \|_{L^{\infty}
    (B (0, d))} \leqslant K (a, d) \| h \|_{C^0 (a)} . \label{turnonthelight}
  \end{equation}
  Now, we write the equation $L_{u_{\lambda}} (\varphi) = h$ as
  \[ - \Delta \varphi - \lambda \varphi = - R (\varphi) + h \]
  where
  \[ R (\varphi) = V \varphi - f_t (x, u_{\lambda}) (\varphi) . \]
  Therefore, by Lemma \ref{simon} we have
  \[ \varphi = I_{\lambda} \ast (R (\varphi) + h) . \]
  By Lemma \ref{doti} we have
  \[ \| I_{\lambda} \ast h \|_{C^1 (a)} \leqslant K (a) \| h \|_{C^0 (a)} \]
  and
  \begin{eqnarray*}
    &  & \| I_{\lambda} \ast R (\varphi) \|_{C^1 (a)}\\
    & \leqslant & K (a) \| R (\varphi) \|_{C^0 (a)}\\
    & \leqslant & K (a, d) \| \varphi \|_{C^1 (B (0, d))} + K (a) \| R
    (\varphi) e^{a | . |} \|_{L^{\infty} (\mathbb{R}^n \backslash B (0, d))} .
  \end{eqnarray*}
  By (\ref{turnonthelight}) we have $\| \varphi \|_{C^1 (B (0, d))} \leqslant
  K (d) \| h \|_{C^0 (a)}$ and by Lemma \ref{bb},
  \[ \| R (\varphi) e^{a | . |} \|_{L^{\infty} (\mathbb{R}^n \backslash B (0,
     d))} = o_{d \rightarrow + \infty} (1) \| \varphi \|_{C^1 (a)} . \]
  Combining these estimates, we have
  \[ \| \varphi \|_{C^1 (a)} \leqslant K (a, d) \| h \|_{C^0 (a)} + o_{d
     \rightarrow + \infty} (1) \| \varphi \|_{C^1 (a)}, \]
  hence taking $d$ large enough, we conclude the proof of this lemma.
\end{proof}

\subsubsection{Construction by a fixed point argument}

Now we look for $U_{\lambda + \varepsilon}$ solution of $\tmop{Eq}_{\lambda +
\varepsilon} (U_{\lambda + \varepsilon}) = 0$ under the decomposition
$U_{\lambda + \varepsilon} = u_{\lambda} + \varphi_{\varepsilon}$, and thus $\varphi_{\varepsilon}$
satisfies
\[ L_{u_{\lambda}} (\varphi_{\varepsilon}) +\mathcal{R}_{u_{\lambda}} (\varphi_{\varepsilon}) =
   \varepsilon u_{\lambda}, \]
where
\[ L_{u_{\lambda}} (\varphi) \assign - \Delta \varphi + V \varphi - \lambda
   \varphi - f_t (x, u_{\lambda}) (\varphi) \]
is the operator studied in the previous subsection, and
\[ \mathcal{R}_{u_{\lambda}} (\varphi) \assign - \varepsilon \varphi - f (x,
   u_{\lambda} + \varphi) + f (x, u_{\lambda}) + f_t (x, u_{\lambda})
   (\varphi) \]
contains the nonlinear terms in $\varphi$.

\begin{lemma}
  \label{csyst}Under the assumptions and notation of Proposition
  \ref{constbranch}, if $\| \varphi \|_{C^1 (a)} \leqslant 1$, then for $a >
  0$ small enough we have
  \[ \| \mathcal{R}_{u_{\lambda}} (\varphi) \|_{C^0 (a)} \leqslant K (a) \|
     \varphi \|_{C^0 (a)} (| \varepsilon | + \| \varphi \|_{C^0 (a)}) . \]
\end{lemma}

\begin{proof}
  This follows from ($F_1$) for the term $- f (x, u_{\lambda} + \varphi) + f
  (x, u_{\lambda}) + f_t (x, u_{\lambda}) (\varphi)$ and immediate estimates
  for the term $- \varepsilon \varphi$.
\end{proof}

\begin{lemma}
  \label{wannafly}Under the assumptions and notation of Proposition
  \ref{constbranch}, there exists $K > 0$ such that, for $| \varepsilon |$
  small enough, there exists $\varphi_{\varepsilon} \in H^1_{\tmop{rad}} \cap
  C^1$ solution of
  \[ L_{u_{\lambda}} (\varphi_{\varepsilon}) +\mathcal{R}_{u_{\lambda}}
     (\varphi_{\varepsilon}) = \varepsilon u_{\lambda} \]
  with
  \[ \| \varphi_{\varepsilon} \|_{C^1 (a)} \leqslant K | \varepsilon | . \]
\end{lemma}

\begin{proof}
  This follows from the fixed point theorem: in the space of functions
  $\varphi$ such that $\| \varphi \|_{C^1 (a)} \leqslant K \varepsilon$ for
  $\varepsilon$ small enough and some fixed $K > 0$, the function
  \[ \varphi \rightarrow L_{u_{\lambda}}^{- 1} (\mathcal{R}_{u_{\lambda}}
     (\varphi) - \varepsilon u_{\lambda}) \]
  is a contraction. This is indeed a consequence of Lemmas \ref{stem} and
  \ref{csyst}, and the fact that
  \[ \| \varepsilon u_{\lambda} \|_{C^1 (a)} \leqslant K | \varepsilon | \]
  if $a > 0$ is small enough from Lemma \ref{bb}.
\end{proof}

This concludes the construction of $U_{\lambda + \varepsilon} = u_{\lambda} +
\varphi_{\varepsilon}$ for $| \varepsilon |$ small. Let us now show some
properties of this function.

\subsubsection{Properties of $U_{\lambda + \varepsilon}$}

First, we show that $U_{\lambda + \varepsilon} = u_{\lambda} +
\varphi_{\varepsilon}$ is in $C^2 (\mathbb{R}^n, \mathbb{R})$ and decay
exponentially fast. This follows from Lemma \ref{bb} and the following result.

\begin{lemma}
  \label{hershey}The function $\varphi_{\varepsilon}$ of Lemma \ref{wannafly},
  solution of $L_{u_{\lambda}} (\varphi_{\varepsilon})
  +\mathcal{R}_{u_{\lambda}} (\varphi_{\varepsilon}) = \varepsilon
  u_{\lambda}$, is in $C^2 (\mathbb{R}^n, \mathbb{R})$ and satisfies for $|
  \varepsilon |$ small enough that
  \[ \| \varphi_{\varepsilon} \|_{C^2 (a)} \leqslant K | \varepsilon | \]
  for some constant $K, a > 0$ depending on $u_{\lambda}$ and $\lambda$.
\end{lemma}

\begin{proof}
  By $(F 1)$, Lemma \ref{bb} and standard elliptic regularity on
  $L_{u_{\lambda}} (\varphi_{\varepsilon}) +\mathcal{R}_{u_{\lambda}}
  (\varphi_{\varepsilon}) = \varepsilon u_{\lambda}$, we deduce that
  $\varphi_{\varepsilon} \in C^2 (\mathbb{R}^n, \mathbb{R})$. The decay
  estimate follows from the equation $L_{u_{\lambda}} (\varphi_{\varepsilon})
  +\mathcal{R}_{u_{\lambda}} (\varphi_{\varepsilon}) = \varepsilon
  u_{\lambda}$, with a proof similar as the one of Lemma \ref{stem}.
\end{proof}

Next, we show that $\varepsilon \rightarrow U_{\lambda + \varepsilon}$ is
$C^1$.

\begin{lemma}
  \label{gosepe}There exists $\varepsilon_0 > 0$ such that, the function
  $\varphi_{\varepsilon}$ of Lemma \ref{wannafly} satisfies that $\varepsilon
  \rightarrow \varphi_{\varepsilon} \in C^1 (] - \varepsilon_0, \varepsilon_0
  [, C^2 (a))$, with
  \[ \| \partial_{\varepsilon} \varphi_{\varepsilon} \|_{C^2 (a)} \leqslant K
  \]
  for some $K, a > 0$. In particular, $\varepsilon \rightarrow E (U_{\lambda +
  \varepsilon})$ and $\varepsilon \rightarrow Q (U_{\lambda + \varepsilon})$
  are $C^1$ functions on $] - \varepsilon_0, \varepsilon_0 [$.
\end{lemma}

\begin{proof}
  First, we infer that for $| \varepsilon |$ small enough, the operator
  $L_{u_{\lambda}} + d_{\varphi} \mathcal{R} (\varphi_{\varepsilon})$ is
  invertible. This is a consequence of the fact that $L_{u_{\lambda}}$ is
  invertible (uniformly in $\varepsilon$) by Lemma \ref{stem}, and that by
  Lemma \ref{bb},
  \[ \| d_{\varphi} \mathcal{R} (\varphi_{\varepsilon}) (\varphi) \|_{C^2 (a)}
     \leqslant K \varepsilon \| \varphi \|_{C^2 (a)} \]
  since $\| \varphi_{\varepsilon} \|_{C^2 (a)} \leqslant K | \varepsilon |$ by
  Lemma \ref{hershey}. Therefore, by the implicit function theorem,
  $\varepsilon \rightarrow \varphi_{\varepsilon}$ is $C^1$ for $| \varepsilon
  |$ small enough with
  \begin{equation}
    \partial_{\varepsilon} \varphi_{\varepsilon} = (L_{u_{\lambda}} +
    d_{\varphi} \mathcal{R} (\varphi_{\varepsilon}))^{- 1} (u_{\lambda} +
    \varphi_{\varepsilon}) . \label{kutaisi}
  \end{equation}
  Indeed, we have
  \[ L_{u_{\lambda}} (\varphi_{\varepsilon}) +\mathcal{R}_{u_{\lambda}}
     (\varphi_{\varepsilon}) = \varepsilon u_{\lambda} \]
  and differentiating with respect to $\varepsilon$ leads to (\ref{kutaisi})
  (recall that $\mathcal{R}_{u_{\lambda}} (\varphi_{\varepsilon})$ contains a
  term $- \varepsilon \varphi_{\varepsilon}$). In particular,
  \[ \| \partial_{\varepsilon} \varphi_{\varepsilon} \|_{C^2 (a)} \leqslant K
     \| u_{\lambda} + \varphi_{\varepsilon} \|_{C^2 (a)} \leqslant K. \]
  Now, recall that
  \[ E (U_{\lambda + \varepsilon}) = \frac{1}{2} \int_{\mathbb{R}^n} | \nabla
     U_{\lambda + \varepsilon} |^2 + V (| x |) | U_{\lambda + \varepsilon} |^2
     -\mathfrak{R}\mathfrak{e} (f (| x |, U_{\lambda + \varepsilon})
     \overline{U_{\lambda + \varepsilon}}) . \]
  From $U_{\lambda + \varepsilon} = u_{\lambda} + \varphi_{\varepsilon}$ and
  $\| \partial_{\varepsilon} \varphi_{\varepsilon} \|_{C^2 (a)} \leqslant K$
  as well as the fact that $f$ is $C^1$ with respect to its second variable,
  we check easily that $E (U_{\lambda + \varepsilon})$ is $C^1$ with respect
  to $\varepsilon$ for $| \varepsilon |$ small enough and the same hold for $Q
  (U_{\lambda + \varepsilon})$.
\end{proof}

\begin{remark}
	The proof of Lemma \ref{gosepe} can be adapted if $f$ is analytic to show
	that $\varepsilon \rightarrow \varphi_{\varepsilon}$ and $\varepsilon
	\rightarrow E (U_{\lambda + \varepsilon}), \varepsilon \rightarrow Q
	(U_{\lambda + \varepsilon})$ are analytic functions for $| \varepsilon |$
	small.
\end{remark}

We conclude this subsection with some identities involving $U_{\lambda +
	\varepsilon}$.

\begin{lemma}
	\label{massnotzero}Under the previous notation and hypotheses, we have the
	following identities:

	\[ \partial_{\varepsilon} (Q (U_{\lambda + \varepsilon}))_{| \varepsilon =
		0 } = \langle L^{- 1}_{u_{\lambda}} (u_{\lambda}), u_{\lambda}
	\rangle, \]
	\[ L_{U_{\lambda + \varepsilon}} (\partial_{\varepsilon} U_{\lambda +
		\varepsilon}) = U_{\lambda + \varepsilon} \]
	and
	\[ \partial_{\varepsilon} (E (U_{\lambda + \varepsilon}) - \lambda Q
	(U_{\lambda + \varepsilon}))_{| \varepsilon = 0 } = 0. \]
\end{lemma}

The second identity shows that $\partial_{\varepsilon} U_{\lambda +
	\varepsilon | \varepsilon = 0 } = L^{- 1}_{u_{\lambda}}
(u_{\lambda})$, which is a good sign that the branch contruscted above plays
an important role in understanding $u_{\lambda}$.

\

\begin{proof}
	We recall that
	\[ Q (U_{\lambda + \varepsilon}) = \frac{1}{2} \int (U_{\lambda +
		\varepsilon})^2 \]
	and thus
	\[ \partial_{\varepsilon} (Q (U_{\lambda + \varepsilon}))_{|
		\varepsilon = 0} = \int \partial_{\varepsilon} U_{\lambda + \varepsilon |
		\varepsilon = 0 } u_{\lambda} . \]
	Now, simply differentiating the equation $\tmop{Eq}_{\lambda + \varepsilon}
	(U_{\lambda + \varepsilon}) = 0$ with respect to $\varepsilon$ shows the
	identity $L_{U_{\lambda + \varepsilon}} (\partial_{\varepsilon} U_{\lambda +
		\varepsilon}) = U_{\lambda + \varepsilon}$ and thus $\partial_{\varepsilon}
	U_{\lambda + \varepsilon | \varepsilon = 0 } = L^{-
		1}_{u_{\lambda}} (u_{\lambda})$, concluding the proof of
	\[ \partial_{\varepsilon} (Q (U_{\lambda + \varepsilon}))_{| \varepsilon = 0
	} = \langle L^{- 1}_{u_{\lambda}} (u_{\lambda}), u_{\lambda}
	\rangle . \]
	Finally, we check by integration by parts that
	\[ \partial_{\varepsilon} (E (U_{\lambda + \varepsilon}) - \lambda Q
	(U_{\lambda + \varepsilon}))_{| \varepsilon = 0 } = \langle
	\tmop{Eq}_{\lambda} (u_{\lambda}), \partial_{\varepsilon} U_{\lambda +
		\varepsilon | \varepsilon = 0 } \rangle = 0. \]
\end{proof}

\subsection{Local uniqueness near the branch}

In this subsection, we show Propositions \ref{coerc2} and \ref{uniquresult}. The identities and the end of Proposition \ref{coerc2} have been shown in
Lemma \ref{massnotzero} above.

\subsubsection{Proof of Proposition \ref{coerc2}}

\begin{proof}
	We suppose here first that for any $\varphi \in H_{\tmop{rad}}^1,$ if
	$\langle \varphi, u_{\lambda} \rangle = \langle \varphi, L^{-
		1}_{u_{\lambda}} (u_{\lambda}) \rangle = 0$ then
	\[ \langle L_{u_{\lambda}} (\varphi), \varphi \rangle \geqslant \kappa \|
	\varphi \|_{H^1}^2 . \]
	Therefore, by standard coercivity methods, without orthogonality conditions
	and up to reducing the value of $\kappa$, we have that
	\[ \langle L_{u_{\lambda}} (\varphi), \varphi \rangle \geqslant \kappa \|
	\varphi \|_{H^1}^2 - \frac{1}{\kappa} \langle \varphi, u_{\lambda}
	\rangle^2 - \frac{1}{\kappa} \langle \varphi, L^{- 1}_{u_{\lambda}}
	(u_{\lambda}) \rangle^2 . \]
	Now, remark that
	\[ L_{U_{\lambda + \varepsilon}} (\varphi) - L_{u_{\lambda}} (\varphi) = -
	\varepsilon \varphi - f_t (x, U_{\lambda + \varepsilon}) (\varphi) + f_t
	(x, u_{\lambda}) (\varphi) \]
	hence, by Lemmas \ref{bb} and \ref{hershey},
	\[ \langle L_{U_{\lambda + \varepsilon}} (\varphi) - L_{u_{\lambda}}
	(\varphi), \varphi \rangle = O_{\varepsilon \rightarrow 0} (\varepsilon)
	\| \varphi \|_{H^1}^2 . \]
	Similarly,
	\[ \langle \varphi, u_{\lambda} \rangle^2 - \langle \varphi, U_{\lambda +
		\varepsilon} \rangle^2 = O_{\varepsilon \rightarrow 0} (\varepsilon) \|
	\varphi \|_{H^1}^2 . \]
	and
	\[ \langle \varphi, L^{- 1}_{u_{\lambda}} (u_{\lambda}) \rangle^2 - \langle
	\varphi, L^{- 1}_{U_{\lambda + \varepsilon}} (U_{\lambda + \varepsilon})
	\rangle^2 = O_{\varepsilon \rightarrow 0} (\varepsilon) \| \varphi
	\|_{H^1}^2 . \]
	We deduce that, for $\varepsilon$ small enough,
	\begin{eqnarray*}
		\langle L_{U_{\lambda + \varepsilon}} (\varphi), \varphi \rangle &
		\geqslant & \kappa \| \varphi \|_{H^1}^2 - \frac{1}{\kappa} \langle
		\varphi, U_{\lambda + \varepsilon} \rangle^2 - \frac{1}{\kappa} - \langle
		\varphi, L^{- 1}_{U_{\lambda + \varepsilon}} (U_{\lambda + \varepsilon})
		\rangle^2 - O_{\varepsilon \rightarrow 0} (\varepsilon) \| \varphi
		\|_{H^1}^2\\
		& \geqslant & \frac{\kappa}{2} \| \varphi \|_{H^1}^2 - \frac{1}{\kappa}
		\langle \varphi, U_{\lambda + \varepsilon} \rangle^2 - \frac{1}{\kappa} -
		\langle \varphi, L^{- 1}_{U_{\lambda + \varepsilon}} (U_{\lambda +
			\varepsilon}) \rangle^2 .
	\end{eqnarray*}
	If we suppose instead that for any $\varphi \in H_{\tmop{rad}}^1,$ if
	$\langle \varphi, u_{\lambda} \rangle = 0$ implies $\langle L_{u_{\lambda}}
	(\varphi), \varphi \rangle \geqslant \kappa \| \varphi \|_{H^1}^2$, we can
	do a similar proof. This concludes the proof of Proposition \ref{coerc2}.
\end{proof}

\subsubsection{Proof of Proposition \ref{uniquresult}}

\begin{proof}
	We will split the proof into two parts, depending on which coercivity result we
	use. But first, remark that
	\[ | Q (\psi) - Q (u_{\lambda}) | \leqslant \| \psi - u_{\lambda} \|_{L^2}
	\| \psi + u_{\lambda} \|_{L^2} \leqslant K \eta . \]
	Then, let us show that $| \lambda - \lambda_0 | = o_{\eta \rightarrow 0}
	(1)$. Since $\tmop{Eq}_{\lambda_0} (\psi) - \tmop{Eq}_{\lambda}
	(u_{\lambda}) = 0$, taking its scalar product with $\psi$ (and taking $\eta$
	small enough) leads to
	\[ (\lambda - \lambda_0) \| \psi \|_{L^2} \leqslant K \| \psi - u_{\lambda}
	\|_{H^1} \leqslant K \eta . \]
	These estimates hold in both cases below.
	
	\
	
	{\emph{Case 1: assuming $(E 2)$ and $\langle u_{\lambda}, L^{-
				1}_{u_{\lambda}} (u_{\lambda}) \rangle \neq 0$ hold.}}
	
	\
	
	We decompose the function $\psi$ for some small $| \varepsilon |$ as
	\[ \psi = U_{\lambda + \varepsilon} + \varphi . \]
	By Lemma \ref{hershey} and the hypotheses of Proposition \ref{uniquresult},
	\begin{equation}
		\| \varphi \|_{H^1} \leqslant \| \psi - u_{\lambda} \|_{H^1} + \|
		u_{\lambda} - U_{\lambda + \varepsilon} \| \leqslant K (\eta + |
		\varepsilon |) . \label{goodpeople}
	\end{equation}
	Furthermore,
	\[ Q (\psi) = Q (U_{\lambda + \varepsilon}) + \langle \varphi, U_{\lambda +
		\varepsilon} \rangle + \| \varphi \|_{L^2}^2 . \]
	Now, we recall from Lemma \ref{massnotzero} and $\langle u_{\lambda}, L^{-
		1}_{u_{\lambda}} (u_{\lambda}) \rangle \neq 0$ that $\partial_{\varepsilon}
	Q (U_{\lambda + \varepsilon})_{| \varepsilon = 0 } \neq 0$ and
	also that
	\[ | Q (\psi) - Q (u_{\lambda}) | \leqslant K \eta . \]
	Therefore, if $\eta$ is small enough, we can choose $\varepsilon$ small such
	that $Q (U_{\lambda + \varepsilon}) = Q (\psi)$, since $Q (U_{\lambda +
		\varepsilon}) = Q (u_{\lambda}) + \varepsilon \partial_{\varepsilon} Q
	(U_{\lambda + \varepsilon})_{| \varepsilon = 0 } + O_{\varepsilon
		\rightarrow 0} (\varepsilon^2) .$ Then, by (\ref{goodpeople}) we deduce that
	\begin{equation}
		| \langle \varphi, U_{\lambda + \varepsilon} \rangle | = \| \varphi
		\|_{L^2}^2 \leqslant \| \varphi \|_{H^1}^2 \leqslant K (\eta + |
		\varepsilon |) \| \varphi \|_{H^1} . \label{sunny}
	\end{equation}
	Now, since $\tmop{Eq}_{\lambda_0} (\psi) = 0$, we have by the
	decomposition $\psi = U_{\lambda + \varepsilon} + \varphi$ that
	\begin{equation}
		L_{U_{\lambda + \varepsilon}} (\varphi) + (\lambda + \varepsilon -
		\lambda_0) U_{\lambda + \varepsilon} + (\lambda + \varepsilon - \lambda_0)
		\varphi + \tmop{NL}_{U_{\lambda + \varepsilon}} (\varphi) = 0
		\label{nosilence}
	\end{equation}
	where
	\[ \tmop{NL}_{U_{\lambda + \varepsilon}} (\varphi) = - f (x, U_{\lambda +
		\varepsilon} + \varphi) + f (x, U_{\lambda + \varepsilon}) + f_t (x,
	U_{\lambda + \varepsilon}) (\varphi) . \]
	Taking the scalar product of (\ref{nosilence}) with $\partial_{\varepsilon}
	U_{\lambda + \varepsilon}$, we have
	\begin{eqnarray*}
		&  & (\lambda + \varepsilon - \lambda_0) \langle \partial_{\varepsilon}
		U_{\lambda + \varepsilon}, U_{\lambda + \varepsilon} \rangle\\
		& = & \langle \partial_{\varepsilon} U_{\lambda + \varepsilon},
		L_{U_{\lambda + \varepsilon}} (\varphi) \rangle + (\lambda + \varepsilon -
		\lambda_0) \langle \varphi, \partial_{\varepsilon} U_{\lambda +
			\varepsilon} \rangle\\
		& + & \langle L_{U_{\lambda + \varepsilon}} (\varphi),
		\partial_{\varepsilon} U_{\lambda + \varepsilon} \rangle .
	\end{eqnarray*}
	We recall that
	\[ L_{U_{\lambda + \varepsilon}} (\partial_{\varepsilon} U_{\lambda +
		\varepsilon}) = U_{\lambda + \varepsilon}, \]
	in particular by Lemma \ref{gosepe}, there exists a constant $K > 0$
	independent of $\varepsilon$ such that
	\[ \| \partial_{\varepsilon} U_{\lambda + \varepsilon} \|_{H^1} \leqslant K
	\]
	for $\varepsilon$ small enough. Also, by Lemma \ref{massnotzero}, $| \langle
	\partial_{\varepsilon} U_{\lambda + \varepsilon}, U_{\lambda + \varepsilon}
	\rangle | = | \partial_{\varepsilon} (Q (U_{\lambda + \varepsilon})) |
	\geqslant K > 0$ where $K$ is independent of $\varepsilon$, provided that $|
	\varepsilon |$ is small enough.
	
	\
	
	Now, by (\ref{sunny}) we have
	\[ | \langle \partial_{\varepsilon} U_{\lambda + \varepsilon}, L_{U_{\lambda
			+ \varepsilon}} (\varphi) \rangle | = | \langle L_{U_{\lambda +
			\varepsilon}} (\partial_{\varepsilon} U_{\lambda + \varepsilon}), \varphi
	\rangle | = | \langle U_{\lambda + \varepsilon}, \varphi \rangle |
	\leqslant K (\eta + | \varepsilon |) \| \varphi \|_{H^1} . \]
	Additionally,
	\[ | (\lambda + \varepsilon - \lambda_0) \langle \varphi,
	\partial_{\varepsilon} U_{\lambda + \varepsilon} \rangle | \leqslant |
	\lambda + \varepsilon - \lambda_0 | \| \varphi \|_{L^2} \|
	\partial_{\varepsilon} U_{\lambda + \varepsilon} \|_{L^2} \leqslant K
	(\eta + | \varepsilon |) \| \varphi \|_{H^1} \]
	and by Lemmas \ref{bb} and \ref{hershey} and (\ref{goodpeople})
	\[ | \langle L_{U_{\lambda + \varepsilon}} (\varphi), \partial_{\varepsilon}
	U_{\lambda + \varepsilon} \rangle | \leqslant K \| \varphi \|_{H^1}^2
	\leqslant K (\eta + | \varepsilon |) \| \varphi \|_{H^1} . \]
	Combining these estimates, we deduce that
	\begin{equation}
		| \lambda + \varepsilon - \lambda_0 | \leqslant K (\eta + | \varepsilon |)
		\| \varphi \|_{H^1} \label{genixremix}
	\end{equation}
	for some constant $K$ independent of $\varepsilon, \eta$.
	
	\
	
	Now, taking the scalar product of (\ref{nosilence}) with $\varphi$, we have
	\begin{eqnarray*}
		0 & = & \langle L_{U_{\lambda + \varepsilon}} (\varphi), \varphi \rangle +
		(\lambda_0 - \lambda + \varepsilon) \langle U_{\lambda + \varepsilon}, \varphi \rangle\\
		& + & (\lambda + \varepsilon - \lambda_0) \langle \varphi, \varphi
		\rangle + \langle \tmop{NL}_{U_{\lambda + \varepsilon}} (\varphi), \varphi
		\rangle .
	\end{eqnarray*}
	By Proposition \ref{coerc2} and (\ref{sunny}), we have
	\[ \langle L_{U_{\lambda + \varepsilon}} (\varphi), \varphi \rangle
	\geqslant \kappa \| \varphi \|_{H^1}^2 - K (\eta + | \varepsilon |)^2 \|
	\varphi \|_{H^1}^2 . \]
	By (\ref{genixremix}) we can write
	\[ | (\lambda + \varepsilon - \lambda_0) \langle U_{\lambda + \varepsilon},
	\varphi \rangle | \leqslant K (\eta + | \varepsilon |) \| \varphi
	\|_{H^1}^2 \]
	and by $(F 1)$, Lemmas \ref{bb} and \ref{hershey} as well as
	(\ref{goodpeople}), taking $\eta$ small enough,
	\[ | \langle \tmop{NL}_{U_{\lambda + \varepsilon}} (\varphi), \varphi
	\rangle | \leqslant K (\eta + | \varepsilon |) \| \varphi \|_{H^1}^2 . \]
	Combining these estimates together, we deduce that
	\[ (\kappa - K (\eta + | \varepsilon |)) \| \varphi \|_{H^1}^2 \leqslant 0,
	\]
	therefore, for $\eta$ small enough (which implies that $| \varepsilon |$ is
	small as well), we have $\| \varphi \|_{H^1} = 0$. This, together with (\ref{genixremix})
   implies that $\lambda_0 = \lambda + \varepsilon$, and this conclude the
	proof of the proposition.
	
	\
	
	{\emph{Case 2: assuming $(E 1)$ hold and $\langle u_{\lambda}, L^{-
				1}_{u_{\lambda}} (u_{\lambda}) \rangle = 0$.}}
	
	\
	
	In this case, we decompose similarly
	\[ \psi = U_{\lambda + \varepsilon} + \varphi . \]
	However now, we choose $\varepsilon$ small such that $\langle \varphi,
	\partial_{\lambda} U_{\lambda + \varepsilon} \rangle = 0$. Let us show that
	this is possible. We have $| \langle \varphi, \partial_{\lambda} U_{\lambda
		+ \varepsilon | \varepsilon = 0 } \rangle | \leqslant K \| \varphi
	\|_{L^2} \leqslant K \eta$ and since $\varphi = \psi - U_{\lambda +
		\varepsilon}$,
	\[ \partial_{\varepsilon} (\langle \varphi, \partial_{\lambda} U_{\lambda +
		\varepsilon} \rangle)_{| \varepsilon = 0 } = - \langle
	\partial_{\lambda} U_{\lambda}, \partial_{\lambda} U_{\lambda} \rangle +
	\langle \varphi, \partial_{\lambda}^2 U_{\lambda + \varepsilon |
		\varepsilon = 0 } \rangle . \]
	We have $\langle \partial_{\lambda} U_{\lambda}, \partial_{\lambda}
	U_{\lambda} \rangle \geqslant K > 0$ while $| \langle \varphi,
	\partial_{\lambda}^2 U_{\lambda + \varepsilon | \varepsilon = 0 }
	\rangle | \leqslant K \eta$, therefore for $\eta$ small enough we can indeed
	choose $| \varepsilon |$ small (with $| \varepsilon | = o_{\eta \rightarrow
		0} (1)$) such that $\langle \varphi, \partial_{\lambda} U_{\lambda +
		\varepsilon} \rangle = 0$.
	
	Differentiating $\tmop{Eq}_{\lambda + \varepsilon} (U_{\lambda +
		\varepsilon}) = 0$ with respect to $\varepsilon$ leads to $L_{U_{\lambda +
			\varepsilon}} (\partial_{\lambda} U_{\lambda + \varepsilon}) = U_{\lambda +
		\varepsilon}$, and therefore with this choice of $\varepsilon$ we just
	showed that
	\[ \langle \varphi, L_{U_{\lambda + \varepsilon}}^{- 1} (U_{\lambda +
		\varepsilon}) \rangle = 0. \]

	We now study the quantity $\langle \varphi, U_{\lambda + \varepsilon}
	\rangle$, which is the second orthogonality needed in this case. We recall
	from Lemma \ref{massnotzero} that $L_{U_{\lambda + \varepsilon}}
	(\partial_{\varepsilon} U_{\lambda + \varepsilon}) = U_{\lambda +
		\varepsilon}$, therefore
	\[ \langle \varphi, U_{\lambda + \varepsilon} \rangle = \langle \varphi,
	L_{U_{\lambda + \varepsilon}} (\partial_{\varepsilon} U_{\lambda +
		\varepsilon}) \rangle = \langle L_{U_{\lambda + \varepsilon}} (\varphi),
	\partial_{\varepsilon} U_{\lambda + \varepsilon} \rangle . \]
	Since $\psi = U_{\lambda + \varepsilon} + \varphi$ and
	$\tmop{Eq}_{\lambda_0} (\psi) = 0$, we have
	\begin{equation}
		0 = \tmop{Eq}_{\lambda_0} (U_{\lambda + \varepsilon}) + L_{U_{\lambda +
				\varepsilon}} (\varphi) + (\lambda + \varepsilon - \lambda_0) \varphi +
		\tmop{NL}_{U_{\lambda + \varepsilon}} (\varphi) \label{passion}
	\end{equation}
	and since $\tmop{Eq}_{\lambda + \varepsilon} (U_{\lambda + \varepsilon}) =
	0$, we have
	\[ \tmop{Eq}_{\lambda_0} (U_{\lambda + \varepsilon}) = (\lambda +
	\varepsilon - \lambda_0) (U_{\lambda + \varepsilon}) . \]
	Taking the scalar product of (\ref{passion}) with $U_{\lambda +
		\varepsilon}$ leads to
	\begin{eqnarray*}
		&  & (\lambda + \varepsilon - \lambda_0) \| U_{\lambda + \varepsilon}
		\|_{L^2}\\
		& + & \langle \varphi, L_{U_{\lambda + \varepsilon}} (U_{\lambda +
			\varepsilon}) \rangle + (\lambda + \varepsilon - \lambda_0) \langle
		\varphi, U_{\lambda + \varepsilon} \rangle + \langle \tmop{NL} (\varphi),
		U_{\lambda + \varepsilon} \rangle\\
		& = & 0
	\end{eqnarray*}
	hence by the estimates on $U_{\lambda + \varepsilon}$ of Lemma
	\ref{hershey}, we deduce that this implies that
	\[ | \lambda + \varepsilon - \lambda_0 | \leqslant K \| \varphi \|_{L^2}
	\leqslant K \eta . \]
	Now, taking the scalar product of (\ref{passion}) with
	$\partial_{\varepsilon} U_{\lambda + \varepsilon}$, we have
	\begin{eqnarray*}
		\langle \varphi, U_{\lambda + \varepsilon} \rangle & = & \langle
		L_{U_{\lambda + \varepsilon}} (\varphi), \partial_{\varepsilon} U_{\lambda
			+ \varepsilon} \rangle\\
		& = & - (\lambda + \varepsilon - \lambda_0) \langle U_{\lambda +
			\varepsilon}, \partial_{\varepsilon} U_{\lambda + \varepsilon} \rangle\\
		& - & (\lambda + \varepsilon - \lambda_0) \langle \varphi,
		\partial_{\varepsilon} U_{\lambda + \varepsilon} \rangle - \langle
		\tmop{NL} (\varphi), \partial_{\varepsilon} U_{\lambda + \varepsilon}
		\rangle .
	\end{eqnarray*}
	We have
	\[ | (\lambda + \varepsilon - \lambda_0) \langle \varphi,
	\partial_{\varepsilon} U_{\lambda + \varepsilon} \rangle | \leqslant K
	\eta \| \varphi \|_{L^2} \]
	and as in the first case,
	\[ | \langle \tmop{NL} (\varphi), \partial_{\varepsilon} U_{\lambda +
		\varepsilon} \rangle | \leqslant K \eta \| \varphi \|_{L^2} . \]
	Finally, we recall that we supposed that $\langle U_{\lambda + \varepsilon},
	\partial_{\varepsilon} U_{\lambda + \varepsilon} \rangle_{| \varepsilon = 0
	} = \langle u_{\lambda}, L^{- 1}_{u_{\lambda}} (u_{\lambda})
	\rangle = 0.$ Since $\varepsilon \rightarrow U_{\lambda + \varepsilon}$ is
	$C^1$ by Lemma \ref{gosepe}, we deduce that
	\[ \langle U_{\lambda + \varepsilon}, \partial_{\varepsilon} U_{\lambda +
		\varepsilon} \rangle = o_{\varepsilon \rightarrow 0} (1). \]
	Combining this $| \lambda + \varepsilon - \lambda_0 | \leqslant K \|
	\varphi \|_{L^2}$, we deduce that
	\[ | \langle \varphi, U_{\lambda + \varepsilon} \rangle | \leqslant K
	(o_{\varepsilon \rightarrow 0} (1) + \eta) \| \varphi \|_{L^2} . \]
	We can now conclude as in the first case by taking the scalar product of
	(\ref{passion}) with $\varphi$.
\end{proof}

\subsection{Proof of Lemmas \ref{< 0}, \ref{=0 is finite}, Propositions \ref{properties of m(c)}, \ref{simo}, \ref{ext 1} and
\ref{ext 2}}

\subsubsection{Proof of Proposition \ref{properties of m(c)}} \label{sec. properties of m(c)}
\begin{proof}
	First we prove (1). Let $c_n \rightarrow c^+ \in ]c_{\ast},c^{\ast}[$, $u_n \in S_{c_n}$, $u \in S_c$. Let $\lambda_n$, $\lambda$ be the Lagrange multipliers corresponding to $u_n$, $u$ respectively. By $(A 3)$, passing to a subsequence if necessary, we may assume that $u_n \rightarrow u_{\infty}$ in $H^1$ and $\lambda_n \rightarrow \lambda_{\infty}$. Using the continuity of the functional $E$ in $H^1,$ we can conclude that $u_{\infty}$ is a NGSS on $S_c$ and $\lambda_{\infty}$ is the Lagrange multiplier. Note that $\sqrt{\frac{c}{c_n}}u_n \in X_c$. Thus
	$$
	E(u) \leq E(\sqrt{\frac{c}{c_n}}u_n),
	$$
	and then
	\begin{eqnarray}
		m(c_n) - m(c) &=& E(u_n) - E(u) \nonumber \\
		&\geq& E(u_n) - E(\sqrt{\frac{c}{c_n}}u_n) \nonumber \\
		&=& \langle D_uE(u_n),(1 - \sqrt{\frac{c}{c_n}})u_n\rangle + o(1 - \sqrt{\frac{c}{c_n}}) \nonumber \\
		&=& 2(\sqrt{c_n} - \sqrt{c})\sqrt{c_n}\lambda_n + o(\sqrt{c_n} - \sqrt{c}),
	\end{eqnarray}
	Similarly, we have
	$$
	E(u_n) \leq E(\sqrt{\frac{c_n}{c}}u),
	$$
	and then
	\begin{eqnarray}\label{2.10}
		m(c_n) - m(c) &=& E(u_n) - E(u) \nonumber \\
		&\leq& E(\sqrt{\frac{c_n}{c}}u) - E(u) \nonumber \\
		&=& \langle D_uE(u),(\sqrt{\frac{c_n}{c}}-1)u\rangle + o(\sqrt{\frac{c_n}{c}}-1) \nonumber \\
		&=& 2(\sqrt{c_n} - \sqrt{c})\sqrt{c}\lambda + o(\sqrt{c_n} - \sqrt{c}).
	\end{eqnarray}
	
	Since $c_n > c$, we have
	\begin{equation} \label{eq5.15}
		\frac{2\sqrt{c_n}}{\sqrt{c_n} + \sqrt{c}}\lambda_n + o_n(1) \leq \frac{m(c_n) - m(c)}{c_n-c} \leq \frac{2\sqrt{c}}{\sqrt{c_n} + \sqrt{c}}\lambda + o_n(1).
	\end{equation}
	We can take $\lambda = \min\mathcal{T}_c$. Then using \eqref{eq5.15} and that $\lambda_n \to \lambda_\infty \geq \min\mathcal{T}_c$, we obtain $\lambda_{\infty} = \min\mathcal{T}_c$ (this ensures that $\min \mathcal{T}_c > -\infty$ and this is well-defined). It follows that
	$$
	\lim_{n \rightarrow \infty}\frac{m(c_n) - m(c)}{c_n-c} = \min\mathcal{T}_c,
	$$
	implying that $m_+'(c) = \min\mathcal{T}_c$. Furthermore, by the arbitrariness of the choice for $\lambda_n$ one gets
	$$
	\lim_{c \to \tilde{c}^+}\min\mathcal{T}_c = \lim_{c \to \tilde{c}^+}\max\mathcal{T}_c = \min\mathcal{T}_{\tilde{c}}.
	$$
	Similarly, we can show that $m'_{-}(c) = \max\mathcal{T}_c$ and
	$$
	\lim_{c \to \tilde{c}^-}\min\mathcal{T}_c = \lim_{c \to \tilde{c}^-}\max\mathcal{T}_c = \max\mathcal{T}_{\tilde{c}}.
	$$
	This implies that $m(c)$ is Lipschitz continuous. Note that $\max \mathcal{T}_c < \Lambda_0 \leq 0$ by $(A 1)$. This yields that $m(c)$ is strictly decreasing and differentable at almost everywhere.
	
	Now we prove (2). The definition of $m(c)$ yields that for any $u_c \in S_c$ the inequality
	\begin{align}
		E(u_c+\epsilon h) \geq m(Q(u_c+\epsilon h)) = m(c + \epsilon\langle u_c,h\rangle + \epsilon^2Q(h))
	\end{align}
	is valid for all $\epsilon$ and all test functions $h$. Moreover,
	\begin{align}
		E(u_c+\epsilon h) & = E(u_c) + \epsilon\lambda_c\langle u_c,h \rangle + \frac{\epsilon^2}{2}\langle L_{u_c}h,h\rangle+ \frac{\epsilon^2}{2}\lambda_c\langle h,h\rangle + o(\epsilon^2) \nonumber \\
		& = m(c) + \epsilon\lambda_c\langle u_c,h \rangle + \frac{\epsilon^2}{2}\langle L_{u_c}h,h\rangle + \frac{\epsilon^2}{2}\lambda_c\langle h,h\rangle + o(\epsilon^2),
	\end{align}
	where $\lambda_c$ is the Lagrange multiplier for $u_c$. We obtain
	\begin{align}
		m(c + \epsilon\langle u_c,h\rangle + \epsilon^2Q(h)) \leq m(c) + \epsilon\lambda_c\langle u_c,h \rangle + \frac{\epsilon^2}{2}\langle L_{u_c}h,h\rangle + \frac{\epsilon^2}{2}\lambda_c\langle h,h\rangle + o(\epsilon^2).
	\end{align}
	Writing the same inequality with $\epsilon \to -\epsilon$ and adding the two yields
	\begin{align}
		& m(c + \epsilon\langle u_c,h\rangle + \epsilon^2Q(h)) + m(c - \epsilon\langle u_c,h\rangle + \epsilon^2Q(h)) \nonumber \\
		& \leq 2m(c) + \epsilon^2\langle L_{u_c}h,h\rangle + \epsilon^2\lambda_c\langle h,h\rangle + o(\epsilon^2).
	\end{align}
	Next we take suitable test function. $L_{u_c}$ has exactly one simple negative eigenvalue, denoted by $\sigma_c < 0$. Let $\xi_c$ be the normalized eigenfunction. Note that $\langle L_{u_c}\phi,\phi\rangle \geq 0$ for any $\langle \phi, u_c \rangle = 0$. This implies that $\langle \xi_c,u_c \rangle \neq 0$. Taking $h := \frac{\xi_c}{\langle \xi_c,u_c \rangle}$, one gets
	\begin{align}
		& m(c + \epsilon + \epsilon^2Q(h)) + m(c - \epsilon + \epsilon^2Q(h)) - 2m(c) \nonumber \\
		& \leq 2\epsilon^2\lambda_cQ(h) + 2\epsilon^2\sigma_cQ(h) + o(\epsilon^2).
	\end{align}
	Using $m'_{+}(c) = \lambda_{min}(c)$, that $m'(c)$ exists almost everywhere and that $m'_{\pm}(c)$ is bounded, we have
	\begin{align}
		m(c + \epsilon + \epsilon^2Q(h)) - m(c + \epsilon) = \int_{c + \epsilon}^{c + \epsilon + \epsilon^2Q(h)}\lambda_{min}(z)dz.
	\end{align}
	Then using the right continunity of $\lambda_{min}$ at $c$, we derive that
	\begin{align}\label{2.18}
		& m(c + \epsilon + \epsilon^2Q(h)) - m(c + \epsilon) - \lambda_{min}(c)\epsilon^2Q(h) \nonumber \\
		= & \int_{c + \epsilon}^{c + \epsilon + \epsilon^2Q(h)}(\lambda_{min}(z) - \lambda_{min}(c))dz = o(\epsilon^2).
	\end{align}
	Similarly, for $\epsilon$ small enough we have
	\begin{align}
		& m(c - \epsilon + \epsilon^2Q(h)) - m(c - \epsilon) - \lambda_{max}(c)\epsilon^2Q(h) \nonumber \\
		= & \int_{c - \epsilon}^{c - \epsilon + \epsilon^2Q(h)}(\lambda_{max}(z) - \lambda_{max}(c))dz = o(\epsilon^2).
	\end{align}
	It follows that for all $\epsilon > 0$ small enough,
	\begin{align}
		& m(c + \epsilon) + m(c - \epsilon) - 2m(c) \nonumber \\
		\leq &  2\epsilon^2\lambda_cQ(h) - (\lambda_{min}(c) + \lambda_{max}(c))\epsilon^2Q(h) + 2\epsilon^2\sigma_cQ(h) + o(\epsilon^2).
	\end{align}
	Taking $\lambda_c = \lambda_{min}(c)$ we have
	\begin{align} \label{eq concave}
		\limsup_{\epsilon \to 0}\frac{m(c + \epsilon) + m(c - \epsilon) - 2m(c)}{\epsilon^2}  \leq 2\sigma_cQ(h) < 0.
	\end{align}
	Using \eqref{eq concave} we conclude that the function $m$ is strictly concave down in a neighborhood of $c$. This in turns shows that $m$ is strictly concave down on all $c \in ]c_{\ast},c^{\ast}[$. Thus $m(c)$ is twice differentiable at almost every $c \in ]c_{\ast},c^{\ast}[$ and $m''(c) < 0$ when it exists. Next we show that $m'_-(c_1) < m'_+(c_2), \forall c_1 > c_2$. Note that $\mathcal{T}_c$ has just one element for almost every $c \in ]c_{\ast},c^{\ast}[$, denoted by $\lambda(c)$. And $\lambda'(c) = m''(c) < 0$ for almost every $c \in ]c_{\ast},c^{\ast}[$. Hence, both $\lambda_{min}(c)$ and $\lambda_{max}(c)$ are strictly decreasing. Then using (1), for any $c_1 > c_2$,
	\begin{align}
		\lambda_{min}(c_2) = \lim_{c \to c_2^+}\lambda_{max}(c) > \lambda_{max}(c_1).
	\end{align}
	It remains to prove the boundedness of Dini derivative. We use \eqref{eq right assum} to prove \eqref{eq right res}. The key idea is to apply the following formulas 
	\begin{align} \label{eq 1}
		2\lambda_{min}(c)c = \int_{\mathbb{R}^n}|\nabla u_c|^2-f(x,u_c)u_c,
	\end{align}
	and
	\begin{align} \label{eq 2}
		2\lambda_{min}(c+\epsilon_k)(c+\epsilon_k) = \int_{\mathbb{R}^n}|\nabla u_{c+\epsilon_k}|^2-f(x,u_{c+\epsilon_k})u_{c+\epsilon_k},
	\end{align}
	\eqref{eq 1} comes from $\tmop{Eq}_{\lambda_{min}(c)}(u_c) = 0$ and $\int_{\mathbb{R}^n}|u_c|^2 = 2c$. A similar way leads to \eqref{eq 2}. \eqref{eq 2} minusing \eqref{eq 1} yields
	\begin{align}
		& 2(\lambda_{min}(c+\epsilon_k) - \lambda_{min}(c))c \nonumber \\
		= & \left\langle L_{u_c}(u_c),u_{c+\epsilon_k}-u_c \right\rangle + \lambda_{min}(c)\left\langle u_c,u_{c+\epsilon_k}-u_c \right\rangle + o(\|u_{c+\epsilon_k}-u_c\|_{H^1}) + O(\epsilon_k),
	\end{align}
	implying that
	\begin{align}
		& 2\limsup_{k \to \infty}\frac{\lambda_{min}(c+\epsilon_k) - \lambda_{min}(c)}{\epsilon_k} \nonumber \\
		= & \limsup_{k \to \infty}\frac{1}{c}\left( \left\langle L_{u_c}(u_c),\frac{u_{c+\epsilon_k}-u_c }{\epsilon_k}\right\rangle + \lambda_{min}(c)\left\langle u_c,\frac{u_{c+\epsilon_k}-u_c}{\epsilon_k} \right\rangle +O(1) \right)  \nonumber \\
		< & \infty,
	\end{align}
	since $\lim_{k \to \infty}\|\frac{u_{c+\epsilon_k} - u_c}{\epsilon_k}\|_{H^1} < \infty$. This shows that \eqref{eq right res} holds true. Similarly we can prove \eqref{eq left res}, which completes the proof.
\end{proof}

\subsubsection{Proof of Lemma \ref{< 0}} \label{sec. < 0}

\begin{proof}
	We just address the case of $\lambda = \lambda_{min}(c)$ here. The proof of the other case is similar. Take $c_k \to c^+$ and $\lambda_k = \lambda_{min}(c_k)$. Let $u_k \in S_{c_k}$ satisfying $\tmop{Eq}_{\lambda_k}(u_k) = 0$. By the proof of Proposition \ref{properties of m(c)}, up to a subsequence, $\lambda_k \to \lambda^-$ and $u_k \to u_c$ strongly in $H^1$. By $(A 2)$, $\tmop{Ker} L_u \cap L^2_{\tmop{rad}} = \{ 0 \}$. Then using Proposition \ref{constbranch}, we construct
	the branch $\varepsilon \rightarrow U_{\lambda +\varepsilon} \in C^1 ([-\varepsilon_0, \varepsilon_0], L^2_{\tmop{rad}})$, solution of $\tmop{Eq}_{\lambda + \varepsilon}
	(U_{\lambda + \varepsilon}) = 0$ for $| \varepsilon |$ small. Moreover, $U_{\lambda} = u_c$. Since $u_k \to u_c$ strongly in $H^1$, using Proposition \ref{uniquresult}, one gets that $u_k = U_{\lambda + \varepsilon_k}$ for $k$ large enough and a family of small $(\varepsilon_k)$ and $\lambda_k = \lambda + \varepsilon_k$. Then direct computation yields that
	\begin{align}
		\langle L_{u_c}^{-1}u_c,u_c\rangle & = 2\partial_\varepsilon Q(U_{\lambda+\varepsilon})|_{\varepsilon = 0} \nonumber \\
		& = 2\lim_{k \to \infty}\frac{Q(u_k)-Q(u_c)}{\lambda_k-\lambda} \nonumber \\
		& = 2\lim_{c_k \to c^+}\frac{c_k-c}{\lambda_{min}(c_k)-\lambda_{min}(c)} \nonumber \\
		& \leq \frac{2}{\liminf_{\tilde{c} \to c^+}\frac{\lambda_{min}(\tilde{c})-\lambda_{min}(c)}{\tilde{c}-c}}.
	\end{align}
	By Proposition \ref{properties of m(c)}, $\frac{\lambda_{min}(\tilde{c})-\lambda_{min}(c)}{\tilde{c}-c} < 0$ and $\liminf_{\tilde{c} \to c^+}\frac{\lambda_{min}(\tilde{c})-\lambda_{min}(c)}{\tilde{c}-c} > -\infty$. This yields $\langle L_{u_c}^{-1}u_c,u_c\rangle < 0$ and completes the proof.
\end{proof}

\subsubsection{Proof of Lemma \ref{=0 is finite}} \label{sec. =0 is finite}

\begin{proof}
	Suppose on the contrary that there exists $(c_k)_{k\in \mathbb{N}}$ such that $\langle L_{u_k}^{-1}u_k,u_k\rangle = 0$ for some $u_k \in S_{c_k}$.
	
	
	Then up to a subsequence, assume that $c_k \to \tilde{c} \in ]c_{\ast},c^{\ast}[$. Then $E(u_k) = m(c_k) \to m(\tilde{c})$ and $Q(u_k) = c_k \to \tilde{c}$. By $(A 3)$, we obtain that $u^\ast \in S_{\tilde{c}}$ such that $u_k \to u^\ast$ in $H^1$. Let $\lambda_k, \lambda^\ast$ be the corresponding Lagrange multipliers of $u_k, u^\ast$ respectively. 
	By $(A 2)$ $u^\ast$ is non-degenerate in $H^1_{rad}$. Then by the Proposition \ref{constbranch}, there exists $\varepsilon \rightarrow U_{\lambda^\ast + \varepsilon} \in C^1 ([-\varepsilon_0, \varepsilon_0], L^2_{\tmop{rad}})$
	such that $U_{\lambda^\ast} = u^\ast$, $U_{\lambda^\ast + \varepsilon}$ is radial and solves the equation
	\[ (-\Delta + V) U_{\lambda^\ast + \varepsilon} = (\lambda^\ast + \varepsilon) U_{\lambda^\ast +
		\varepsilon} + f(x,U_{\lambda^\ast + \varepsilon}). \]
	Using Proposition \ref{uniquresult}, one gets that $u_k = U_{\lambda_k}$ for large $n$. Let $\varepsilon_k = \lambda_k - \lambda^\ast \to 0$. We assume that $\sup_{k}|\varepsilon_k| < \varepsilon_0$. Then direct computation yields that $\partial_\varepsilon Q(U_{\lambda^\ast+\varepsilon})|_{\varepsilon = \varepsilon_k} = \langle L_{U_{\lambda_k}}^{-1}U_{\lambda_k},U_{\lambda_k}\rangle = 0$. Also, using the strong convergence of $U_{\lambda_k} \to U_{\lambda^\ast}$ we obtain $\partial_\varepsilon Q(U_{\lambda^\ast+\varepsilon})|_{\varepsilon = 0} = \langle L_{U_{\lambda^\ast}}^{-1}U_{\lambda^\ast},U_{\lambda^\ast}\rangle = 0$.
	
	Now let us set $M(\varepsilon) = Q(U_{\lambda^\ast+\varepsilon})$. We find a sequence $(\varepsilon_k)_{k \in \mathbb{N}}$ with $\varepsilon_k \to 0$ such that $M(\varepsilon_k) \neq M(\varepsilon_j) \neq M(0)$ for any $k \neq j$ and $M'(\varepsilon_k) = M'(0) = 0$. However, note that $M$ is a real-analytic function on $(-\varepsilon_0,\varepsilon_0)$. This leads to a contradiction. 
\end{proof}

\subsubsection{Proof of Proposition \ref{simo}}

\begin{proof}
	We recall the notations of Proposition \ref{simo}. Take $\tilde{c} \in ]c_{\ast},c^{\ast}[$ and
	denote $(u_{j,min})_{j \in [1 \ldots N_1]}, N_1 \geqslant 1$ the elements in $S_{\tilde{c}}$ such that $\tmop{Eq}_{\lambda_{min}(\tilde{c})}u_{j,min} = 0$. We infer that this set is finite. Indeed, if it is infinite, we can take a sequence which is strongly convergent. Then using Propositions \ref{constbranch} and \ref{uniquresult} we obtain a contradiction since elements in this sequence have the same Lagrange multiplier. By Proposition
	\ref{constbranch}, for any $j \in [1 \ldots N_1]$ we can extend $u_{j,min}$ to a
	branch $U_{j, \varepsilon,min}$ for $| \varepsilon | \leqslant
	\mathbf{\epsilon}_j$ for some $\mathbf{\epsilon}_j > 0$, with $U_{j, 0,min}
	= u_{j,min}$. We define $\varepsilon_0 \assign \inf_{j \in [1 \ldots N_1]}
	\mathbf{\epsilon}_j$ and $\varepsilon_0 > 0$ since there are a finite
	number of branches.
	
	First, we infer that for any $\eta > 0$ there exists $\tau_0 > 0$ such that if
	$c
	\in ]\tilde{c} , \tilde{c} + \tau_0]$, then
	\[ S_c \subset \bigcup_{j \in [1 \ldots N_1]} B_{H^1} (u_{j,min}, \eta) . \]
	Here, $B_{H^1} (u_{j,min}, \eta)$ is the ball in $H^1 (\mathbb{R}^n)$ centered at
	$u_{j,min}$ and of radius $\eta$. This is a direct consequence of the proof of Proposition \ref{properties of m(c)}.
	
	We apply this for $\eta > 0$ small for which we can apply Proposition
	\ref{uniquresult} for all the $u_{j,min}$. For this $\eta$, we deduce that there
	exists $\tau_0 > 0$ such that if $| c - \tilde{c} | \leqslant \tau_0$, then any $v
	\in S_c$ satisfy $\tmop{Eq}_{\lambda} (v) = 0$ for some $\lambda < \Lambda_0$ by $(A
	1)$, and there exists $j \in [1 \ldots N_1]$ such that $\| v - u_{j,min} \|_{H^1
		(\mathbb{R}^n)} \leqslant \eta$. By Proposition \ref{uniquresult}, we deduce
	that there exists $\varepsilon_j$ such that $v = U_{j, \varepsilon_j}$. This
	concludes the proof of the case of $c\in ]\tilde{c} , \tilde{c} + \tau_0]$. Similarly, we can prove the case when $c
	\in [\tilde{c} - \tau_0 , \tilde{c}[$ and reach the desired conclusions of this proposition.
\end{proof}

\subsubsection{Proof of Proposition \ref{ext 1}}\label{sec. ext 1}

\begin{proof}
	By Proposition \ref{simo} there
	exists $\tau_0 > 0$ such that for $c \in ]\tilde{c}, \tilde{c} + \tau_0]$,
	$S_c \subset \{ U_{\varepsilon, min} \}$ for some small $| \varepsilon |$. By Lemma \ref{< 0}, $\langle L^{-1}_{u_{min}} (u_{min}), u_{min} \rangle < 0$. By Proposition \ref{coerc2}, $\partial_{\varepsilon} (Q (U_{\varepsilon,min}))_{| \varepsilon = 0} < 0$, thus $Q (U_{\varepsilon,min})$ is strictly decreasing near $0$. This implies that $| S_c | \leqslant 1$. By $(A 1)$ we have $| S_c | \geqslant 1$, hence $| S_c | = 1$. By the monotonicity of $\mathcal{T}_c$ (see Proposition \ref{properties of m(c)}), we can take $\varepsilon_0 > 0$ such that $U_{-\varepsilon_0,min} = u_{\tilde{c} + \tau_0}$ and $S_c \subset \{ U_{\varepsilon, min} \}$ for $\varepsilon \in [-\varepsilon_0,0[$. On the other hand, for any $\varepsilon \in [-\varepsilon_0,0[$, $Q(U_{\varepsilon,min}) \in ]\tilde{c}, \tilde{c} + \tau_0]$ and for any $\varepsilon \neq \tilde{\varepsilon} \in [-\varepsilon_0,0[$ $Q(U_{\varepsilon,min}) \neq Q(U_{\tilde{\varepsilon},min})$. this implies that $\{U_{\varepsilon, min}, \varepsilon \in [-\varepsilon_0,0[\} \subset \{u_c, c \in ]\tilde{c}, \tilde{c} + \tau_0]\}$. Hence we complete the proof when $c \in ]\tilde{c}, \tilde{c} + \tau_0]$. The proof when $c \in [\tilde{c} - \tau_0, \tilde{c}[$ is similar and we complete the proof.
\end{proof}

\subsubsection{Proof of Proposition \ref{ext 2}}

\begin{proof}
	Let $(u_{j,min})_{j \in [1 \ldots N_1]}, N_1 \geqslant 1$ be elements in $S_{c}$ such that $\tmop{Eq}_{\lambda_{min}(c)}(u_{j,min}) = 0$. Recall that $U_{j, \varepsilon, min}, j \in [1 \ldots N_1]$ are the associated branches as defined in Proposition \ref{constbranch}. These branches are divided into two categories. 
	The first one (that might be empty) is then $u_{j,min}$ is isolated in $S := \{u \in S_c, c \in ]c_{\ast},c^{\ast}[\}$. The other one is when there exists a sequence of elements in $S$ converging to $u_{j,min}$. Note that the latter is non-empty and we can take such convergent sequence on the branch stemming from $u_{j,min}$ by Proposition \ref{simo}. Let $(u_{j,min})_{j \in [1 \ldots L]}, 1 \leq L \leq N_1$ be in this one, and there exist negative sequences $(\varepsilon_k^{(j)})_k, j \in [1 \ldots L]$ such that $Q(U_{j,\varepsilon_k^{(j)},min}) \in S_{c_k^{(j)}}$ with $c_k^{(j)} \to c^+$. Since $f(x,u)$ is real analytic w.r.t. $u$, we know that both $Q(U_{j,\varepsilon,min})$ and $E(U_{j,\varepsilon,min})$ are real analytic w.r.t. $\varepsilon \in [-\varepsilon_0, \varepsilon_0]$ for some small $\varepsilon_0 > 0$ using the implicit function theorem. Thus $\Phi_{\lambda+\varepsilon}(U_{j,\varepsilon,min})$ is also real analytic w.r.t. $\varepsilon \in [-\varepsilon_0, \varepsilon_0]$ where $\Phi_{\lambda} := E - \lambda Q$. By \cite[Appendix C]{Hajaiej-Song-unique} we know any $u \in S_c$ with $\tmop{Eq}_\lambda(u) = 0$ is GSS for $\Phi_{\lambda}$ (a minimizer constrainted on Nehari manifold). For $i \neq j \in [1 \ldots L]$, we have
	$$
	\Phi_{\lambda_{min}(c)+\varepsilon_k^{(j)}}(U_{j,\varepsilon_k^{(j)},min}) \leq \Phi_{\lambda_{min}(c)+\varepsilon_k^{(j)}}(U_{i,\varepsilon_k^{(j)}},min),
	$$
	$$
	\Phi_{\lambda_{min}(c)+\varepsilon_k^{(i)}}(U_{i,\varepsilon_k^{(i)}},min) \leq \Phi_{\lambda_{min}(c)+\varepsilon_k^{(i)}}(U_{j,\varepsilon_k^{(i)}},min).
	$$
	This implies that $\Phi_{\lambda_{min}(c)+\varepsilon}(U_{j,\varepsilon,min}) \equiv \Phi_{\lambda_{min}(c)+\varepsilon}(U_{i,\varepsilon,min})$ on $[-\varepsilon_0, \varepsilon_0]$. Also using \cite[Appendix C]{Hajaiej-Song-unique} we know that $Q(U_{j,\varepsilon,min}) \equiv Q(U_{i,\varepsilon,min}), E(U_{j,\varepsilon,min}) \equiv E(U_{i,\varepsilon,min})$. We can take $\tau_0 > 0$ small enough such that for $c \in ]\tilde{c}, \tilde{c} + \tau_0]$, $S_c \cap E(U_{j,\varepsilon,min}) = \emptyset$ for any $j \in [L+1,N_1]$. Note that $Q (U_{j,\varepsilon,min})$ is strictly decreasing in $[-\varepsilon_0, \varepsilon_0]$ by its real analyticity for any $j \in [1 \ldots L]$. Then similar to the proof of Proposition \ref{ext 1}, we can take $U_{j,-\varepsilon_0,min} \in S_{\tilde{c} + \tau_0}$ for any $j \in [1 \ldots L]$ and show that $\{u \in S_c, c \in ]\tilde{c}, \tilde{c} + \tau_0]\} = \{U_{j, \varepsilon, min}, \varepsilon \in [-\varepsilon_0,0[, j \in [1 \ldots L]\}$. This completes the case when $c \in ]\tilde{c}, \tilde{c} + \tau_0]$. The proof when $c \in [\tilde{c} - \tau_0, \tilde{c}[$ is similar.
\end{proof}

\section{Assumptions on $V$ and $f$ such that properties $(A)$, $(B)$, $(C)$ and $(N)$ are satisfied} \label{sec. conditions}

Now we start our second part, showing that under some assumptions on $V$ and $f$, properties $(A)$, $(B)$, $(C)$ and $(N)$ are satisfied, and giving some examples of classical equations that this covers.

\
\subsection{Hypotheses on $f$ and $V$}

To show $(A 3)$ we assume

$(H 1)$ $(h_1)$-$(h_5)$ and $(h_8)$ hold; or $(h_1)$, $(h_5)$-$(h_7)$, $(h_8)$ hold with $\gamma = -n-2$; or $(h_1)$, $(h_5)$-$(h_7)$, $(h_9)$ hold with $\gamma = -n-\Gamma$ where $\Gamma = \min\{2,\frac{n}{q},-\alpha\}$, $q = \frac{n}{2}$ when $n \geq 3$, $q > 1$ when $n = 1,2$; \\
where

$(h_1)$ $\lim_{t \to 0}\frac{f(x,t)}{t} = 0$ and $\lim_{|t|\to\infty}\frac{f(x,t)}{|t|^{1+4/n}} = 0$ uniformly for $x \in \mathbb{R}^n$.

$(h_2)$ $f(x,t) \to \bar{f}(t)$ as $|x| \to \infty$ uniformly for $t$ bounded, and $f(x,t) \geq \bar{f}(t)$.

$(h_3)$ $\exists \xi > 0$ such that $\bar{F}(\xi) := \int_0^\xi\bar{f}(\tau)d\tau > 0$.

$(h_4)$ $\liminf_{t \to 0}\frac{\bar{F}(t)}{t^{2+4/n}} = +\infty$.

$(h_5)$ For any $k > 1$, $F(x,kt) > k^2F(x,t), \forall x \in \mathbb{R}^n, t \neq 0$.

$(h_6)$ $\forall \epsilon > 0, \exists K > 0, \forall |x| > K, |F(x,t)| \leq \epsilon(|t|^2 + |t|^{2+4/n})$.

$(h_7)$ $\liminf_{\tau \to 0^+}\tau^{\gamma}F(\frac{x}{\tau},\tau^{\frac{n}{2}}t) = +\infty, \forall x \in \mathbb{R}^n, \forall t \neq 0$, where $\gamma$ will be determined in the following results.

$(h_8)$ $V(x) \leq \lim_{|x| \to \infty}V(x) = 0$.

$(h_9)$ There exists $\alpha < 0$ such that $V^+(tr) \lesssim t^\alpha V^+(r)$ for $tr > 1$ where $V^+ = \max\{V,0\}$.

\

To show $(A 1)$ we further assume

$(H 2)$ $(h_{10})$-$(h_{12})$ hold; \\
where

$(h_{10})$ $f(x,t)$ is radial w.r.t. $x \in \mathbb{R}^n$ and nonincreasing w.r.t. $|x| > 0$, $V(x)$ is radial w.r.t. $x \in \mathbb{R}^n$ and nondecreasing w.r.t. $|x| > 0$.

$(h_{11})$ $\Lambda_0 = 0$, $2V+x\cdot \nabla_xV \leq 0$, and $2nF(x,t) + 2x\cdot \nabla_xF(x,t) - (n-2)f(x,t)t > 0$ for any $t > 0$; 

$(h_{12})$ $V \equiv 0$ and $f(x,t) = f(t)$; or there exists $d \in \mathbb{R}$ such that $f_t(x,t)t - f(x,t) + d(2\lambda t - x\cdot\nabla_xV t + 2f(x,t) + x\cdot\nabla_xf(x,t))$ does not change sign for all $t > 0$, $\lambda < 0$.
\

To show $(A 2)$ we further assume

$(H 3)$ $(h_{10})$ and $(h_{13})$ hold; \\
where

$(h_{13})$ $f_t(x,t)t > f(x,t)$ for any $t > 0$. And there exists $k \in \mathbb{R}$ such that $f_t(x,t) + (2k -1)\frac{f(x,t)}{t} + k\frac{x\cdot\nabla_x(f(x,t)+Vt)}{t}$ or $kf_t(x,t) + (2-k)\frac{f(x,t)}{t} + \frac{x\cdot\nabla_x(f(x,t)+Vt)}{t}$ is non-increasing in $r = |x|$ and increasing in $t > 0$ or non-decreasing in $r = |x|$ and decreasing in $t > 0$.

\

To show $(B)$ we assume

$(H 4)$ $(h_{10})$, $(h_{14})$, $(h_{15})$ hold and $V \equiv 0$; \\
where

$(h_{14})$ $\lim_{\tau \to 0}\frac{p\tau^{-\bar{p}}F(\tau^{\frac{2-\bar{p}}{2}}x,\tau t)}{|x|^\theta |t|^p} = 1, \forall x \neq 0, t \neq 0$, where $\bar{p}, p, \theta$ satisfy
$$
\theta \leq 0, \quad \bar{p} = p + \frac{2-\bar{p}}{2}\theta, \quad 2 < p \leq \bar{p} < 2+4/n.
$$
Moreover, for $|\tau| \leq 1$, there exist $2 < r \leq s < 2 + 4/n$ such that
\begin{equation} \label{eq near 0}
	|F(\tau^{\frac{2-\bar{p}}{2}}x,\tau t)| \lesssim \left\{
	\begin{array}{cc}
		|\tau|^{\bar{p}}(|t|^r+|t|^s), & |\tau^{\frac{2-\bar{p}}{2}}x| \geq 1, \\
		|\tau|^{p}(|t|^r+|t|^s), & |\tau^{\frac{2-\bar{p}}{2}}x| < 1,
	\end{array}
	\right.
\end{equation}
and $\theta \geq \frac{\max\{r,s\}-2}{2}n-2$ if $n \geq 2$, $\theta > -1$ if $n = 1$. 

$(h_{15})$ $|F(x,t)| \lesssim |t|^r+|t|^s$ for some $2 < r \leq s < 2 + 4/n$, and $$\lim_{|\tau| \to \infty}\frac{q\tau^{-q}F(\tau^{\frac{2-q}{2}}x,\tau t)}{|t|^q} = 1, \forall x \in \mathbb{R}^n, t \neq 0,$$ where $2 < q < 2+4/n$.

\

To show $(C)$ we assume

$(H 5)$ $(h_{16})$-$(h_{18})$ hold or $(h_{19})$ holds or $(h_{20})$ holds. \\
where

$(h_{16})$ $n \geq 2$, $f(x,t) = f(t), V \equiv 0$ and for any $\lambda < 0$, there exist $b > 0$ and $\tilde{b} \in (b,\infty]$ such that $f(t) + \lambda t < 0$ on $(0, b)$,  $f(t) + \lambda t > 0$ on $(b,\tilde{b})$ and $f'(b) + \lambda > 0$. If $\tilde{b} < \infty$, then $f(t) + \lambda t < 0$ on $(\tilde{b}, \infty)$ and $f'(\tilde{b}) + \lambda < 0$.

$(h_{17})$ For any $\lambda < 0$, there exists $\xi > b$ such that $\int_0^\xi(f(t)+\lambda t)dt > 0$.

$(h_{18})$ For any $\lambda < 0$, $\frac{f'(t)t+\lambda t}{f(t) + \lambda t}$ is no more than $1$ on $(0,b)$ and decreasing on $(b,\tilde{b})$.

$(h_{19})$ $n \geq 2$, $f(u) = \Sigma a_i|u|^{p_i-2}u - \Sigma b_j|u|^{q_j-2}u$ where $a_i, b_j$ are positive and $2 < q_j < p_i \leq \frac{2n-2}{n-2}$.

$(h_{20})$ $f(x,u) = h(|x|)|u|^{p-2}u$ with $(n - 2s)p < 2(n + \theta)$, $h(r) > 0$, $h(r)$ and $\frac{rh'(r)}{h(r)}$ are non-increasing in $(0, \infty)$, $\theta := \lim_{r\to \infty}\frac{rh'(r)}{h(r)} > -\infty$.

\begin{proposition}
	\label{dubai}Suppose that conditions $(H 1)$-$(H 3)$ are satisfied, then properties $(A)$ hold.
	
	Furthermore, if conditions $(H 4)$ are satified, then property $(B)$ holds.
	
	Finally, if conditions $(H 5)$ are satified, then property $(C)$ holds.
\end{proposition}

	We give some examples here.
	
\begin{example}[Focusing nonlinearties]
	Let:
	\begin{itemizeminus}
		\item $V\equiv 0, f(u) = \Sigma_{i=1}^Na_i|u|^{p_i-2}u$, $a_i > 0$, $p_i \in ]2,2+\frac{4}{n}[$.
    \end{itemizeminus}
	
		It is easy to see that $(h_1)$-$(h_5)$, $(h_8)$ and $(h_{10})$-$(h_{13})$ hold. Thus $(H 1)$-$(H 3)$ are satisfied. Then by Proposition \ref{dubai} the properties $(A)$ hold. (In this case, $]c_{\ast},c^{\ast}[ = ]0,\infty[$.)
			
		Furthermore, we can verify that $(h_{14})$ and $(h_{15})$ hold with $p = \min p_i$, $\theta = 0$ and $q = \max p_i$. Hence, property $(B)$ holds.
		
		Additionally, it is not difficult to verify $(N)$ and that $f(u)$ is real analytic w.r.t. $u$, we don't need property $(C)$ and all assumptions in Theorems \ref{mainthing0} and \ref{mainthing2} are satisfied.
		
		Moreover, properties $(A)$ and $(B)$ hold when some $p_i$ is $2+\frac{4}{n}$ and others belong to $]2,2+\frac{4}{n}[$. We set $f(u) = |u|^{p-2}u + |u|^{\frac{4}{n}}u, p \in ]2,2+\frac{4}{n}[$ for simplicity. In this case, $]c_{\ast},c^{\ast}[ = ]0,c^{\ast}[$ where $c^{\ast} = \int_{\mathbb{R}^n}|Q|^2 < \infty$, where $Q$ is the unique positive radial-decreasing solution to the equation $\Delta Q + |Q|^{p-2}Q - Q = 0$. Note that the uniqueness of NGSS when $c \to c^{\ast}$ needs additional arguments. We briefly give the key steps here. First let us show that $\lim_{c \to c^{\ast}}m(c) = -\infty$. Let $u_\lambda$ be the a GSS (or positive solution) with $\tmop{Eq}_\lambda(u_\lambda) = 0$. Next we prove that $u_\lambda$ is unique as $\lambda \to 0^-$ and $Q(u_\lambda) \to 0$; $u_\lambda$ is unique as $\lambda \to -\infty$, $Q(u_\lambda) \to c^{\ast}$ and $\partial_\lambda Q(u_\lambda) < 0$. Then we argue by contradiction and assume that there exists $u_n \in S_{c_n}$ with $c_n \to c^{\ast}$ and $\lambda_n \nrightarrow -\infty$ where $\lambda_n$ is the Lagrange multiplier. Up to a subsequence, we can assume that $\lambda_n \to \lambda_\infty \in (-\infty,0)$. This implies that $m(c_n)$ is bounded, in a contradiction with $\lim_{c \to c^{\ast}}m(c) = -\infty$. Hence we complete the proof. More details of the proof are left to interested readers.	
			
	Finally, once $p_i > 2+\frac{4}{n}$, then $m(c) = -\infty$ and NGSS does not exist.
\end{example}

\begin{example}[Non-autonomous focusing nonlinearities, I]
	Let:
	\begin{itemizeminus}
		\item $V\equiv 0, f(x,u) = \Sigma_{i=1}^Na_i(|x|)|u|^{p_i-2}u$, $p_i \in ]2,2+\frac{4}{n}[$, $a_i(r) > 0$ is $C^1$ with respect to $r \in \mathbb{R}^+$, $a_i(r) \to a_i(\infty) > 0$ as $r \to \infty$, $a_i'(r) \leq 0$ and $2n + 2\inf_{r>0}\frac{ra_i'(r)}{a_i(r)} > (n-2)p_i$.
    \end{itemizeminus}
	
		It is easy to see that $(h_1)$-$(h_5)$, $(h_8)$ and $(h_{10})$-$(h_{13})$ hold. Thus $(H 1)$-$(H 3)$ are satisfied. Then by Proposition \ref{dubai}, properties $(A)$ hold. (In this case, $]c_{\ast},c^{\ast}[ = ]0,\infty[$.)
		
		Furthermore, we can verify that $(h_{14})$ and $(h_{15})$ hold with $p = \min p_i$, $\theta = 0$ and $q = \max p_i$. Hence, property $(B)$ holds.
		
		Additionally, it is not difficult to verify $(N)$ and that $f(x,u)$ is real analytic w.r.t. $u$, we don't need property $(C)$ and all assumptions in Theorems \ref{mainthing0} and \ref{mainthing2} are satisfied.
\end{example}

\begin{example}[Nonautonomous focusing nonlinearities, II]
	Let:
	\begin{itemizeminus}
		\item $V\equiv 0, f(x,u) = \Sigma_{i=1}^Na_i(|x|)|u|^{p_i-2}u$, $p_i \in ]2,2+\frac{4}{n}[$, $a_i(r) > 0$ is $C^1$ w.r.t. $r \in \mathbb{R}^+$, $a_i(r) \to 0$ as $r \to \infty$, $a_i'(r) \leq 0$, $\liminf_{\tau \to 0^+}\tau^{\frac{np_i}{2}-n-2}a_i(\frac{r}{\tau}) = +\infty$ and $2n + 2\inf_{r>0}\frac{ra_i'(r)}{a_i(r)} > (n-2)p_i$.
    \end{itemizeminus}
	
		It is easy to see that $(h_1)$, $(h_5)$-$(h_7)$, $(h_8)$ and $(h_{10})$-$(h_{13})$ hold with $\gamma = -n-2$. Thus $(H 1)$-$(H 3)$ are satisfied. Then by Proposition \ref{dubai} the properties $(A)$ hold. (In this case, $]c_{\ast},c^{\ast}[ = ]0,\infty[$.)
		
		To show $(B)$, we further assume that: $\frac{a_i(r)}{r^{\theta_i}} \to 1$ as $r \to \infty$ with $0 > \theta_i \geq \frac{\max\{r,s\}-2}{2}n-2$ if $n \geq 2$, $0 > \theta_i > -1$ if $n = 1$. Then $(h_{14})$ and $(h_{15})$ hold with $p = p_k$, $\theta = \theta_k$ where $k$ is such that $\frac{2(p_i+\theta_i)}{2+\theta_i}$ achieves minimum and $q = \max p_i$. Hence, the property $(B)$ holds.
		
		Additionally, it is not difficult to verify $(N)$ and $f(x,u)$ is real analytic w.r.t. $u$, we don't need property $(C)$ and all assumptions in Theorems \ref{mainthing0} and \ref{mainthing2} are satisfied.		
\end{example}

\begin{example}[Focusing-defocusing nonlinearities]
	Let:
	\begin{itemizeminus}
		\item $V\equiv 0, f(u) = |u|^{p-2}u - |u|^{q-2}u$, $2 < p < q$, $p < 2+\frac{4}{n}$, $n \geq 2$.
\end{itemizeminus}	

It is easy to see that $(h_1)$, $(h_2)$, $(h_4)$, $(h_8)$, $(h_{10})-(h_{12})$ hold. $(h_3)$ was verified in \cite{Lewin-Nodari-2020}. Note that $(h_5)$ is not satisfied. However, $(h_5)$ is used to show the strict decreasiness of $m(c)$. In Appendix \ref{app proof}, we prove this result without using $(h_5)$. Hence, along the same lines in this paper, we can show $(A 1)$ and $(A 3)$. By \cite[Theorem 1]{Lewin-Nodari-2020} $(A 2)$ and $(C)$ hold. As for $(B)$, it can be obtained by \cite[Theorem 7]{Lewin-Nodari-2020}. In fact, since $(h_{14})$ is satisfied, the uniqeness for masses close to $0$ can be proved using the method in our paper. However, it needs additional discussions like \cite[Theorem 7]{Lewin-Nodari-2020} for large masses. Also, we believe that there is an approach like the one in our paper to address large masses.

When $2 < p = 2+\frac{4}{n} < q$, we can get $(A)$, $(B)$ and $(C)$ using similar arguments and references to the case that $p < 2+\frac{4}{n}$. In this case, $]c_\ast,c^\ast[ = ]c_\ast,\infty[$ where $c_\ast = \int_{\mathbb{R}^n}|Q|^2dx > 0$, $Q$ is the unique positive radial-decreasing solution to the equation $\Delta Q + |Q|^{p-2}Q - Q = 0$.

When $2 < p < q$ and $2+\frac{4}{n} < p < \frac{2n}{n-2}$ is more complicated. In fact, the existence range is $[c_\ast,\infty[$ for some $c_\ast > 0$, and we don't need the uniqueness of NGSS when $c \to c_\ast$. Moreover, $(A)$, $(C)$ and the uniqueness of NGSS when $c \to \infty$ can be verified using similar arguments and references to the case that $p < 2+\frac{4}{n}$. 

\end{example}


\begin{remark}[Examples with non-constant potentials]
	We can cover the cases where a potential appears and the first simple eigenvalue exists with a positive eigenfunction. For example, we can provide examples such as:
	\begin{itemizeminus}
		\item $\Lambda_0$ is the first simple eigenvalue of $-\Delta + V$ with a positive eigenfunction and $f(u) = |u|^{p-2}u$, $p \in ]2,2+\frac{4}{n}[$. And more conditions on $V$ should be given to satisfy some hypotheses.
	\end{itemizeminus}
	Nonlinearities having more complicated forms can also be provided. The additional features we need to do is to show that NGSS with small mass is unique. In fact, we can show $Q(u_\lambda) \to \infty$ as $\lambda \to -\infty$. Let $u_c \in S_c$ and $\lambda_c$ be the Lagrange multiplier. Then $\lambda_c$ is bounded when $c \to 0$. Then by $(h_1)$ and
	\begin{align}
		\int_{\mathbb{R}^n}|\nabla u_c|^2 + V|u_c|^2 = \lambda \int_{\mathbb{R}^n}|u_c|^2 + \int_{\mathbb{R}^n}f(x,u_c)u_c
	\end{align}
	and that $\int_{\mathbb{R}^n}|u_c|^2 \to 0$, we have $u_c \to 0$ in $H^1$ and can use classical bifurcation theory to show $\lambda_c \to \Lambda_0$. Moreover, we can show $\partial_\lambda Q(u_\lambda) \neq 0$ as $\lambda \to \Lambda_0$, then we deduce the uniqueness of $u_c$ with $c$ small enough. All the techniques were initiated in \cite{Song-Hajaiej-partial confinement} and we don't give more details here.
	
We assume the boundedness of $V$ for simplicity in this paper. In fact, after slight modifications, we can cover the case $V(x) = |x|^2$, the harmonic potential. In this case ($V(x) = |x|^2$, $f(u) = |u|^{p-2}u$, $p \in ]2,2+\frac{4}{n}[$), $]c_{\ast},c^{\ast}[ = ]0,\infty[$ and $(A)$, $(B)$, $(C)$ can be verified.
\end{remark}


\subsection{Properties of the minimizers}

We summarize here the steps of the proof of Proposition \ref{dubai}. The
results of this subsections are proven in section \ref{sono}.

\subsubsection{Sketch of the Proof of $(A 1)$ and $(A 3)$ }

First we show that $(A 3)$ holds under conditions $(H 1)$. The fact that $(h_1)$-$(h_5)$, $(h_8)$ hold was proved in \cite{Yang-Qi-Zou-2022}. We collect some their results and omit the proofs.

\begin{lemma} \label{lem compact1}
	Assume that $(h_1)$-$(h_5)$, $(h_8)$ hold. Then
	\begin{itemize}
		\item[(i)] $-\infty < m(c) < 0$.
		\item[(ii)] $m(c)$ is strictly decreasing and continuous for $c \in ]0,\infty[$.
		\item[(iii)] Any sequence $\{u_k\} \subset H^1$ having the properties
		$$
		E(u_k) \rightarrow m(c) \ and \ Q(u_k) \rightarrow c
		$$
		has a converging subsequence in $H^1$.
	\end{itemize}
\end{lemma}

Next we consider other cases in $(H 1)$. The first step is to show $m(c) \in (-\infty,0)$.

\begin{lemma}
	\label{lem estimate of m(c)} Assume that $(h_1)$, $(h_7)$, $(h_8)$ hold with $\gamma = -n-2$; or that $(h_1)$, $(h_7)$, $(h_9)$ hold with $\gamma = -n-\Gamma$ where $\Gamma = \min\{2,\frac{n}{q},-\alpha\}$, $q = \frac{n}{2}$ when $n \geq 3$, $q > 1$ when $n = 1,2$. Then $-\infty < m(c) < 0$.
\end{lemma}

The second step is to prove the strict decreasiness and continuity of $m(c)$ w.r.t. $c \in ]0,\infty[$.
\begin{lemma}
	\label{lem decreasedness}  Under the assumptions of Lemma \ref{lem estimate of m(c)}. Further assume that $(h_5)$ holds.
	\begin{itemize}
		\item[(i)] $m(c)$ is strictly decreasing and continuous with respect to $c \in ]c_{\ast},c^{\ast}[$.
		\item[(ii)] Let $c_1, c_2 > 0$ be such that $c = c_1 + c_2$. Then $m(c) \leq m(c_1) + m(c_2)$. Furthermore, if $m(c_1)$ or $m(c_2)$ can be achieved, then $m(c) < m(c_1) + m(c_2)$.
	\end{itemize}
\end{lemma}

Now we can prove the following compactness result.

\begin{lemma}
	\label{lem coompactness} Assume that $(h_1)$, $(h_5)$-$(h_7)$, $(h_8)$ hold with $\gamma = -n-2$; or that $(h_1)$, $(h_5)$-$(h_7)$, $(h_9)$ hold with $\gamma = -n-\Gamma$ where $\Gamma = \min\{2,\frac{n}{q},-\alpha\}$, $q = \frac{n}{2}$ when $n \geq 3$, $q > 1$ when $n = 1,2$. Then any sequence $\{u_k\} \subset H^1$ having the properties
	$$
	E(u_k) \rightarrow m(c) \ and \ Q(u_k) \rightarrow c
	$$
	has a converging subsequence in $H^1$.
\end{lemma}

Now we have shown $(A 3)$. Then by some standard discussions we can prove $(A 1)$.

\begin{proposition} \label{A1 and A3}
	$(A 1)$ and $(A 3)$ hold under the conditions $(H 1)$ and $(H 2)$.
\end{proposition}

See subsection \ref{sec. A1 A3} for its proof.

\subsubsection{Sketch of the Proof of $(A 2)$}

The proof of $(A 1)$ shows that $u = u(r)$ is radial, positive and decreasing in $r > 0$. Hence we need to establish $(A_2)$ in the radially symmetric space $H^1_{rad}$. We argue by contradiction. Since $L_u$ has exactly one simple negative eigenvalue, $0$ is its second eigenvalue. Let $w = w(r)$ be the corresponding eigenfunction. Some standard arguments yield that $w(r)$ changes its sign exactly once. Using $\left\langle v,L_uv\right\rangle \geq 0$ for any $v \in H_{rad}^1$ satisfying $v \bot u$ and $\left\langle u,L_uu\right\rangle < 0$, we obtain $\int_{\mathbb{R}^n}u wdx = 0$. We can prove $\int_{\mathbb{R}^n}(f(x,u) - f_t(x,u)u)wdx = 0$ and $\int_{\mathbb{R}^n}(2f(x,u) + x\cdot \nabla_x(f(x,u) - Vu))wdx = 0$. Then we can find a contradiction using $(h_{13})$. See more details in subsection \ref{sec. A2}.

\begin{proposition} \label{prop abs non degeneracy}
	Suppose that conditions $(A 1)$, $(A 3)$ and $(H 3)$ hold. Let $u \in S_c$ and $\lambda$ be the corresponding Lagrange multiplier. Then $\ker L_u|_{L_{rad}^2} = \{0\}$.
\end{proposition}

\subsubsection{Sketch of the Proof of $(B)$}

The key idea is to use a suitable scaling transform and the uniqueness, non-degeneracy of the NGSS of the limit equation. 

\begin{proposition}
	Under assumption $(H 4)$, there exist $c_+ > c_- > 0$ such that for $0 <
	c \leqslant c_-$ and $c \geqslant c_+$, the set $S_c$ contains
	at most one element (up to the natural invariances of the problem).
\end{proposition}

\begin{proof}
	This is a corollary obtained by combining Proposition \ref{abs uniqueness} and Lemma \ref{assumptions hold}. 
\end{proof}

\subsubsection{Sketch of the Proof of $(C)$}

Note that any $u \in S_c$ is radial and positive,  $(C)$ is a consequence of the uniqueness of positive radial solution (c.f. \cite{Adachi-Shibata-Watanabe} when $(h_{16})$-$(h_{18})$ hold, and c.f. \cite{Chen-Lin} when $(h_{19})$ holds).

When $(h_{20})$ holds, any $u \in S_c$ is a GSS (see \cite[Appendix C]{Hajaiej-Song-unique}) and $(C)$ is corollary of \cite[Theorem 5.16]{Hajaiej-Song-gss}.

\section{General results on the variational problem}\label{sono}

This section is devoted to the proof of Proposition \ref{dubai}.

\subsection{Proof of $(A 1)$ and $(A 3)$ } \label{sec. A1 A3}

\subsubsection{Proof of Lemma \ref{lem estimate of m(c)}}

\begin{proof} Assumption $(h_1)$ implies that for any $\epsilon > 0$, there exists $C_\epsilon > 0$ such
	that
	\begin{equation*}
		|\int_{\mathbb{R}^n}F(x,u)dx| \leq \epsilon |u|^2 + C_\epsilon|u|^{2+\frac4n}.
	\end{equation*}
	The classical Gagliardo-Nirenberg inequality shows that
	\begin{equation}
		\int_{\mathbb{R}^n}|u|^{2+4/n}dx \leq C(n)\int_{\mathbb{R}^n}|\nabla u|^2dx\left( \int_{\mathbb{R}^n}|u|^2dx\right)^{\frac{2}{n}}.
	\end{equation}
	Choose $\epsilon > 0$ such that $\epsilon C(n)c^{\frac{2}{n}} = \frac14$. Then, for any $u \in S_c$, we have
	\begin{equation} \label{eq5.2}
		E(u) \geq \frac14\int_{\mathbb{R}^n}|\nabla u|^2dx - Cc.
	\end{equation}
	We get that
	$$
	m(c) = \inf_{u \in S_c}E(u) > -\infty.
	$$
	
	Now we aim to prove that $m(c) < 0$. For any $u \in S_c$, let
	$$
	u_\tau = \tau^{\frac{n}{2}}u(\tau x) \in S_c.
	$$
	
	\underline{If $(h_8)$ hold:}
	When $V \leq 0$, we have
	\begin{align}
		E(u_\tau) \leq & \frac{\tau^2}{2}\int_{\mathbb{R}^n}|\nabla u|^2dx - \tau^{-N}\int_{\mathbb{R}^n}F(\frac{x}{\tau},\tau^{\frac{n}{2}}u)dx \nonumber \\
		= & \tau^2\left(\frac{1}{2}\int_{\mathbb{R}^n}|\nabla u|^2dx - \tau^{-n-2}\int_{\mathbb{R}^n}F(\frac{x}{\tau},\tau^{\frac{n}{2}}u)dx\right).
	\end{align}
	By $(h_7)$ and Fatou's lemma we obtain
	\begin{align}
		\liminf_{\tau \to 0^+}\tau^{-n-2}\int_{\mathbb{R}^n}F(\frac{x}{\tau},\tau^{\frac{n}{2}}u)dx \geq \int_{\mathbb{R}^n}\liminf_{\tau \to 0^+}\tau^{-n-2}F(\frac{x}{\tau},\tau^{\frac{n}{2}}u)dx = +\infty.
	\end{align}
	Hence, $E(u_\tau) < 0$ if $\tau$ is small enough, implying that $m(c) < 0$.

\

	\underline{The case of $(h_9)$ holding:}
	By $(h_9)$, we obtain
	\begin{align}
		& \int_{|x|\leq \tau}V^+(\frac{x}{\tau})|u|^2dx = \int_{|x| \leq \tau}V^+(\frac{x}{\tau})|u|^2dx + \int_{|x| > \tau}V^+(\frac{x}{\tau})|u|^2dx\nonumber \\
		\lesssim & (\int_{|x| \leq \tau}|V^+(\frac{x}{\tau})|^{q}dx)^{\frac1q}(\int_{|x| \leq \tau}|u|^{2q^\ast}dx)^{\frac{1}{q^\ast}} + \tau^{-\alpha}\int_{|x| > \tau}V^+(x)|u|^2dx \nonumber \\
		\lesssim & \tau^{\frac{n}{q}}(\int_{B_1}|V^+|^{q}dx)^{\frac1q} + \tau^{-\alpha}\int_{\mathbb{R}^n}V^+(x)|u|^2dx,	
	\end{align}
	where $1/q + 1/q^\ast = 1$ and thus $2q^\ast \leq 2^\ast$. Therefore,
	\begin{align}
		E(u_\tau) \lesssim \frac{\tau^2}{2}\int_{\mathbb{R}^n}|\nabla u|^2dx & + \frac12\tau^{\frac{n}{q}}(\int_{B_1}|V^+|^{q}dx)^{\frac1q} \nonumber \\
		& +\frac12\tau^{-\alpha}\int_{\mathbb{R}^n}V^+(x)|u|^2dx - \tau^{-N}\int_{\mathbb{R}^n}F(\frac{x}{\tau},\tau^{\frac{n}{2}}u)dx \nonumber \\
		= \tau^\Gamma(\frac{1}{2}\tau^{2-\Gamma}\int_{\mathbb{R}^n}|\nabla u|^2dx & + \frac12\tau^{\frac{n}{q}-\Gamma}(\int_{B_1}|V^+|^{q}dx)^{\frac1q} \nonumber \\
		& +\frac12\tau^{-\alpha-\Gamma}\int_{\mathbb{R}^n}V^+(x)|u|^2dx - \tau^{-n-\Gamma}\int_{\mathbb{R}^n}F(\frac{x}{\tau},\tau^{\frac{n}{2}}u)dx). 	
	\end{align}
	where $\Gamma = \min\{2,\frac{n}{q},-\alpha\}$. By $(h_7)$ and Fatou's lemma we obtain
	\begin{align}
		\liminf_{\tau \to 0^+}\tau^{-n-\Gamma}\int_{\mathbb{R}^n}F(\frac{x}{\tau},\tau^{\frac{n}{2}}u)dx \geq \int_{\mathbb{R}^n}\liminf_{\tau \to 0^+}\tau^{-n-\Gamma}F(\frac{x}{\tau},\tau^{\frac{n}{2}}u)dx = +\infty.
	\end{align}
	Hence,  $E(u_\tau) < 0$ if $\tau$ is small enough, implying that $m(c) < 0$. The proof is complete.
\end{proof}

\subsubsection{Proof of Lemma \ref{lem decreasedness}}

\begin{proof} For any $u \in S_c$ and $k > 1$, using $(h_5)$ we have
	\begin{equation}
		E(ku) = \frac{k^{2}}{2}\int_{\mathbb{R}^n}(|\nabla u|^2+V|u|^2)dx - \int_{\mathbb{R}^n}F(x,ku)dx < k^{2}E(u).
	\end{equation}
	Let $\{u_k\}$ be a minimizing sequence for $m(c)$. Then we know that $m(k^2c) \leq k^2m(c)$. Since $m(c) < 0$, we have $m(k^2c) \leq k^2m(c) < m(c)$ for any $k > 1$, this completes the proof of the decreasiness. We can use standard arguments to prove the continuity of $m(c)$, see the proof of \cite[Lemma 4.2]{Yang-Qi-Zou-2022} or \cite[Theorem 3.3]{Esfahani-Hajaiej-Luo-Song}, and we omit the details. This completes the proof of (i).
	
	Next we prove (ii). WLOG, we may assume that $c_1 \leq c_2$. Then
	\begin{align}
		m(c) & = m(\frac{c}{c_2}c_2) \nonumber \\
		& \leq \frac{c}{c_2}m(c_2) = m(c_2) + \frac{c_1}{c_2}m(c_2) \nonumber \\
		& = m(c_2) + \frac{c_1}{c_2}m(\frac{c_2}{c_1}c_1) \leq m(c_2) + m(c_1).
	\end{align}
	If $m(c_1)$ can be achieved by $u^\ast$, using $(h_5)$ we obtain
	\begin{align}
		m(k^2c_1) \leq E(ku^\ast) & = \frac{k^{2}}{2}\int_{\mathbb{R}^n}(|\nabla u^\ast|^2+V|u^\ast|^2)dx - \int_{\mathbb{R}^n}F(x,ku^\ast)dx \nonumber \\
		& < k^{2}E(u^\ast) = m(c_1),
	\end{align}
	yielding that $m(c) < m(c_2) + m(c_1)$. Similarly, if $m(c_2)$ can be achieved, then $m(c) < m(c_2) + m(c_1)$. The proof is now complete.
\end{proof}

\subsubsection{Proof of Lemma \ref{lem coompactness}}

\begin{proof}  By the proof of Lemma \ref{lem estimate of m(c)}, $E$ is coercive on $S_c$. Hence, if $E(u_k) \rightarrow m(c)$, we know that $\int_{\mathbb{R}^n}|\nabla u_k|^2dx$ is bounded. Therefore, $u_k$ is bounded in $H^1$. Up to a subsequence, we may assume that there exists $u \in H^1$ such that $u_k \rightharpoonup u$ in $H^1$, $u_k \rightharpoonup u$ in $L^2$ and $u_k \rightarrow u$ in $L_{loc}^p, 2 \leq p < 2^\ast$.

\

	Step 1: $u \neq 0$.

\

	Suppose on the contrary that $u \equiv 0$. By $(h_6)$, $V(x) \to 0$ as $|x| \to \infty$ and the strong convergence of $u_k$ in $L_{loc}^p, 2 \leq p < 2^\ast$, one can verify that
	$$
	\int_{\mathbb{R}^n}V|u_k|^2dx \to 0, \quad \int_{\mathbb{R}^n}F(x,u_k)dx \to 0.
	$$
	Hence,
	$$
	E(u_k) = \frac12\int_{\mathbb{R}^n}|\nabla u_k|^2dx + o_n(1),
	$$
	implying that $m(c) \geq 0$. This contradicts $m(c) < 0$.

\

	Step 2: $Q(u) = c$.

\

	By the weak lower semi-continuity of the norm, we can assume that $Q(u) = c_1 \in (0,c]$. Arguing by contradiction, assume $c_1 < c$. By Br\'{e}zis-Lieb lemma, one gets
	$$
	E(u_k) = E(u) + E(u_k-u) + o_n(1) \geq m(c_1) + E(u_k-u) + o_n(1).
	$$
	Similar to the Step 1, we can prove that $E(u_k-u) \geq o_n(1)$. Thus,
	$$
	m(c) + o_n(1) = E(u_k) \geq m(c_1) + o_n(1),
	$$
	implying that $m(c) \geq m(c_1)$, which is a contradiction.

\

	Step 3: Strong convergence.

\

	The weak convergence of $u_k$ in $L^2$ and the fact that $Q(u_k) \rightarrow Q(u)$ show that $u_k$ strongly converges to $u$ in $L^2$. Since $Q(u) = c$, we have
	$$
	m(c) \leq E(u) = E(u_k) - E(u_k-u) + o_n(1) \leq m(c) + o_n(1)
	$$
	showing that $E(u_k-u) \to 0$. This yields to the strong convergence of $u_k$ in $H^1$.
\end{proof}

\subsubsection{Proof of Proposition \ref{A1 and A3}}

\begin{proof}
	Lemmas \ref{lem compact1}, \ref{lem estimate of m(c)} and \ref{lem coompactness} show that if $(H 1)$ holds, then for any $c \in ]0,\infty[$, $- \infty < m (c) < 0$ and the set $S_c$ contains at least one element. By Lagrange multiplier principle, $u \in S_c$ satisfies $\tmop{Eq}_\lambda(u) = 0$ for some $\lambda \in \mathbb{R}$. We now prove $\lambda < \Lambda_0$ under the conditions $(H 2)$. Note that $u$ satisfies the following Pohozaev identity (see \cite[Lemma 5.6]{Hajaiej-Song-gss})
	\begin{align}
		& (n-2)\int_{\mathbb{R}^n}|\nabla u|^2 + n\int_{\mathbb{R}^n}V|u|^2 + \int_{\mathbb{R}^n}x\cdot \nabla_xV|u|^2 \nonumber \\
		= & \lambda n \int_{\mathbb{R}^n}|u|^2 + 2n\int_{\mathbb{R}^n}F(x,u) + 2\int_{\mathbb{R}^n}x\cdot \nabla_xF(x,u).
	\end{align}
	Then we have
	\begin{align}
		& \int_{\mathbb{R}^n}(2V+x\cdot \nabla_xV)|u|^2-2\lambda \int_{\mathbb{R}^n}|u|^2 \nonumber \\
		= & \int_{\mathbb{R}^n}2nF(x,u) + 2x\cdot \nabla_xF(x,u) - (n-2)f(x,u)u.
	\end{align}
	By $(h_{11})$ we obtain $\lambda < 0 = \Lambda_0$. 
Using $u \mapsto |u|$ and that $E(|u|) = E(u)$ we may assume $|u|$ is also a solution. By the strong maximum principle, $|u|$ is positive, thus, $u = |u|$ is positive since $u$ is continuous. Then using the moving plane method or Schwarz rearrangement (c.f. \cite{Naito-2000} or \cite{Burchard-Hajaiej-2006}) we know that $u = u(r)$ is radial, positive and decreasing in $r > 0$. This is standard and we omit the details.
	
	Finally we show that $L_u$ has exactly one simple eigenvalue below $0$. Since $u$ is a NGSS of (\ref{solarstone}) on $X_{c}$ and the codimension of $X_{c}$ in $H^1$ is $1$, we know that $L_u$ has at most one simple eigenvalue below $0$. It is sufficient to show that the first eigenvalue of $L_u|_{L_{rad}^2}$ is indeed below $0$. When $V \equiv 0$ and $f(x,u) = f(u)$, we use spherical harmonic decomposition. Let $\mathcal{H}_0 = L^2_{rad}(\mathbb{R}^n)$ and
	$$
	\mathcal{H}_k := \text{span}\{w(r)Y_k, w(r) \in \mathcal{H}_0, -\Delta_{\mathbb{S}^{n-1}}Y_k = k(k+n-2)Y_k\}, k = 1, 2, \cdots
	$$
	We get an orthogonal decomposition
	$$
	L^2(\mathbb{R}^n) = \oplus_{k = 0}^\infty\mathcal{H}_k.
	$$
	Let $L_k = L_u|_{\mathcal{H}_k}$. Observe that $L_k$ can be viewed as an operator acting on the subspace of radial functions $L^2_{rad}(\mathbb{R}^n)$, through the formula
	$$
	L_k = -\partial_{rr} - \frac{n-1}{r}\partial_r - \lambda - f'(u) + \frac{k(k+n-2)}{r^2}.
	$$
	Direct computation yields to $L_1(u') = 0$ and
	$$
	\langle L_0(u'),u' \rangle = \langle L_1(u'),u' \rangle - (n-1)\int_0^\infty (u')^2r^{n-3}dr < 0.
	$$
	This shows that $L_u|_{L_{rad}^2}$ has at least one simple eigenvalue below $0$.
	
	In the other case when there exists $d \in \mathbb{R}$ such that $f_t(x,t)t - f(x,t) + d(2\lambda t - x\cdot\nabla_xV t + 2f(x,t) + x\cdot\nabla_xf(x,t))$ does not change sign for all $t > 0$, $\lambda < 0$, we argue by contradiction. Let $0$ be the first eigenvalue of $L_u$ and $w$ be the corresponding eigenfunction which is radial and positive. Note that
	$$
	L_u(u) = f(x,tu) - f_t(x,u)u,
	$$
	$$
	L_u(x\cdot u) = 2\lambda u - x\cdot\nabla_xV u + 2f(x,u) + x\cdot\nabla_xf(x,u).
	$$
	We know that
	$$
	\int_{\mathbb{R}^n}(f_t(x,u)u - f(x,u) + d(2\lambda u - x\cdot\nabla_xV u + 2f(x,u) + x\cdot\nabla_xf(x,u)))w = 0.
	$$
	However, this contradicts that both $w$ and $f_t(x,u)u - f(x,u) + d(2\lambda u - x\cdot\nabla_xV u + 2f(x,u) + x\cdot\nabla_xf(x,u))$ do not change sign. The proof of $(A 1)$ and $(A 3)$ is complete.
\end{proof}

\subsection{Proof of $(A 2)$} \label{sec. A2}

We address the problem on non-degeneracy in this subsection.
If the first eigenvalue of an operator of the form $-\Delta + W(|x|)$ is simple and the first eigenfunction is radial and positive (this result can be proved by Perron-Frobenius arguments), then the second eigenfunction changes its sign at least once in order to satisfy the orthogonality condition. Using Courant’s nodal domain theorem, the second eigenfunction has at most $2$ nodal domains, implying that it changes its sign exactly once on $r = |x|$. Hence, we can obtain

\begin{lemma}
	\label{formerA2}
	Let $u \in S_{c}$ be radial and $W = V - f_t(.,u)$. Then the first eigenvalue of $-\Delta + W$ is simple and the corresponding eigenfunction is in $H_{rad}^1$.  Furthermore, the second eigenvalue of $(-\Delta + W)|_{L_{rad}^2}$ (if it exists) has an eigenfunction with exactly one change of sign.
\end{lemma}

\subsubsection{Proof of Proposition \ref{prop abs non degeneracy}}

\begin{proof} We argue by contradiction. Suppose that there exists $w \in H_{rad}^1$ with $\|w\| = 1$ such that $L_uw = 0$.
	
	\
	
	Step 1: $\left\langle v,L_uv\right\rangle \geq 0$ for any $v \in H_{rad}^1$ satisfying $v \bot u$.
	
	\
	
	Since $u$ is a NGSS of (\ref{solarstone}) on $X_{c}$, one gets that $\left\langle v,L_uv\right\rangle \geq 0$ for any $v \in T_uX_{c}$, where $T_uX_{c}$ is the tangent space of $X_{c}$ at $u$ in $H^1$. Noticing that $T_uX_{c} = \{v \in H^1: v \bot u\}$, we complete the proof of Step 1.
	
	\
	
	Step 2: $u \bot w$, i.e. $\int_{\mathbb{R}^n}u wdx = 0$.
	
	\
	
	Noticing that $w - \|u\|_{L^2}^{-2}\left\langle w,u\right\rangle u \bot u$, we have that
	\begin{eqnarray}
		0 &\leq& \left\langle w - \|u\|_{L^2}^{-2}\left\langle w,u\right\rangle u,L_u(w - \|u\|_{L^2}^{-2}\left\langle w,u\right\rangle u)\right\rangle \nonumber \\
		&=& \|u\|_{L^2}^{-4}\left\langle w,u\right\rangle^2\left\langle u,L_uu\right\rangle.
	\end{eqnarray}
	By $(h_{13})$, $\left\langle u,L_uu\right\rangle < 0$, it follows that $\left\langle w,u\right\rangle = 0$.
	
	\
	
	Step 3: $f(x,u) - f_t(x,u)u \bot w$, i.e. $\int_{\mathbb{R}^n}(f(x,u) - f_t(x,u)u)wdx = 0$.
	
	\
	
	Since
	$$
	L_uu = f(x,u) - f_t(x,u)u,
	$$
	we know that
	\begin{equation}
		\left\langle w,f(x,u) - f_t(x,u)u\right\rangle = \left\langle w,L_uu\right\rangle = \left\langle L_uw,u\right\rangle = 0.
	\end{equation}
	
	\
	
	Step 4: $2f(x,u) + x\cdot \nabla_x(f(x,u) - Vu) \bot w$, i.e. $\int_{\mathbb{R}^n}(2f(x,u) + x\cdot \nabla_x(f(x,u) - Vu))wdx = 0$.

\

	Direct computation yields
	\begin{align}
		L_u(x\cdot\nabla u) & = -2\Delta u + x\cdot \nabla_x(f(x,u) - Vu) \nonumber \\
		& = 2\lambda u + 2f(x,u) + x\cdot \nabla_x(f(x,u) - Vu).
	\end{align}
	Then using $u \bot w$ we obtain $2f(x,u) + x\cdot \nabla_x(f(x,u) - Vu) \bot w$.
	
	\
	
	Step 5: We find a contradiction to show that $\ker L_u|_{L_{rad}^2} = \{0\}$.
	
	\
	
	Note that $L_u|_{L_{rad}^2}$ has exactly one simple eigenvalue below $0$. Hence, $0$ is the second eigenvalue of $L_u|_{L_{rad}^2}$ and $w$ is the corresponding eigenfunction. Thanks to Lemma \ref{formerA2}, we may assume that $w = w(r)$, $w(\widehat{r}) = 0$ and $w(r) \geq 0$ on $(0,\widehat{r})$, $w(r) \leq 0$ on $(\widehat{r},+\infty)$ for some $\widehat{r} > 0$. By $(h_{13})$ we consider the case that $kf_t(x,t) + (2-k)\frac{f(x,t)}{t} + \frac{x\cdot\nabla_x(f(x,t)+Vt)}{t}$ is non-increasing in $r = |x|$ and increasing in $t > 0$ for some $k \in \mathbb{R}$. The proof of other cases is similar. Set
	\begin{equation}
		\varphi = \varphi(r) = C_0u - (kf_t(r,u)u + (2-k)f(r,u) + rf_r(r,u)-rV'(r)u),
	\end{equation}
	where
	$$
	C_0 = kf_t(\widehat{r},u(\widehat{r})) + (2-k)\frac{f(\widehat{r},u(\widehat{r}))}{u(\widehat{r})} + \frac{\widehat{r}(f_r(\widehat{r},u(\widehat{r})) )}{u(\widehat{r})}- \widehat{r}V'(\widehat{r}),
	$$
	ensuring that $\varphi(\widehat{r}) = 0$. In the proof of $(A 1)$ we know $u$ is decreasing with respect to $r > 0$. Note that $kf_t(r,u) + (2-k)\frac{f(r,u)}{u} + \frac{rf_r(r,u)}{u} -rV'(r)$ is decreasing with respect to $r > 0$. Hence, $\varphi(r) > 0$ on $(0,\widehat{r})$, $\varphi(r) < 0$ on $(\widehat{r},+\infty)$. Then we derive that
	\begin{equation}
		\left\langle w,\varphi\right\rangle = \int_{\mathbb{R}^n}w(C_0u - (kf_t(r,u)u + (2-k)f(r,u) + rf_r(r,u)-rV'(r)u))dx > 0,
	\end{equation}
	contradicting Steps 2, 3, 4.
	The proof of is complete.
\end{proof}

\begin{remark}
	Provided that some conditions are satisfied, our discussions are applicable to other operators (including nonlocal operators) and radial domains, such as problems on a ball with Dirichlet boundary condition.
\end{remark}

\subsection{Proof of $(B)$} \label{sec. B}

We have the following assumptions:

$(B_1)$ There exists $T: \mathbb{R}_+ \rightarrow B(H^1,H^1)$ such that $T(c)$ is invertible, $Q(T(c)u) = 1/cT(u)$.

$(B_2)$ There exist $\rho: \mathbb{R}_+ \rightarrow \mathbb{R}_+$ and $E_c \in C^{2}(H^1, \mathbb{R})$ such that $E_c(T(c)u) = \rho(c)E(u)$.

By $(B_1)$ and $(B_2)$, $u_c \in X_c$ is a NGSS for $E$ $\Leftrightarrow$ $T(c)u_c \in X_1$ is a NGSS for $E_c$. We aim to show the uniqueness of NGSS for $E$ when $c$ is small or large. It is sufficient to prove the uniqueness for $E_c$. Let $\tau(c) = \inf_{u \in X_1}E_c(u)$. We assume the following:

$(B_3)$ $\lim_{c \rightarrow 0^+}E_c = E_0$, $E_0 \in C^{2}(H^1, \mathbb{R})$ has a unique, radial, non-degenerate NGSS $v_0$ on $X_1$. Furthermore, for any $c_k \rightarrow 0^+$ and $v_k \in X_1$ such that $E_{c_k}(v_k) = \tau(c_k)$, up to a subsequence and translations, $v_k \rightarrow v_0$ in $H^1$.

$(B_4)$ $\lim_{c \rightarrow \infty}E_c = E_\infty$, $E_\infty \in C^{2}(H^1, \mathbb{R})$ has a unique, radial, non-degenerate NGSS $v_\infty$ on $X_1$. Furthermore, for any $c_k \rightarrow \infty$ and $v_k \in X_1$ such that $E_{c_k}(v_k) = \tau(c_k)$, up to a subsequence and translations, $v_k \rightarrow v_\infty$ in $H^1$.

\begin{proposition} \label{abs uniqueness}
	$(i)$ Assume that $(B_1)$-$(B_3)$ hold. If all the NGSS of \eqref{solarstone} are radial, then there exists $c_0 > 0$ such that \eqref{solarstone} has at most one NGSS when $c < c_0$.
	
	$(ii)$ Assume that $(B_1)$, $(B_2)$, $(B_4)$ hold. If all the NGSS of \eqref{solarstone} are radial, then there exists $c_\infty > 0$ such that \eqref{solarstone} has at most one NGSS when $c > c_\infty$.
\end{proposition}

\begin{proof} We give the proof of $(i)$ and similar arguments show $(ii)$. First, we use the implicit function theorem for
	$$
	K: H_{rad}^1 \times \mathbb{R}^+ \times \mathbb{R} \to H_{rad}^{1}, \ K(v,c,\lambda) = u - (-\Delta+V-\lambda)^{-1}f(x,u).
	$$
	By $(B_3)$, let $v_0 \in S_1$ be the unique, radial, non-degenerate NGSS for $E_0$ and $\lambda_0$ be the Lagrange multiplier. Then $K(v_0,0,\lambda_0) = 0$ and $K_v(v_0,0,\lambda_0)$ is invertible from $H_{rad}^1$ to $H_{rad}^1$. Using the implicit function theorem, we know that there exist $\epsilon > 0, \delta > 0, c_0 > 0$ such that $K(v,c,\lambda) = 0$ has a unique solution $v(c,\lambda) \in H_{rad}^1$ with $\|v(c,\lambda)-v_0\|_{H^1} \leq \epsilon$ when $c \in [0,c_0)$ and $\lambda \in (\lambda_0-\delta, \lambda_0-\delta)$. Moreover, $v(0,\lambda_0) = v_0$.
	
	We suppose by contradiction that there exist two NGSS $u_k, w_k$ on $S_{c_k}$ with $c_k \to 0$. Then $u_k, w_k$ are radial. Let $\tilde{u}_k = T(c)u_k, \tilde{w}_k = T(c)w_k$. By $(B_1)$ and $(B_2)$, $\tilde{u}_k, \tilde{w}_k \in S_1$ satisfy  $E_{c_k}(\tilde{u}_k) = E_{c_k}(\tilde{w}_k) = \tau(c_k)$. By $(B_3)$, up to a subsequence, $\tilde{u}_k \to v_0, \tilde{w}_k \to v_0$ strongly in $H^1$. Let $\lambda_k^1, \lambda_k^2$ be the Lagrange multipliers of $\tilde{u}_k, \tilde{w}_k$ respectively. The strong convergence implies that $\lambda_k^1 \to \lambda_0, \lambda_k^2 \to \lambda_0$. For large $k$, this contradicts the uniqueness derived from the implicit function theorem. Therefore, the proof is completed.
\end{proof}






\begin{lemma} \label{assumptions hold}
	(i) Assume that $(h_{10})$, $(h_{14})$ hold and $V \equiv 0$. Then $(B_1)$-$(B_3)$ hold true.
	
	(ii) Assume that $(h_{10})$, $(h_{15})$ hold. Then $(B_1)$, $(B_2)$, $(B_4)$ hold true.
\end{lemma}

\begin{proof} Set
	$$
	T(c)u(x) = c^\alpha u(c^\beta x), \alpha = \frac{2}{(\bar{p}-2)n-4}, \beta = \frac{\bar{p}-2}{(\bar{p}-2)n-4},
	$$
	$$
	\rho(c) = c^{\bar{p}\alpha - n\beta},
	$$
	$$
	E_c(v) = \frac{1}{2}\int_{\mathbb{R}^n}|\nabla v|^2dx - c^{\bar{p}\alpha}\int_{\mathbb{R}^n}F(c^\beta x,c^{-\alpha}v)dx,
	$$
	$$
	E_0(v) = \frac{1}{2}\int_{\mathbb{R}^n}|\nabla v|^2dx - \frac1p\int_{\mathbb{R}^n}|x|^\theta|v|^pdx.
	$$	
	
	$(B_1)$ and $(B_2)$ can be verified directly. Now we show that $(B_3)$ holds true. Equation $-\Delta v = \lambda v + |x|^\theta |v|^{p-2}v, \theta \in (-2,0]$ has a unique, positive solution in $H^1_{rad}(\mathbb{R}^n)$ which is non-degenerate for any fixed $\lambda < 0$ and has no positive solution in $H^1_{rad}(\mathbb{R}^n)$ for $\lambda \geq 0$. Under the assumption $(h_{10})$, note that any NGSS of $E_0$ is positive and $E_0$ is scaling invariant. We know $E_0$ has a unique, radial, non-degenerate NGSS $v_0$ on $S_1$. Moveover, $v_0 \in L^\infty$. Take $c_k \rightarrow 0^+$ and $v_k \in S_1$ such that $E_{c_k}(v_k) = \tau(c_k)$, we aim to show that up to a subsequence and up to translations, $v_k \rightarrow v_0$ in $H^1$. The proof will be divided into three steps.
	
\

	Step 1: $\{v_k\}$ is bounded in $H^1$.

\

	Using \eqref{eq near 0} one gets
	\begin{align}
		|c_k^{\bar{p}\alpha}F(c_k^\beta x,c_k^{-\alpha}v_0))| \lesssim \left\{
		\begin{array}{cc}
			(|v_0|^r+|v_0|^s), & |c_k^{-\frac{2-\bar{p}}{2}\alpha}x| \geq 1, \\
			c_k^{(\bar{p}-p)\alpha}(|v_0|^r+|v_0|^s), & |c_k^{-\frac{2-\bar{p}}{2}\alpha}x| < 1.
		\end{array}
		\right.
	\end{align}
	Note that $\alpha < 0$ and $-n < \theta = \frac{2(\bar{p}-p)}{2-\bar{p}} \leq 0$. For large $n$ we know $c_k^{-\frac{2-\bar{p}}{2}\alpha} > 1$. We take a control function
	\begin{align}
		\psi(x) = \left\{
		\begin{array}{cc}
			K(|v_0|^r+|v_0|^s), & |x| \geq 1, \\
			K|x|^{\theta}(|v_0|^r+|v_0|^s), & |x| < 1,
		\end{array} \right.
	\end{align}
	where $K$ is a positive constant large enough. Then using Lebesgue control convergence theorem we obtain
	\begin{align}
		\lim_{n \to \infty}\int_{\mathbb{R}^n}(\frac1p|x|^\theta|v_0|^p-c_k^{\bar{p}\alpha}F(c_k^\beta x,c_k^{-\alpha}v_0))dx = 0.
	\end{align}
	Then together with
	\begin{align}
		\tau(c_k) \leq E_{c_k}(v_0) = E_0(v_0) + \int_{\mathbb{R}^n}(\frac1p|x|^\theta|v_0|^p-c_k^{\bar{p}\alpha}F(c_k^\beta x,c_k^{-\alpha}v_0))dx,
	\end{align}
	one gets $\limsup_{n \rightarrow \infty}\tau(c_k) \leq \tau(0)$.
	
	Similar to the discussions above, by $(h_{14})$ we have $|c_k^{\bar{p}\alpha}F(c_k^\beta x,c_k^{-\alpha}v_k)| \leq \psi_k(x)$ where
	\begin{align}
		\psi_k(x) = \left\{
		\begin{array}{cc}
			K(|v_k|^r+|v_k|^s), & |x| \geq 1, \\
			K|x|^{\theta}(|v_k|^r+|v_k|^s), & |x| < 1.
		\end{array} \right.
	\end{align}
	Thus using H\"{o}lder inequality and Gagliardo-Nirenberg inequality we have
	\begin{align}
		&\int_{\mathbb{R}^n}|c_k^{\bar{p}\alpha}F(c_k^\beta x,c_k^{-\alpha}v_k)|dx \nonumber \\
		\lesssim & \int_{\mathbb{R}^n}(|v_k|^r+|v_k|^s)dx + (\int_{B_1}|x|^{b^\ast\theta}dx)^{\frac{1}{b^\ast}}(\int_{\mathbb{R}^n}|v_k|^{br}dx)^{\frac1b} \nonumber \\
		& + (\int_{B_1}|x|^{a^\ast\theta}dx)^{\frac{1}{a^\ast}}(\int_{\mathbb{R}^n}|v_k|^{as}dx)^{\frac1a} \nonumber \\
		\lesssim & (\int_{\mathbb{R}^n}|\nabla v_k|^2dx)^{\frac{n(r-2)}{4}} + (\int_{\mathbb{R}^n}|\nabla v_k|^2dx)^{\frac{n(s-2)}{4}} \nonumber \\
		& + (\int_{\mathbb{R}^n}|\nabla v_k|^2dx)^{\frac{n(br-2)}{4b}} + (\int_{\mathbb{R}^n}|\nabla v_k|^2dx)^{\frac{n(as-2)}{4a}},
	\end{align}
	where $$\frac{1}{a} + \frac{1}{a^\ast} = 1, \quad \frac{1}{b} + \frac{1}{b^\ast} = 1,$$ $$\frac{n(br-2)}{4b} < 1, \quad \frac{n(as-2)}{4a} < 1,$$ $$a^\ast\theta < -n, \quad b^\ast\theta < -n.$$
	We obtain
	\begin{align}
		E_{c_k}(v_k) = & \frac{1}{2}\int_{\mathbb{R}^n}|\nabla v_k|^2dx - c_k^{\bar{p}\alpha}\int_{\mathbb{R}^n}F(c_k^\beta x,c_k^{-\alpha}v_k)dx \nonumber \\
		\gtrsim & \int_{\mathbb{R}^n}|\nabla v_k|^2dx - (\int_{\mathbb{R}^n}|\nabla v_k|^2dx)^{\frac{n(r-2)}{4}} - (\int_{\mathbb{R}^n}|\nabla v_k|^2dx)^{\frac{n(s-2)}{4}} \nonumber \\
		& - (\int_{\mathbb{R}^n}|\nabla v_k|^2dx)^{\frac{n(br-2)}{4b}} - (\int_{\mathbb{R}^n}|\nabla v_k|^2dx)^{\frac{n(as-2)}{4a}}.
	\end{align}
	Note that
	$$
	\frac{n(r-2)}{4} < 1, \quad \frac{n(s-2)}{4} < 1, \quad \frac{n(br-2)}{4b} < 1, \quad \frac{n(as-2)}{4a} < 1.
	$$
	Using $\limsup_{n \rightarrow \infty}\tau(c_k) \leq \tau(0)$ we have $\limsup_{n \rightarrow \infty}\int_{\mathbb{R}^n}|\nabla v_k|^2dx < \infty$, implying that $\{v_k\}$ is bounded in $H^1$.

\

	Step 2: $\tau(c_k) \rightarrow \tau(0)$.

\

	By $(h_{10})$, using the symmetric radial decreasing rearrangement of $v_k$, up to translations, we can suppose that $v_k$ are nonnegative and radially symmetric. Up to a subsequence, we may assume that $v_k$ converges strongly in $L^d, 2 < d < 2^\ast$.
	
	In Step 1 we have proved that $\limsup_{n \rightarrow \infty}\tau(c_k) \leq \tau(0)$. Therefor, it is sufficient to show that $\liminf_{n \rightarrow \infty}\tau(c_k) \geq \tau(0)$. Note that
	$$
	\tau(0) \leq E_{0}(v_k) = E_{c_k}(v_k) - \int_{\mathbb{R}^n}(\frac1p|x|^\theta|v_k|^p-c_k^{\bar{p}\alpha}F(c_k^\beta x,c_k^{-\alpha}v_k))dx.
	$$
	Similar to Step 1, by $(h_{14})$, the strong convergence of $v_k$ in $L^d, 2 < d < 2^\ast$ and Lebesgue control convergence theorem, one can prove that $$\int_{\mathbb{R}^n}(\frac1p|x|^\theta|v_k|^p-c_k^{\bar{p}\alpha}F(c_k^\beta x,c_k^{-\alpha}v_k))dx \rightarrow 0.$$ Hence, $\liminf_{n \rightarrow \infty}\tau(c_k) \geq \tau(0)$.

\

	Step 3: Up to a subsequence and up to translations in autonomous cases, $v_k \rightarrow v_0$ in $H^1$.

\

	Step 2 shows that
	$$
	E_0(v_k) = E_{c_k}(v_k) + o_n(1) \to \tau(0).
	$$
	Moreover, $\int_{\mathbb{R}^n}|v_k|^2dx \equiv 1$. Then by a standard process, up to a subsequence and up to translations in the autonomous cases, one obtains that $v_k \to v_0$ in $H^1$.

	Now we prove $(ii)$. Let
	$$
	T(c)u(x) = c^\alpha u(c^\beta x), \alpha = \frac{2}{(q-2)n-4}, \beta = \frac{q-2}{(q-2)n-4},
	$$
	$$
	\rho(c) = c^{q\alpha - n\beta},
	$$
	$$
	E_c(v) = \frac{1}{2}\int_{\mathbb{R}^n}|\nabla v|^2dx + \frac{1}{2}c^{2\beta}\int_{\mathbb{R}^n}V(c^\beta x)|v|^2dx  - c^{q\alpha}\int_{\mathbb{R}^n}F(c^\beta x,c^{-\alpha}v)dx,
	$$
	$$
	E_\infty(v) = \frac{1}{2}\int_{\mathbb{R}^n}|\nabla v|^2dx - \frac1q\int_{\mathbb{R}^n}|v|^qdx.
	$$	
	$(B_1)$ and $(B_2)$ can be verified directly. Now we show that $(B_3)$ holds true. It is well-known that quation $-\Delta v = \lambda v + |v|^{q-2}v$ has a unique, radial, positive solution in $H^1(\mathbb{R}^n)$ which is non-degenerate for any fixed $\lambda < 0$ and has no positive solution in $H^1(\mathbb{R}^n)$ for $\lambda \geq 0$. Note that any NGSS of $E_\infty$ is positive and $E_\infty$ is scaling invarient. We know that $E_\infty$ has a unique, radial, non-degenerate NGSS $v_\infty$ on $S_1$. Moveover, $v_\infty \in L^\infty$. Take $c_k \rightarrow \infty$ and $v_k \in S_1$ such that $E_{c_k}(v_k) = \tau(c_k)$, we aim to show that up to a subsequence and up to translations, $v_k \rightarrow v_\infty$ in $H^1$. Note that $\beta < 0$ and $V \in L^\infty$, thus,
	$$
	|c_k^{2\beta}\int_{\mathbb{R}^n}V(c^\beta x)|v_k|^2dx| \leq o_n(1)\int_{\mathbb{R}^n}|v_k|^2dx.
	$$
	Then similar to the proof of $(i)$, we can complete the proof.
\end{proof}

\appendix

\section{Proof of strict decreasiness of $m(c)$ without $(h_5)$} \label{app proof}

\begin{lemma} \label{lem proof}
Let $V\equiv 0, f(u) = |u|^{p-2}u - |u|^{q-2}u$. If $2 < p < q$ and $p < 2 + \frac{4}{n}$, then $m(c) \in ]-\infty,0[$ and $m(c)$ is strictly decreasing.
\end{lemma}

\begin{proof}
	That $m(c) > -\infty$ can be deduced by Gagliardo-Nirenberg inequality, see \cite[Theorem 1]{Song-Hajaiej-partial confinement}. Using $u^t(x) = t^{\frac{n}{2}}u(tx)$ and $$E(u^t) = \frac{t^2}{2}\int_{\mathbb{R}^n}|\nabla u|^2 - \frac{t^{\frac{p-2}{2}n}}{p}\int_{\mathbb{R}^n}|u|^p + \frac{t^{\frac{q-2}{2}n}}{q}\int_{\mathbb{R}^n}|u|^q,$$ we know that $E(u^t) < 0$ for small and thus $m(c) < 0$ for all $c$. Then we know that $m(c) \in ]-\infty,0[$. Let $c = c_1 + c_2$. Take the minimizing sequences $(u_k)_{k \in \mathbb{N}}$, $(v_k)_{k \in \mathbb{N}}$ for $m(c_1)$ and $m(c_2)$. We can also assume that $u_k$ and $v_k$ have compact support. Set $w_k = u_k(\cdot) + v_k(\cdot + z_k)$ where $z_k$ large such that $\text{supp}u_k \cap \text{supp}v_k = \emptyset$. Then $w_k \in X_c$ and thus $E(w_k) \geq m(c)$. We obtain $m(c) \leq E(u_k) + E(v_k) \leq m(c_1) + m(c_2) + o_k(1)$, implying $m(c) \leq m(c_1) + m(c_2)$. Note that $m(c_2) < 0$. We obtain $m(c) < m(c_1)$ for any $c_1 < c$. This completes the proof.
\end{proof}

\end{document}